\documentclass[a4paper]{article}

\usepackage{amssymb}

\def\GAP{\textsf{GAP}}
\def\ATLAS{\textsc{Atlas}}

\def\GL{{\rm GL}}

\begin{document}

\title{Constructing the ordinary character tables of some Atlas groups using character theoretic methods}

\author{\textsc{Thomas Breuer} \\[0.5cm]
\textit{Lehrstuhl D f{\"u}r Mathematik} \\
\textit{RWTH, 52056 Aachen, Germany}}

\date{May 8th, 2016}

\maketitle

\abstract{%
We compute the character tables of the following groups
with character theoretic methods,
using known information about the conjugacy classes
and about the character tables of some subgroups and factor groups:
$Th$, $J_4$, $2.{}^2E_6(2)$, $2.{}^2E_6(2).2$, and $2.B$.}

\textwidth16cm
\oddsidemargin0pt

\parskip 1ex plus 0.5ex minus 0.5ex
\parindent0pt

\tableofcontents

\section{Overview}

The character tables of
$Th$, $J_4$, $2.{}^2E_6(2)$, $2.{}^2E_6(2).2$, and $2.B$
are shown in the {\ATLAS} of Finite Groups~\cite{CCN85}, and these {\ATLAS}
tables are contained in the character table library~\cite{CTblLib}
of the computer algebra system {\GAP}~\cite{GAP483}.
Except for the case of $Th$ (see~\cite{BMO16}),
I am not aware of published proofs of the correctness
of these character tables.

In the following sections,
we show how one can compute the tables in question with {\GAP},
using character theoretic methods.
We will assume the character tables of certain proper subgroups
and factor groups;
except in the case of $2.B$, these tables have been verified either
by direct computations with the group in question
or by character theoretic methods.
I am not aware of such a verification for the character table of $B$,
thus the construction of the character table of $2.B$ in
Section~\ref{section2B} relies on the correctness of the {\ATLAS} table
of $B$.

The main tools for the computation of the irreducible characters
will be LLL reduction (see~\cite{LLL82})
and the enumeration of orthogonal embeddings (see~\cite{Ple90}).
Several {\GAP} library functions from the Chapter
``Maps Concerning Character Tables'' of the {\GAP} Reference Manual
will be used without comments.
We will use the {\GAP} Character Table Library,
thus we have to load this {\GAP} package.

\begin{verbatim}
    gap> LoadPackage( "ctbllib" );
    true
\end{verbatim}

\section{The character table of $Th$ (December 30th, 2015)}

The character table of the sporadic simple Thompson group $Th$
has been published in~\cite[pp.~162--163]{Smi76c} without proof.
This table (with some rows and columns permuted) is shown in the
{\ATLAS} of Finite Groups (see~\cite[p.~176]{CCN85}),
and this {\ATLAS} table is contained in the character table library
~\cite{CTblLib}
of the computer algebra system {\GAP}~\cite{GAP483}.

We start with the description of the conjugacy classes of $Th$
as given in~\cite{Par77} (see Section\ref{sectclassesTh}),
then we compute the character table of a subgroup of type $2^5:L_5(2)$
in $Th$ (the so-called Dempwolff group, see Section~\ref{sectsubgroupTh}),
then we write down the power maps of $Th$
(see Section~\ref{sectpowermapsTh}),
and finally we compute the irreducible characters of $Th$
(see Section~\ref{sectirreduciblesTh}).

\subsection{The conjugacy classes of $Th$}\label{sectclassesTh}

The conjugacy classes of elements of order different from
$1$, $19$, and $31$ in $Th$ are listed in~\cite[Table~I]{Par77},
and~\cite[(6.3)]{Par77} states that there are
one class of elements of order $19$
and two (nonreal) classes of elements of order $31$;
note that Case~I~($\alpha$) holds by the proof of (6.1).

Up to a permutation of classes,
this description of the classes agrees with the element orders
and centralizer orders in the character table that is claimed for $Th$
in {\GAP}'s Character Table Library~\cite{CTblLib}.

\begin{verbatim}
    gap> lib:= CharacterTable( "Th" );;
    gap> parrottnames:= [
    >      "1A", "z", "c2", "c3", "c1", "r1", "v", "b", "zc1", "zc2",
    >      "zc3", "a", "us1", "w", "f1", "f3", "f2", "zb",
    >      "r1c2", "(r1c2)^-1", "r1c3", "vc1", "l", "za",
    >      "c1b", "(c1b)^-1", "zf1", "zf2", "19A", "vb", "c2a2",
    >      "us1c2", "(us1c2)^-1", "wc1", "(wc1)^-1",
    >      "f4", "f5", "(f5)^-1", "r1a", "zbc1", "(zbc1)^-1",
    >      "31A", "31B", "r1f1", "s1f1", "(s1f1)^-1",
    >      "c2l", "(c2l)^-1" ];;
    gap> orders:= OrdersClassRepresentatives( lib );;
    gap> centralizers:= SizesCentralizers( lib );;
    gap> descr:= TransposedMat( [ parrottnames, orders, centralizers ] );;
    gap> for entry in descr do
    >      Print( String( entry[1], -12 ),
    >             String( entry[2], 2 ), "  ",
    >             StringPP( entry[3] ), "\n" );
    >    od;
    1A           1  2^15*3^10*5^3*7^2*13*19*31
    z            2  2^15*3^4*5*7
    c2           3  2^6*3^7*7*13
    c3           3  2^3*3^10
    c1           3  2^4*3^7*5
    r1           4  2^11*3^3*7
    v            4  2^9*3*5
    b            5  2^3*3*5^3
    zc1          6  2^4*3^3*5
    zc2          6  2^6*3^3
    zc3          6  2^3*3^4
    a            7  2^3*3*7^2
    us1          8  2^7*3
    w            8  2^5*3
    f1           9  2^3*3^6
    f3           9  3^6
    f2           9  2*3^4
    zb          10  2^3*3*5
    r1c2        12  2^5*3^2
    (r1c2)^-1   12  2^5*3^2
    r1c3        12  2^2*3^3
    vc1         12  2^3*3
    l           13  3*13
    za          14  2^3*7
    c1b         15  2*3*5
    (c1b)^-1    15  2*3*5
    zf1         18  2^3*3^2
    zf2         18  2*3^2
    19A         19  19
    vb          20  2^2*5
    c2a2        21  3*7
    us1c2       24  2^3*3
    (us1c2)^-1  24  2^3*3
    wc1         24  2^3*3
    (wc1)^-1    24  2^3*3
    f4          27  3^3
    f5          27  3^3
    (f5)^-1     27  3^3
    r1a         28  2^2*7
    zbc1        30  2*3*5
    (zbc1)^-1   30  2*3*5
    31A         31  31
    31B         31  31
    r1f1        36  2^2*3^2
    s1f1        36  2^2*3^2
    (s1f1)^-1   36  2^2*3^2
    c2l         39  3*13
    (c2l)^-1    39  3*13
\end{verbatim}

We create a new character table object for a group of the given order
and with the given element orders and centralizer orders.

\begin{verbatim}
    gap> th:= rec( UnderlyingCharacteristic:= 0,
    >              OrdersClassRepresentatives:= orders,
    >              SizesCentralizers:= centralizers,
    >              Size:= centralizers[1] );;
    gap> ConvertToCharacterTableNC( th );;
\end{verbatim}

\subsection{The subgroup $2^5.L_5(2)$}\label{sectsubgroupTh}

The group $Th$ contains a subgroup $D$ that is a non-split extension
of an elementary abelian group of order $2^5$ by the general linear group
$\GL(5,2)$, see~\cite{Smi76c}.
This subgroup is uniquely determined up to isomorphism by these properties
(see~\cite{Dem72}), it is usually called the \emph{Dempwolff group}.

The character table of $D$ is available in {\GAP}'s character table library,
but we recompute it anew from a permutation representation of $D$,
in order to make the construction of the character table of $Th$
self-contained.

Since {\GAP}'s current default algorithm requires more than 4~GB of space
--too much for my small notebook--
we use some character theoretic methods for computing the irreducible
characters of $D$.

First we compute the irreducible characters that are inflated from the
factor group $\GL(5,2)$.


\begin{verbatim}
    gap> g:= AtlasGroup( "2^5.L5(2)" );;
    gap> bl:= Blocks( g, MovedPoints( g ) );;
    gap> Length( bl[1] );
    2
    gap> acthom:= ActionHomomorphism( g, bl, OnSets );;
    gap> img:= Image( acthom );;
    gap> Size( g ) / Size( img );
    32
    gap> sm:= SmallerDegreePermutationRepresentation( img );;
    gap> NrMovedPoints( Image( sm ) );
    31
    gap> f:= CharacterTable( Image( sm ) );;
    gap> d:= CharacterTable( g );;
    gap> fus:= List( ConjugacyClasses( d ),
    >            c -> PositionProperty( ConjugacyClasses( f ),
    >                   cc -> ( Representative( c )^acthom )^sm in cc ) );;
    gap> infl:= List( Irr( f ), x -> x{ fus } );;
\end{verbatim}

Next we compute the characters induced from all linear characters of
cyclic subgroups of $D$,
and the permutation character of the given permutation representation.
Reducing these characters with the known irreducibles yields one
faithful irreducible character.
Then we form tensor products of the known irreducible characters
with the faithful irreducible character, reduce them with the known
irreducibles,
and apply the LLL algorithm to the reducible characters which we have;
this yields five new irreducible characters.

\begin{verbatim}
    gap> indcyc:= InducedCyclic( d, [ 2 .. NrConjugacyClasses( d ) ], "all" );;
    gap> nat:= NaturalCharacter( g );;
    gap> red:= ReducedOrdinary( d, infl, Concatenation( indcyc, [ nat ] ) );;
    gap> Length( red.irreducibles );
    1
    gap> faithirr:= ShallowCopy( red.irreducibles );;
    gap> ten:= Set( Tensored( infl, faithirr ) );;
    gap> ten:= Reduced( d, faithirr, ten );;
    gap> lll:= LLL( d, Concatenation( red.remainders, ten.remainders ) );;
    gap> Length( lll.irreducibles );
    5
    gap> Append( faithirr, lll.irreducibles );
\end{verbatim}

Next we compute symmetrization of the known irreducible characters,
reduce them, and apply LLL again --four new irreducibles.

\begin{verbatim}
    gap> sym2:= Symmetrizations( d, faithirr, 2 );;
    gap> sym3:= Symmetrizations( d, faithirr, 3 );;
    gap> irr:= Concatenation( infl, faithirr );;
    gap> sym:= Reduced( d, irr, Concatenation( sym2, sym3 ) );;
    gap> lll:= LLL( d, Concatenation( lll.remainders, sym.remainders ) );;
    gap> Length( lll.irreducibles );
    4
    gap> Append( irr, lll.irreducibles );
\end{verbatim}

Next we compute the possible orthogonal embeddings of the
four-dimensional LLL-reduced lattice into the four-dimensional
standard lattice.
we get two solutions for the missing irreducibles of $D$.

\begin{verbatim}
    gap> gram:= MatScalarProducts( d, lll.remainders, lll.remainders );;
    gap> emb:= OrthogonalEmbeddings( gram );;                           
    gap> Length( emb.solutions );
    3
    gap> dec:= List( emb.solutions,
    >                x -> Decreased( d, lll.remainders, emb.vectors{ x } ) );;
    gap> dec:= Filtered( dec, x -> x <> fail );;
    gap> Length( dec );
    2
\end{verbatim}

The second solution is not compatible with the given $2$-nd power map of $D$.
Thus the first solution is the correct one.

\begin{verbatim}
    gap> sym:= Symmetrizations( d, [ dec[2].irreducibles[1] ], 2 );;
    gap> ForAll( dec[2].irreducibles,
    >            x -> IsInt( ScalarProduct( d, sym[1], x ) ) ); 
    false
    gap> SetIrr( d, Concatenation( irr, dec[1].irreducibles ) );
\end{verbatim}

Finally, we show that we have really computed the character table of
a group which contains an elementary abelian normal subgroup $N$
of order $2^5$, \ldots

\begin{verbatim}
    gap> nsg:= ClassPositionsOfNormalSubgroups( d );
    [ [ 1 ], [ 1, 2 ], [ 1 .. 41 ] ]
    gap> SizesConjugacyClasses( d ){ nsg[2] };
    [ 1, 31 ]
    gap> OrdersClassRepresentatives( d ){ nsg[2] };
    [ 1, 2 ]
\end{verbatim}

\ldots that the extension is non-split, \ldots

\begin{verbatim}
    gap> f:= d / nsg[2];;
    gap> PossibleClassFusions( f, d );
    [  ]
\end{verbatim}

\ldots that the factor group by $N$ is isomorphic with $\GL(5,2)$, \ldots

\begin{verbatim}
    gap> n:= PCore( g, 2 );
    <permutation group with 5 generators>
    gap> Size( n );
    32
    gap> IsomorphismGroups( g / n, GL(5,2) ) = fail;
    false
\end{verbatim}

\ldots and that the computed character table is equivalent to the table
which can fetched from {\GAP}'s character table library via the call
\verb|CharacterTable( "2^5.L5(2)" )|.

\begin{verbatim}
    gap> libsub:= CharacterTable( "2^5.L5(2)" );;
    gap> IsRecord( TransformingPermutationsCharacterTables( d, libsub ) );
    true
    gap> d:= libsub;;
\end{verbatim}

\subsection{The power maps of $Th$}\label{sectpowermapsTh}

Before we can compute the irreducible characters of $Th$ by inducing
characters from $D$ and from cyclic subgroups of $Th$,
we determine the power maps of $Th$.

In addition to the conditions that are imposed by the representative orders
and centralizer orders,
we use the following information from~\cite{Par77};
each entry \verb|[ c, p, i ]| means that the \verb|p|-th power of the class with
the name \verb|c| is the class with the name \verb|i|.

\begin{verbatim}
    gap> powinfo:= [
    >     [ "zc3", 2, "c3" ],     # c3 commutes with z
    >     [ "zf2", 2, "f2" ],     # f2 commutes with z
    >     [ "r1c3", 2, "zc3" ],   # r1 commutes with c3
    >     [ "vc1", 2, "zc1" ],    # v squares to z and commutes with c1
    >     [ "wc1", 2, "vc1" ],    # w squares to v and commutes with c1
    >     [ "(wc1)^-1", 2, "vc1" ],
    >     [ "us1c2", 2, "r1c2" ], # us1 squares to r1
    >     [ "(us1c2)^-1", 2, "(r1c2)^-1" ],
    >     [ "us1", 2, "r1" ],     # (5.1)
    >     [ "w", 2, "v" ],        # (5.1)
    >     [ "zbc1", 2, "c1b" ],
    >     [ "(zbc1)^-1", 2, "(c1b)^-1" ],
    >     [ "f1", 3, "c3" ],
    >     [ "f2", 3, "c3" ],
    >     [ "f3", 3, "c3" ],
    >     [ "vc1", 3, "v" ],      # v commutes with c1
    >     [ "zf1", 3, "zc3" ],
    >     [ "zf2", 3, "zc3" ],
    >     [ "us1c2", 3, "us1" ],
    >     [ "(us1c2)^-1", 3, "us1" ],
    >     [ "wc1", 3, "w" ],
    >     [ "(wc1)^-1", 3, "w" ],
    >     [ "f4", 3, "f3" ],
    >     [ "f5", 3, "f3" ],
    >     [ "(f5)^-1", 3, "f3" ],
    >     [ "r1f1", 3, "r1c3" ],
    >     [ "s1f1", 3, "r1c3" ],
    >     [ "(s1f1)^-1", 3, "r1c3" ],
    >     ];;
\end{verbatim}

Next we enter information about Galois conjugation;
we will need $p$-th power maps for primes $p$ up to the largest
element order in $Th$.

\begin{verbatim}
    gap> maxorder:= Maximum( OrdersClassRepresentatives( th ) );
    39
    gap> primes:= Filtered( [ 1 .. maxorder ], IsPrimeInt );;
\end{verbatim}

The classes of \verb|r1f1| and \verb|f4| (element orders $36$ and $27$, respectively)
are rational.

\begin{verbatim}
    gap> for p in primes do
    >      if 36 mod p <> 0 then
    >        Add( powinfo, [ "r1f1", p, "r1f1" ] );
    >      fi;
    >      if 27 mod p <> 0 then
    >        Add( powinfo, [ "f4", p, "f4" ] );
    >      fi;
    >    od;
\end{verbatim}

For the non-rational classes, it is more suitable to set the power map entries
directly, instead of using the list \verb|powinfo|.
Thus we first initialize the power maps and then evaluate this list.

\begin{verbatim}
    gap> powermaps:= [];;
    gap> for p in primes do
    >      powermaps[p]:= InitPowerMap( th, p );
    >    od;
    gap> for entry in powinfo do
    >      p:= entry[2];
    >      pow:= powermaps[p];
    >      src:= Position( parrottnames, entry[1] );
    >      trg:= Position( parrottnames, entry[3] );
    >      if IsInt( pow[ src ] ) then
    >        if pow[ src ] <> trg then
    >          Error( "contradiction!" );
    >        fi;
    >      elif not trg in pow[ src ] then
    >        Error( "contradiction!" );
    >      else
    >        pow[ src ]:= trg;
    >      fi;
    >    od;
    gap> SetComputedPowerMaps( th, powermaps );
\end{verbatim}

Any non-rational class $C$, say, of $Th$ has exactly one Galois conjugate
class,
that is, the character values attained on $C$ generate a quadratic
extension field of the rationals.
Let $x$ be a generating element of this extension.
If $p$ does not divide the order of the elements in $C$
then the $p$-th power map swaps $C$ and its Galois conjugate
if and only if $x$ differs from the Galois image $x^{\ast p}$,
otherwise the $p$-th power map fixes the two classes.
The following small function is suitable for setting the
power map info in these cases.

\begin{verbatim}
    gap> setGaloisInfo:= function( powermaps, classes, orders, primes, x )
    >    local ord, p;
    >    ord:= orders[ classes[1] ];
    >    for p in primes do
    >      if ord mod p <> 0 then
    >        if GaloisCyc( x, p ) = x then
    >          powermaps[p]{ classes }:= classes;
    >        else
    >          powermaps[p]{ classes }:= classes{ [ 2, 1 ] };
    >        fi;
    >      fi;
    >    od;
    >    end;;
\end{verbatim}

The elements \verb|c1b| and \verb|(c1b)^-1| have order $15$,
so the character values on these elements lie in a
non-real quadratic subfield of the field of $15$-th roots of unity,
that is, in the extension by $\sqrt{-3}$ or $\sqrt{-15}$.
Since the subgroup $D$ contains exactly two Galois conjugate classes
of element order $15$,
we can decide from the character table of $D$ which case occurs.

The same argument holds for the classes of the elements \verb|zbc1| and \verb|(zbc1)^-1|,
which have order $30$.

\begin{verbatim}
    gap> pos:= Positions( OrdersClassRepresentatives( d ), 15 );
    [ 22, 24 ]
    gap> f:= Field( List( Irr( d ), x -> x[ pos[1] ] ) );
    NF(15,[ 1, 2, 4, 8 ])
    gap> Sqrt( -15 ) in f;
    true
    gap> pos:= Positions( orders, 15 );
    [ 25, 26 ]
    gap> setGaloisInfo( powermaps, pos, orders, primes, Sqrt( -15 ) );
    gap> pos:= Positions( OrdersClassRepresentatives( d ), 30 );
    [ 23, 25 ]
    gap> f:= Field( List( Irr( d ), x -> x[ pos[1] ] ) );
    NF(15,[ 1, 2, 4, 8 ])
    gap> pos:= Positions( orders, 30 );
    [ 40, 41 ]
    gap> setGaloisInfo( powermaps, pos, orders, primes, Sqrt( -15 ) );
\end{verbatim}

The character values on \verb|f5|, \verb|(f5)^-1| lie in the unique
quadratic subfield of the field of $27$-th roots of unity,
which is generated by $\sqrt{-3}$.
Analogously, the character values on the two classes of element order $31$
lie in the extension by $\sqrt{-31}$.

\begin{verbatim}
    gap> setGaloisInfo( powermaps,
    >        List( [ "f5", "(f5)^-1" ], x -> Position( parrottnames, x ) ),
    >        orders, primes, Sqrt( -3 ) );
    gap> setGaloisInfo( powermaps, Positions( orders, 31 ), orders, primes,
    >                   Sqrt( -31 ) );
\end{verbatim}

Concerning the two classes of element order $39$,
with representatives \verb|c2l| and \verb|(c2l)^-1|,
the field of character values is generated either by $\sqrt{-3}$
or $\sqrt{-39}$.
We try both possibilities,
only with the second one the characters induced from the cyclic subgroup
of order $39$ have integral norms.

\begin{verbatim}
    gap> pos:= Positions( orders, 39 );
    [ 47, 48 ]
    gap> setGaloisInfo( powermaps, pos, orders, primes, Sqrt( -3 ) );
    gap> indcyc:= InducedCyclic( th, [ pos[1] ], "all" );;
    gap> ForAll( indcyc, x -> IsInt( ScalarProduct( th, x, x ) ) );
    false
    gap> setGaloisInfo( powermaps, pos, orders, primes, Sqrt( -39 ) );
    gap> indcyc:= InducedCyclic( th, [ pos[1] ], "all" );;
    gap> ForAll( indcyc, x -> IsInt( ScalarProduct( th, x, x ) ) );
    true
\end{verbatim}

The elements \verb|s1f1|, \verb|(s1f1)^-1| have order $36$,
the field of character values on their classes is generated
by $\sqrt{-3}$ or $\sqrt{-1}$.
Only the first candidate is compatible with induced characters.

\begin{verbatim}
    gap> pos:= Positions( orders, 36 );
    [ 44, 45, 46 ]
    gap> parrottnames{ pos };
    [ "r1f1", "s1f1", "(s1f1)^-1" ]
    gap> setGaloisInfo( powermaps, [ 45, 46 ], orders, primes, Sqrt( -3 ) );
    gap> indcyc:= InducedCyclic( th, [ 45 ], "all" );;    
    gap> ForAll( indcyc, x -> IsInt( ScalarProduct( th, x, x ) ) );
    true
    gap> setGaloisInfo( powermaps, [ 45, 46 ], orders, primes, Sqrt( -1 ) );
    gap> indcyc:= InducedCyclic( th, [ 45 ], "all" );;
    gap> ForAll( indcyc, x -> IsInt( ScalarProduct( th, x, x ) ) );
    false
    gap> setGaloisInfo( powermaps, [ 45, 46 ], orders, primes, Sqrt( -3 ) );
\end{verbatim}

The elements \verb|wc1|, \verb|(wc1)^-1| have order $24$,
the field of character values on their classes is a non-real
quadratic subfield of the field of $24$-th roots of unity,
the generators to check are $\sqrt{-3}$, $\sqrt{-1}$, $\sqrt{-2}$,
and $\sqrt{-6}$.
Only the last candidate is compatible with induced characters.

\begin{verbatim}
    gap> List( [ "wc1", "(wc1)^-1" ], x -> Position( parrottnames, x ) );
    [ 34, 35 ]
    gap> vals:= [ Sqrt( -3 ), Sqrt( -1 ), Sqrt( -2 ), Sqrt( -6 ) ];
    [ E(3)-E(3)^2, E(4), E(8)+E(8)^3, E(24)+E(24)^11-E(24)^17-E(24)^19 ]
    gap> good:= [];;
    gap> for val in vals do
    >      setGaloisInfo( powermaps, [ 34, 35 ], orders, primes, val );
    >      indcyc:= InducedCyclic( th, [ 34 ], "all" );
    >      if ForAll( indcyc, x -> IsInt( ScalarProduct( th, x, x ) ) ) then
    >        Add( good, val );
    >      fi;
    >    od;
    gap> good;
    [ E(24)+E(24)^11-E(24)^17-E(24)^19 ]
    gap> setGaloisInfo( powermaps, [ 34, 35 ], orders, primes, good[1] );
\end{verbatim}

The only classes for which the information about powers is missing
are those of the elements \verb|r1c2| and \verb|(r1c2)^-1| (of order $12$)
and their roots \verb|us1c2|, \verb|(us1c2)^-1| (of order $24$).
We use that the subgroup $D$ of $Th$ contains elements from these classes,
which can be seen from the fact that the elements of order $12$ in $D$
lie in the classes of \verb|r1c2| and \verb|(r1c2)^-1|.
(The class fusion from $D$ to $Th$ will be determined below,
here we need only the compatibility of element orders and centralizer orders
of a class of the subgroup and the possible image classes in the overgroup.)

Note that the classes of \verb|us1c2|, \verb|(us1c2)^-1| must have the same
field of character values as their squares.

\begin{verbatim}
    gap> parrottnames{ [ 19, 20, 32, 33 ] };
    [ "r1c2", "(r1c2)^-1", "us1c2", "(us1c2)^-1" ]
    gap> fus:= InitFusion( d, th );;
    gap> pos:= Positions( OrdersClassRepresentatives( d ), 12 );
    [ 12, 15, 16 ]
    gap> fus{ pos };
    [ [ 19, 20, 21, 22 ], [ 19, 20 ], [ 19, 20 ] ]
    gap> List( pos, x -> Field( List( Irr( d ), chi -> chi[x] ) ) );
    [ Rationals, CF(3), CF(3) ]
    gap> Sqrt( -3 ) in CF(3);
    true
    gap> setGaloisInfo( powermaps, [ 19, 20 ], orders, primes, Sqrt( -3 ) );
    gap> setGaloisInfo( powermaps, [ 32, 33 ], orders, primes, Sqrt( -3 ) );
\end{verbatim}

\subsection{The irreducible characters of $Th$}\label{sectirreduciblesTh}

We start with creating some characters of $Th$ by inducing all linear
characters of cyclic subgroups.

\begin{verbatim}
    gap> indcyc:= InducedCyclic( th, [ 2 .. NrConjugacyClasses( th ) ], "all" );;
\end{verbatim}

In order to induce characters from the subgroup $D$,
we have to determine the class fusion from $D$ to $Th$.
For that, we use element orders, centralizer orders, power maps,
and the fact that the restrictions of the known characters of $Th$
are characters of $D$.
Since the table of $Th$ will have a table automorphism that swaps
exactly the classes of element order $31$,
we may choose a fixed such class as the image of one class of element order
$31$ in $D$.

\begin{verbatim}
    gap> fus:= InitFusion( d, th );
    [ 1, 2, 2, 6, [ 6, 7 ], [ 6, 7 ], 13, [ 13, 14 ], [ 13, 14 ], 5, 9, 
      [ 19, 20, 21, 22 ], 3, 10, [ 19, 20 ], [ 19, 20 ], [ 9, 10 ], 
      [ 32, 33, 34, 35 ], [ 32, 33, 34, 35 ], 8, 18, [ 25, 26 ], [ 40, 41 ], 
      [ 25, 26 ], [ 40, 41 ], 12, 24, 24, 39, 31, 12, 24, 24, 39, 31, [ 42, 43 ], 
      [ 42, 43 ], [ 42, 43 ], [ 42, 43 ], [ 42, 43 ], [ 42, 43 ] ]
    gap> Positions( OrdersClassRepresentatives( d ), 31 );
    [ 36, 37, 38, 39, 40, 41 ]
    gap> fus[36];
    [ 42, 43 ]
    gap> fus[36]:= 42;;
    gap> TestConsistencyMaps( ComputedPowerMaps( d ), fus,
    >        ComputedPowerMaps( th ) );
    true
    gap> possfus:= FusionsAllowedByRestrictions( d, th, Irr( d ), indcyc, fus,
    >      rec( maxlen:= 10, minamb:= 1, maxamb:= 10^6, quick:= false,
    >           contained:= ContainedPossibleCharacters ) );;
    gap> possfus:= RepresentativesFusions( d, possfus, Group( () ) );
    [ [ 1, 2, 2, 6, 7, 6, 13, 14, 13, 5, 9, 22, 3, 10, 19, 20, 10, 33, 32, 8, 18, 
          25, 40, 26, 41, 12, 24, 24, 39, 31, 12, 24, 24, 39, 31, 42, 43, 42, 42, 
          43, 43 ], 
      [ 1, 2, 2, 6, 7, 7, 13, 14, 14, 5, 9, 22, 3, 10, 19, 20, 10, 33, 32, 8, 18, 
          25, 40, 26, 41, 12, 24, 24, 39, 31, 12, 24, 24, 39, 31, 42, 43, 42, 42, 
          43, 43 ] ]
\end{verbatim}

We get two solutions, up to symmetries of the character table of $D$.
Using that characters induced from $D$ must have integral norm,
one of these candidates gets excluded.

\begin{verbatim}
    gap> indd:= InducedClassFunctionsByFusionMap( d, th, Irr( d ), possfus[1] );;
    gap> ForAll( indd, x -> IsInt( ScalarProduct( th, x, x ) ) );
    false
    gap> indd:= InducedClassFunctionsByFusionMap( d, th, Irr( d ), possfus[2] );;
    gap> ForAll( indd, x -> IsInt( ScalarProduct( th, x, x ) ) );
    true
\end{verbatim}

We initialize the list of known irreducibles (with the trivial character),
and reduce the induced characters.
Applying the LLL algorithm to the reduced characters yields four new
irreducibles.

\begin{verbatim}
    gap> irr:= [ TrivialCharacter( th ) ];;
    gap> red:= ReducedOrdinary( th, irr, Concatenation( indcyc, indd ) );;
    gap> lll:= LLL( th, red.remainders );;
    gap> Length( lll.irreducibles );
    4
    gap> Append( irr, lll.irreducibles );
\end{verbatim}

We create symmetrizations and tensor products of the known irreducibles,
and apply LLL again.
This yields three new irreducibles.

\begin{verbatim}
    gap> sym:= Concatenation( List( [ 2, 3, 4, 5 ],
    >              p -> Symmetrizations( th, irr, p ) ) );;
    gap> sym:= ReducedOrdinary( th, irr, sym );;
    gap> ten:= Set( Tensored( irr, irr ) );;
    gap> ten:= ReducedOrdinary( th, irr, ten );;
    gap> lll:= LLL( th, Concatenation( lll.remainders, sym.remainders,
    >                               ten.remainders ) );;
    gap> Length( lll.irreducibles );
    3
    gap> Append( irr, lll.irreducibles );
    gap> DimensionsMat( irr );
    [ 8, 48 ]
\end{verbatim}

The missing $40$ irreducibles are found in one step,
by computing the possible orthogonal embeddings of the LLL-reduced
lattice of virtual characters into the $40$-dimensional standard lattice.

(In order to accelerate these computations,
we create a new LLL-reduced lattice.)


\begin{verbatim}
    gap> indcyc:= ReducedOrdinary( th, irr, indcyc );;
    gap> indd:= ReducedOrdinary( th, irr, indd );;
    gap> sym:= ReducedOrdinary( th, irr, sym.remainders );;
    gap> ten:= ReducedOrdinary( th, irr, ten.remainders );;
    gap> lll:= LLL( th, Concatenation( indcyc.remainders, indd.remainders,
    >                       sym.remainders, ten.remainders ) );;
    gap> gram:= MatScalarProducts( th, lll.remainders, lll.remainders );;
    gap> emb:= OrthogonalEmbeddings( gram, 40 );;
    gap> Length( emb.solutions );
    4
\end{verbatim}

Two of the four solutions do not satisfy the condition that the standard
basis vectors are irreducible characters.

\begin{verbatim}
    gap> dec:= List( emb.solutions,
    >                x -> Decreased( th, lll.remainders, emb.vectors{ x } ) );;
    gap> dec:= Filtered( dec, x -> x <> fail );;
    gap> Length( dec );
    2
\end{verbatim}

We check whether the first solution yields a character table for $Th$
that is permutation equivalent to the character table in
{\GAP}'s library of character tables,
which is equal to the table that is shown in~\cite[p.~176]{CCN85}.

\begin{verbatim}
    gap> SetIrr( th, List( Concatenation( irr, dec[1].irreducibles ),
    >                      x -> Character( th, x ) ) );
    gap> IsRecord( TransformingPermutationsCharacterTables( th, lib ) );
    true
\end{verbatim}

We do the same for the second solution.

\begin{verbatim}
    gap> ResetFilterObj( th, HasIrr );
    gap> SetIrr( th, List( Concatenation( irr, dec[2].irreducibles ),
    >                      x -> Character( th, x ) ) );
    gap> IsRecord( TransformingPermutationsCharacterTables( th, lib ) );
    true
\end{verbatim}

We see that both solutions are permutation equivalent
to the {\ATLAS} table of $Th$.

\section{The character table of $J_4$ (January 1st, 2016)}

The character table of the sporadic simple Janko group $J_4$
is shown in the {\ATLAS} of Finite Groups (see~\cite[pp.~188--189]{CCN85}),
and this {\ATLAS} table is contained in
the character table library~\cite{CTblLib}
of the computer algebra system {\GAP}~\cite{GAP483}.

We start with the description of the conjugacy classes of $J_4$
as given in~\cite{Jan76} (see Section~\ref{sectclassesJ4}),
then we compute the character table of a subgroup of the type $2^{11}:M_{24}$
in $J_4$ (see Section~\ref{sectsubgroupJ4}),
then we write down the power maps of $J_4$
(see Section~\ref{sectpowermapsJ4}),
and finally we compute the irreducible characters of $J_4$
(see Section~\ref{sectirreduciblesJ4}).

\subsection{The conjugacy classes of $J_4$}\label{sectclassesJ4}

The conjugacy classes of elements of $J_4$ are listed in~\cite[Table~I]{Jan76}.

Up to the permutation $(x_{14}, x_{15})(x_{21}, x_{22})$,
this description of the classes agrees with the element orders
and centralizer orders in the character table that is claimed for $J_4$
in {\GAP}'s Character Table Library~\cite{CTblLib}.

(The permutation is obvious for $x_{14}$ and $x_{15}$ because of the
centralizer orders.
For $x_{21}$ and $x_{22}$,
the flip is necessary in order to get a table with the same third power map
as in the {\GAP} table.)

\begin{verbatim}
    gap> lib:= CharacterTable( "J4" );;
    gap> pos:= [ 1 .. NrConjugacyClasses( lib ) ];;
    gap> orders:= OrdersClassRepresentatives( lib );;
    gap> centralizers:= SizesCentralizers( lib );;
    gap> descr:= TransposedMat( [ pos, orders, centralizers ] );;
    gap> for entry in descr do
    >      Print( String( entry[1], 2 ), "  ",
    >             String( entry[2], 2 ), "  ",
    >             StringPP( entry[3] ), "\n" );
    >    od;
     1   1  2^21*3^3*5*7*11^3*23*29*31*37*43
     2   2  2^21*3^3*5*7*11
     3   2  2^19*3^2*5*7*11
     4   3  2^8*3^3*5*7*11
     5   4  2^15*3*5*11
     6   4  2^15*3
     7   4  2^11*3*7
     8   5  2^6*3*5*7
     9   6  2^8*3^3*5*7*11
    10   6  2^8*3^2
    11   6  2^8*3^2
    12   7  2^3*3*5*7
    13   7  2^3*3*5*7
    14   8  2^8*5
    15   8  2^8*3
    16   8  2^9
    17  10  2^6*3*5
    18  10  2^4*5
    19  11  2^3*3*11^3
    20  11  2*11^2
    21  12  2^6*3
    22  12  2^6*3
    23  12  2^4*3
    24  14  2^2*3*7
    25  14  2^2*3*7
    26  14  2^3*7
    27  14  2^3*7
    28  15  2*3*5
    29  16  2^5
    30  20  2^5*5
    31  20  2^5*5
    32  21  2*3*7
    33  21  2*3*7
    34  22  2^3*3*11
    35  22  2*11
    36  23  23
    37  24  2^4*3
    38  24  2^4*3
    39  28  2^2*7
    40  28  2^2*7
    41  29  29
    42  30  2*3*5
    43  31  31
    44  31  31
    45  31  31
    46  33  2*3*11
    47  33  2*3*11
    48  35  5*7
    49  35  5*7
    50  37  37
    51  37  37
    52  37  37
    53  40  2^3*5
    54  40  2^3*5
    55  42  2*3*7
    56  42  2*3*7
    57  43  43
    58  43  43
    59  43  43
    60  44  2^2*11
    61  66  2*3*11
    62  66  2*3*11
\end{verbatim}

We create a new character table object for a group of the given order
and with the given element orders and centralizer orders.

\begin{verbatim}
    gap> j4:= rec( UnderlyingCharacteristic:= 0,
    >              OrdersClassRepresentatives:= orders,
    >              SizesCentralizers:= centralizers,
    >              Size:= centralizers[1] );;
    gap> ConvertToCharacterTableNC( j4 );;
\end{verbatim}

\subsection{The subgroup $2^{11}:M_{24}$}\label{sectsubgroupJ4}

By~\cite[Theorem~A~(4)]{Jan76},
the group $J_4$ contains a subgroup $U$ that is a split extension
of an elementary abelian group $N$ of order $2^{11}$ by the Mathieu group
$M_{24}$.

E.~O'Brien has computed the character table of this subgroup
with the MAGMA system~\cite{Magma} from a representation of $J_4$.
This character table is permutation equivalent to the character table
that is available in {\GAP}'s character table library,
with the name \verb|"J4M1"|, which stands for ``the first class of maximal
subgroups in $J_4$''.

\begin{verbatim}
    gap> u:= CharacterTable( "J4M1" );
    CharacterTable( "mx1j4" )
\end{verbatim}

(Trying to recompute the character table with {\GAP}'s default algorithm
from a permutation representation on $2^{11}$ points failed
on my small notebook, due to space limitations.)

\subsection{The power maps of $J_4$}\label{sectpowermapsJ4}

Before we can compute the irreducible characters of $J_4$ by inducing
characters from $U$ and from cyclic subgroups of $J_4$,
we determine the power maps of $J_4$.

In addition to the conditions that are imposed by the representative orders
and centralizer orders,
we use the following information from~\cite[Table~I]{Jan76};
each entry \verb|[ c, p, i ]| means that the \verb|p|-th power of the \verb|c|-th class
is the \verb|i|-th class.
(Note that the classes $14$ and $21$ must be swapped with $15$ and $22$,
respectively, in our character table,
compared with the classes description from~\cite[Table~I]{Jan76}.)

\begin{verbatim}
    gap> powinfo:= [
    >     [  5,  2,  2 ],
    >     [  6,  2,  2 ],
    >     [  7,  2,  3 ],
    >     [ 10,  3,  2 ],
    >     [ 11,  3,  3 ],
    >     [ 15,  2,  6 ],
    >     [ 16,  2,  6 ],
    >     [ 17,  5,  2 ],
    >     [ 18,  5,  3 ],
    >     [ 21,  3,  5 ],
    >     [ 22,  3,  6 ],
    >     [ 23,  3,  7 ],
    >     [ 24,  7,  2 ],
    >     [ 25,  7,  2 ],
    >     [ 26,  7,  3 ],
    >     [ 27,  7,  3 ],
    >     [ 29,  2, 16 ],
    >     [ 34, 11,  2 ],
    >     [ 35, 11,  3 ],
    >     ];;
\end{verbatim}

Next we enter information about Galois conjugation;
we will need $p$-th power maps for primes $p$ up to the largest
element order in $J_4$.

\begin{verbatim}
    gap> maxorder:= Maximum( OrdersClassRepresentatives( j4 ) );
    66
    gap> primes:= Filtered( [ 1 .. maxorder ], IsPrimeInt );;
\end{verbatim}

The three classes of element order $6$ are rational;
note that the classes of $x_{10}$ and $x_{11}$ have the same size,
but their third powers are different.
Analogously, the three classes of element order $12$ are rational,
because their third powers lie in the three different rational classes
of element order $4$.

\begin{verbatim}
    gap> pos:= Union( Positions( orders, 6 ), Positions( orders, 12 ) );
    [ 9, 10, 11, 21, 22, 23 ]
    gap> for p in primes do
    >      if 6 mod p <> 0 then
    >        for i in pos do
    >          Add( powinfo, [ i, p, i ] );
    >        od;
    >      fi;
    >    od;
\end{verbatim}

Without loss of generality,
we may choose $x_{24}^2 \sim x_{12}$, $x_{25}^2 \sim x_{13}$,
$x_{26}^2 \sim x_{12}$, and $x_{27}^2 \sim x_{13}$.
Analogously, we may choose
$x_{55}^2 \sim x_{32}$, and $x_{56}^2 \sim x_{33}$.

The squares of $x_{39}$ and $x_{40}$ can be chosen as conjugates of
$x_{26}$ and $x_{27}$, respectively,
since the $14$-th powers of $x_{39}$ and $x_{40}$ are conjugate to $x_3$.

When we choose $x_{55}^3 \sim x_{25}$ and $x_{56}^3 \sim x_{24}$,
we have to choose $x_{32}^3 \sim x_{13}$ and $x_{33}^3 \sim x_{12}$,
in order to get compatible $6$-th powers.
(We choose these powers in order to get the same third power map
as in the {\GAP} table.)

\begin{verbatim}
    gap> Add( powinfo, [ 24, 2, 12 ] );
    gap> Add( powinfo, [ 25, 2, 13 ] );
    gap> Add( powinfo, [ 26, 2, 12 ] );
    gap> Add( powinfo, [ 27, 2, 13 ] );
    gap> Add( powinfo, [ 39, 2, 26 ] );
    gap> Add( powinfo, [ 40, 2, 27 ] );
    gap> Add( powinfo, [ 55, 2, 32 ] );
    gap> Add( powinfo, [ 56, 2, 33 ] );
    gap> Add( powinfo, [ 55, 3, 25 ] );
    gap> Add( powinfo, [ 56, 3, 24 ] );
    gap> Add( powinfo, [ 32, 3, 13 ] );
    gap> Add( powinfo, [ 33, 3, 12 ] );
\end{verbatim}

For the non-rational classes, it is more suitable to set the power map entries
directly, instead of using the list \verb|powinfo|.
Thus we first initialize the power maps and then evaluate this list.

\begin{verbatim}
    gap> powermaps:= [];;
    gap> for p in primes do
    >      powermaps[p]:= InitPowerMap( j4, p );
    >    od;
    gap> for entry in powinfo do
    >      p:= entry[2];
    >      pow:= powermaps[p];
    >      src:= entry[1];
    >      trg:= entry[3];
    >      if IsInt( pow[ src ] ) then
    >        if pow[ src ] <> trg then
    >          Error( "contradiction!" );
    >        fi;
    >      elif not trg in pow[ src ] then
    >        Error( "contradiction!" );
    >      else
    >        pow[ src ]:= trg;
    >      fi;
    >    od;
    gap> SetComputedPowerMaps( j4, powermaps );
\end{verbatim}

There are two Galois conjugate conjugacy classes of elements of order $7$
in $J_4$,
that is, the character values attained on these classes lie in the unique
quadratic extension field $F$, say, of the rationals.
The field $F$ is generated by $\sqrt{-7}$.
If $p$ is a prime different from $7$
then the $p$-th power map swaps the two classes
if and only if $\sqrt{-7}$ differs from the Galois image $\sqrt{-7}^{\ast p}$,
otherwise the $p$-th power map fixes the two classes.
The same holds for the $p$-th powers of the elements of order
$n \in \{ 14, 21, 28, 35, 42 \}$ in $J_4$ if $p$ does not divide $n$.
The function \verb|setGaloisInfo| that has been introduced in
Section~\ref{sectpowermapsTh} will be used for setting the power map info
in these cases.

(As we have used already above,
the distribution of the four classes of element order $14$ into two pairs
is determined by the known seventh powers.)

\begin{verbatim}
    gap> x:= Sqrt( -7 );;
    gap> pos:= Positions( orders, 7 );
    [ 12, 13 ]
    gap> setGaloisInfo( powermaps, pos, orders, primes, x );
    gap> setGaloisInfo( powermaps, [ 24, 25 ], orders, primes, x );
    gap> setGaloisInfo( powermaps, [ 26, 27 ], orders, primes, x );
    gap> pos:= Positions( orders, 21 );
    [ 32, 33 ]
    gap> setGaloisInfo( powermaps, pos, orders, primes, x );
    gap> pos:= Positions( orders, 28 );
    [ 39, 40 ]
    gap> setGaloisInfo( powermaps, pos, orders, primes, x );
    gap> pos:= Positions( orders, 35 );
    [ 48, 49 ]
    gap> setGaloisInfo( powermaps, pos, orders, primes, x );
    gap> pos:= Positions( orders, 42 );
    [ 55, 56 ]
    gap> setGaloisInfo( powermaps, pos, orders, primes, x );
\end{verbatim}

Without loss of generality,
we choose the representatives $x_{48}$, $x_{49}$ in such a way that
their $5$-th powers are $x_{13}$ and $x_{12}$, respectively.

\begin{verbatim}
    gap> powermaps[5]{ [ 48, 49 ] }:= [ 13, 12 ];;
\end{verbatim}

There are four possibilities for the two classes of element order $33$.
Either they are rational, or the character values lie in the quadratic
number field generated by $\sqrt{-3}$, $\sqrt{-11}$, or $\sqrt{33}$.
We check for which of them the induction of characters from the cyclic
subgroup can yield characters of $J_4$.

\begin{verbatim}
    gap> pos:= Positions( orders, 33 );
    [ 46, 47 ]
    gap> setGaloisInfo( powermaps, pos, orders, primes, 1 );
    gap> ind:= InducedCyclic( j4, [ 46 ], "all" );;
    gap> ForAll( ind, x -> IsInt( ScalarProduct( j4, x, x ) ) );
    false
    gap> setGaloisInfo( powermaps, pos, orders, primes, Sqrt( -3 ) );
    gap> ind:= InducedCyclic( j4, [ 46 ], "all" );;
    gap> ForAll( ind, x -> IsInt( ScalarProduct( j4, x, x ) ) );
    false
    gap> setGaloisInfo( powermaps, pos, orders, primes, Sqrt( -11 ) );
    gap> ind:= InducedCyclic( j4, [ 46 ], "all" );;
    gap> ForAll( ind, x -> IsInt( ScalarProduct( j4, x, x ) ) );
    false
    gap> setGaloisInfo( powermaps, pos, orders, primes, Sqrt( 33 ) );
    gap> ind:= InducedCyclic( j4, [ 46 ], "all" );;
    gap> ForAll( ind, x -> IsInt( ScalarProduct( j4, x, x ) ) );
    true
\end{verbatim}

We see that the character values must lie in the field generated by
$\sqrt{33}$.
This implies that also the two classes of element order $66$ are
Galois conjugate, and the character values lie in the same field.
Moreover, we may choose $x_{46}$ and $x_{47}$ as the squares of $x_{61}$
and $x_{62}$, respectively.

\begin{verbatim}
    gap> setGaloisInfo( powermaps, pos, orders, primes, Sqrt( 33 ) );
    gap> pos:= Positions( orders, 66 );
    [ 61, 62 ]
    gap> setGaloisInfo( powermaps, pos, orders, primes, Sqrt( 33 ) );
    gap> powermaps[2]{ pos }:= [ 46, 47 ];;
\end{verbatim}

We show that the two classes of element order $20$ cannot be real,
in the same way as we did for the elements of order $33$.
Since the subgroup $U$ of $J_4$ contains two Galois conjugate classes
of element order $20$, with character values in the quadratic number field
generated by $\sqrt{5}$,
the two classes of element order $20$ in $J_4$ have the same property.

\begin{verbatim}
    gap> pos:= Positions( orders, 20 );
    [ 30, 31 ]
    gap> setGaloisInfo( powermaps, pos, orders, primes, 1 );
    gap> ind:= InducedCyclic( j4, [ 30 ], "all" );;
    gap> ForAll( ind, x -> IsInt( ScalarProduct( j4, x, x ) ) );
    false
    gap> u:= CharacterTable( "J4M1" );
    CharacterTable( "mx1j4" )
    gap> pos:= Positions( OrdersClassRepresentatives( u ), 20 );
    [ 60, 61 ]
    gap> flds:= List( pos, i -> Field( List( Irr( u ), x -> x[i] ) ) );
    [ NF(5,[ 1, 4 ]), NF(5,[ 1, 4 ]) ]
    gap> x:= Sqrt(5);;
    gap> ForAll( flds, f -> x in f );
    true
    gap> setGaloisInfo( powermaps, Positions( orders, 20 ), orders, primes, x );
\end{verbatim}

As a consequence, also the two classes of element order $40$ must be
Galois conjugate,
the character values on these classes must lie in the same quadratic field,
and we may choose $x_{31}$, $x_{30}$ as the squares of $x_{53}$, $x_{54}$,
respectively.
(Again, we choose these images in order to get the same power maps as in
the {\GAP} table.)


\begin{verbatim}
    gap> pos:= Positions( orders, 40 );
    [ 53, 54 ]
    gap> setGaloisInfo( powermaps, pos, orders, primes, Sqrt( 5 ) );
    gap> powermaps[2]{ pos }:= [ 31, 30 ];;
\end{verbatim}

For each element order $p \in \{ 31, 37, 43 \}$,
there are three Galois conjugate conjugacy classes in $J_4$.
Thus the character values lie in the unique cubic subfield of the
field of $p$-th roots of unity in each of these cases.
Without loss of generality,
we choose the representatives in such a way that

\[
   \begin{array}{cc}
      x_{43}^{\ast 5} \sim x_{44}, & x_{44}^{\ast 5} \sim x_{45}, \\
      x_{50}^{\ast 2} \sim x_{51}, & x_{51}^{\ast 2} \sim x_{52}, \\
      x_{57}^{\ast 6} \sim x_{58}, & x_{58}^{\ast 6} \sim x_{59}
   \end{array}
\]

hold.
Note that the three cubic field extensions are generated by the
algebraic integer $c_p$ (see~\cite[p.~xxvii]{CCN85}),
which can be created with the function \verb|EC| in {\GAP}.

\begin{verbatim}
    gap> x:= EC( 31 );;
    gap> classes:= [ 43 .. 45 ];;
    gap> vals:= List( [ 1, 5, 25 ], k -> GaloisCyc( x, k ) );;
    gap> for p in primes do
    >   if p mod 31 <> 0 then
    >     for i in [ 1 .. 3 ] do
    >       powermaps[p][ classes[i] ]:=
    >         classes[ Position( vals, GaloisCyc( vals[i], p ) ) ];
    >     od;
    >   fi;
    > od;
    gap> x:= EC( 37 );;
    gap> classes:= [ 50 .. 52 ];;
    gap> vals:= List( [ 1, 2, 4 ], k -> GaloisCyc( x, k ) );;
    gap> for p in primes do
    >   if p mod 37 <> 0 then
    >     for i in [ 1 .. 3 ] do
    >       powermaps[p][ classes[i] ]:=
    >         classes[ Position( vals, GaloisCyc( vals[i], p ) ) ];
    >     od;
    >   fi;
    > od;
    gap> x:= EC( 43 );;
    gap> classes:= [ 57 .. 59 ];;
    gap> vals:= List( [ 1, 6, 36 ], k -> GaloisCyc( x, k ) );;
    gap> for p in primes do
    >   if p mod 43 <> 0 then
    >     for i in [ 1 .. 3 ] do
    >       powermaps[p][ classes[i] ]:=
    >         classes[ Position( vals, GaloisCyc( vals[i], p ) ) ];
    >     od;
    >   fi;
    > od;
\end{verbatim}

What information is missing now?

\begin{verbatim}
    gap> pos:= PositionsProperty( powermaps[2], IsList );
    [ 21, 22, 23, 35, 37, 38 ]
    gap> orders{ pos };
    [ 12, 12, 12, 22, 24, 24 ]
\end{verbatim}

The square of $x_{35}$ is not yet determined.
Since the centralizer order of each of the two classes of element order $11$
is even, both $x_{19}$ and $x_{20}$ must have square roots of order $22$,
and since $x_{34}^2$ is conjugate to $x_{19}$,
we know that $x_{35}^2$ is conjugate to $x_{20}$.

\begin{verbatim}
    gap> powermaps[2]{ [ 34, 35 ] };
    [ 19, [ 19, 20 ] ]
    gap> powermaps[2][35]:= 20;;
\end{verbatim}

Since $x_{23}$ is the only representative of order $12$
whose $6$-th power is conjugate to $x_3$,
and since $x_{11}$ is the only representative of order $6$
whose cube is conjugate to $x_3$,
we know that the square of $x_{23}$ is conjugate to $x_{11}$.
Conversely, $x_{11}$ cannot be conjugate to the squares of $x_{21}$
or $x_{22}$.

\begin{verbatim}
    gap> powermaps[2]{ [ 21, 22, 23 ] }:= [ [ 9, 10 ], [ 9, 10 ], 11 ];;
\end{verbatim}

Because of $x_{37}^3 \sim x_{38}^3 \sim x_{14}$ and $x_{14}^2 \sim x_6$,
and because the only representative of order $12$ whose third power is
conjugate to $x_6$ is $x_{21}$, we conclude that the squares of
$x_{37}$ and $x_{38}$ are conjugate to $x_{21}$.

\begin{verbatim}
    gap> powermaps[2]{ [ 37, 38 ] }:= [ 22, 22 ];;
\end{verbatim}

Now the three open questions are
whether the classes of element order $24$ are rational or not,
and whether the squares of $x_{21}$ and $x_{22}$ are conjugate to
$x_9$ or $x_{10}$.

First we show that the squares of both $x_{21}$ and $x_{22}$ are conjugate
to $x_{10}$,
since the other three possibilities do not admit a class fusion from $U$.

\begin{verbatim}
    gap> poss:= [];;
    gap> for cand in [ [ 9, 9 ], [ 9, 10 ], [ 10, 9 ] ] do
    >   powermaps[2]{ [ 21, 22 ] }:= cand;
    >   fus:= InitFusion( u, j4 );
    >   TestConsistencyMaps( ComputedPowerMaps( u ), fus, powermaps );
    >   indcyc:= InducedCyclic( j4, [ 21, 22 ], "all" );
    >   possfus:= FusionsAllowedByRestrictions( u, j4, Irr( u ), indcyc, fus,
    >        rec( maxlen:= 10, minamb:= 1, maxamb:= 10^6, quick:= false,
    >             contained:= ContainedPossibleCharacters ) );
    >   Add( poss, Length( possfus ) );
    > od;
    gap> poss;
    [ 0, 0, 0 ]
    gap> powermaps[2]{ [ 21, 22 ] }:= [ 10, 10 ];;
\end{verbatim}

Next we show that the classes of $x_{37}$, $x_{38}$ cannot be rational.

\begin{verbatim}
    gap> pos:= Positions( orders, 24 );
    [ 37, 38 ]
    gap> setGaloisInfo( powermaps, pos, orders, primes, 1 );
    gap> indcyc:= InducedCyclic( j4, pos, "all" );;
    gap> ForAll( indcyc, x -> IsInt( ScalarProduct( j4, x, x ) ) );
    false
\end{verbatim}

Note that $U$ contains two Galois conjugate classes of element order $24$,
with character values in the field generated by $\sqrt{3}$.
Thus also the character values in the classes of $x_{37}$, $x_{38}$ of $J_4$
lie in this field.

\begin{verbatim}
    gap> setGaloisInfo( powermaps, pos, orders, primes, Sqrt( 3 ) );
\end{verbatim}

\subsection{The irreducible characters of $J_4$}\label{sectirreduciblesJ4}

We start with creating some characters of $J_4$ by inducing all linear
characters of cyclic subgroups.

\begin{verbatim}
    gap> indcyc:= InducedCyclic( j4, [ 2 .. NrConjugacyClasses( j4 ) ], "all" );;
\end{verbatim}

In order to induce characters from the subgroup $U$,
we have to determine the class fusion of $U$ in $J_4$.
For that, we use element orders, centralizer orders, power maps,
and the fact that the restrictions of the known characters of $J_4$
are characters of $U$.

Note that we may use also the character table automorphisms of $U$
in those cases where a pair of Galois conjugate classes of $U$ is mapped to
a pair of Galois conjugate classes of $J_4$.

\begin{verbatim}
    gap> u:= CharacterTable( "J4M1" );
    CharacterTable( "mx1j4" )
    gap> fus:= InitFusion( u, j4 );
    [ 1, [ 2, 3 ], 2, [ 2, 3 ], [ 2, 3 ], [ 2, 3 ], [ 2, 3 ], 4, 4, 5, [ 5, 6 ], 
      [ 5, 6 ], 7, [ 5, 6 ], [ 5, 6, 7 ], [ 5, 6, 7 ], [ 5, 6, 7 ], [ 5, 6, 7 ], 
      [ 5, 6, 7 ], [ 5, 6, 7 ], [ 5, 6, 7 ], [ 5, 6, 7 ], 8, 9, [ 9, 10, 11 ], 
      [ 9, 10, 11 ], [ 9, 10, 11 ], [ 9, 10, 11 ], [ 9, 10, 11 ], [ 9, 10, 11 ], 
      [ 9, 10, 11 ], [ 12, 13 ], [ 12, 13 ], 15, 16, [ 14, 15, 16 ], 
      [ 14, 15, 16 ], [ 14, 15, 16 ], [ 14, 15, 16 ], 17, [ 17, 18 ], [ 17, 18 ], 
      [ 17, 18 ], [ 19, 20 ], [ 21, 22 ], [ 21, 22 ], [ 21, 22, 23 ], 
      [ 21, 22, 23 ], [ 21, 22, 23 ], [ 21, 22, 23 ], [ 21, 22, 23 ], 
      [ 21, 22, 23 ], [ 26, 27 ], [ 26, 27 ], [ 24, 25, 26, 27 ], 
      [ 24, 25, 26, 27 ], 28, 28, 29, [ 30, 31 ], [ 30, 31 ], [ 32, 33 ], 
      [ 32, 33 ], [ 34, 35 ], 36, 36, [ 37, 38 ], [ 37, 38 ], [ 39, 40 ], 
      [ 39, 40 ], 42, 42 ]
    gap> Print( AutomorphismsOfTable( u ), "\n" );
    Group( [ (67,68), (65,66), (60,61), (57,58)(71,72), (57,58)(67,68)(71,72), 
      (57,58)(60,61)(71,72), (32,33)(53,54)(55,56)(62,63)(69,70), 
      ( 5, 6)(15,16)(21,22)(30,31)(42,43)(49,50)(51,52) ] )
    gap> fus{ [ 60, 61, 67, 68, 69, 70 ] };
    [ [ 30, 31 ], [ 30, 31 ], [ 37, 38 ], [ 37, 38 ], [ 39, 40 ], [ 39, 40 ] ]
    gap> fus[60]:= 30;;
    gap> fus[67]:= 37;;
    gap> fus[69]:= 40;;
    gap> TestConsistencyMaps( ComputedPowerMaps( u ), fus,
    >        ComputedPowerMaps( j4 ) );
    true
\end{verbatim}

We get $1440$ possible class fusions with the property that
all those characters of $J_4$ that are induced from cyclic subgroups
restrict to characters of $U$,
in $720$ orbits under table automorphisms of $U$.


\begin{verbatim}
    gap> irr:= [ TrivialCharacter( j4 ) ];;
    gap> indcyc:= ReducedOrdinary( j4, irr, indcyc );;
    gap> possfus:= FusionsAllowedByRestrictions( u, j4, Irr( u ),
    >      indcyc.remainders, fus,
    >      rec( maxlen:= 10, minamb:= 1, maxamb:= 10^3, quick:= false,
    >           contained:= ContainedPossibleCharacters ) );;
    gap> Length( possfus );
    1440
    gap> reps:= RepresentativesFusions( u, possfus, Group(()) );;
    gap> Length( reps );
    720
\end{verbatim}

Only two of these fusion candidates have the property that all characters
induced from $U$ to $J_4$ have integral norms.
All scalar products between these induced characters are integral for
exactly one candidate.
Thus we have determined the class fusion up to table automorphisms of $U$.

\begin{verbatim}
    gap> reps:= Filtered( reps,
    >      map -> ForAll( InducedClassFunctionsByFusionMap( u, j4, Irr(u), map ),
    >                     x -> IsPosInt( ScalarProduct( j4, x, x ) ) ) );;
    gap> Length( reps );
    2
    gap> inds:= List( reps,
    >      map -> InducedClassFunctionsByFusionMap( u, j4, Irr( u ), map ) );;
    gap> ForAll( Flat( MatScalarProducts( j4, inds[1], inds[1] ) ), IsInt );
    false
    gap> ForAll( Flat( MatScalarProducts( j4, inds[2], inds[2] ) ), IsInt );
    true
\end{verbatim}

We reduce the induced characters with the trivial character.
Applying the LLL algorithm to the reduced characters yields $29$ new
irreducibles.

\begin{verbatim}
    gap> ind:= ReducedOrdinary( j4, irr, inds[2] );;
    gap> Length( ind.irreducibles );
    0
    gap> lll:= LLL( j4, Concatenation( indcyc.remainders, ind.remainders ) );;
    gap> Length( lll.irreducibles );
    29
    gap> Append( irr, lll.irreducibles );
\end{verbatim}

The remaining LLL-reduced lattice is spanned by vectors of norm $2$.
We compute irreducible characters from sublattices of the types $D_4$
and $D_5$.

\begin{verbatim}
    gap> lll.norms;
    [ 2, 2, 2, 2, 2, 2, 2, 2, 2, 2, 2, 2, 2, 2, 2, 2, 2, 2, 2, 2, 2, 2, 2, 2, 2, 
      2, 2, 2, 2, 2, 2, 2 ]
    gap> dn:= DnLatticeIterative( j4, lll.remainders );;
    gap> Length( dn.irreducibles );
    28
    gap> Append( irr, dn.irreducibles );
\end{verbatim}

Now just four irreducibles are missing.
We compute the possible orthogonal embeddings of the LLL-reduced
lattice of virtual characters into the $4$-dimensional standard lattice.

\begin{verbatim}
    gap> gram:= MatScalarProducts( j4, dn.remainders, dn.remainders );;
    gap> emb:= OrthogonalEmbeddings( gram, 4 );;
    gap> Length( emb.solutions );
    3
\end{verbatim}

One of the three solutions does not satisfy the condition that the standard
basis vectors are irreducible characters.

\begin{verbatim}
    gap> dec:= List( emb.solutions,
    >                x -> Decreased( j4, dn.remainders, emb.vectors{ x } ) );;
    gap> dec:= Filtered( dec, x -> x <> fail );;
    gap> Length( dec );
    2
\end{verbatim}

The first solution is not compatible with the $2$-nd power map.

\begin{verbatim}
    gap> possirr:= List( dec, x -> x.irreducibles );;
    gap> chi:= possirr[1][1];;
    gap> sym:= Symmetrizations( j4, [ chi ], 2 );;
    gap> ForAll( sym, x -> IsInt( ScalarProduct( j4, x, chi ) ) );
    false
\end{verbatim}

Thus we are left with one solution.
We check whether it yields a character table for $J_4$
that is permutation equivalent to the character table in
{\GAP}'s library of character tables,
which is equal to the table that is shown in~\cite[pp.~188--189]{CCN85}.

\begin{verbatim}
    gap> SetIrr( j4, List( Concatenation( irr, possirr[2] ),
    >                      x -> Character( j4, x ) ) );
    gap> IsRecord( TransformingPermutationsCharacterTables( lib, j4 ) );
    true
\end{verbatim}

We see that the table which we have computed is permutation equivalent
to the {\ATLAS} table of $J_4$.

\section{The character table of $2.{}^2E_6(2)$ (February 29th, 2016)}%
\label{sect22e62}

In the following, we compute the character table of the double cover $2.G$
of the group $G = {}^2E_6(2)$ with character theoretic methods.

For that, we assume the character table of $G$,
the existence of several subgroups of $2.G$,
and the knowledge of their character tables, see Section~\ref{assumptions2G}.
The first step is to determine an approximation of the conjugacy classes
of $2.G$ (see Section~\ref{theclasses}).
Next we determine the class fusions of the subgroups in question
(see Sections~\ref{f42subs} to~\ref{2o10m2subs}),
and finally we use standard techniques
to compute the faithful irreducible characters of $2.G$
from induced characters (see Section~\ref{thecharacters}).

\subsection{Assumptions}\label{assumptions2G}

We assume that the outer automorphism group of $G = {}^2E_6(2)$
is isomorphic to the symmetric group of order six.
Let $\alpha$ be an outer automorphism of order three
and $\beta$ be an outer automorphism of order two.

We assume that a perfect central extension $H$ of the structure $2^2.G$
exists,
and that the action of $\langle \alpha, \beta \rangle$
lifts to an action on $H$, which permutes the three central involutions
transitively.
Thus the lift of $\beta$ fixes a unique central subgroup of order two,
and we choose the factor group of $H$ by this subgroup
as the double cover $2.G$ whose character table we are going to construct.
We denote the natural epimorphism from $2.G$ to $G$ by $\pi$.

(The split extension of the chosen double cover $2.G$ by $\beta$ will be
the group $2.G.2$ whose character table will be considered in
Section~\ref{section22e622}.)

We assume that $G$ contains subgroups of the following structures.

\begin{itemize}
\item
    $F_4(2)$,
\item
    $Fi_{22}$,
\item
    $3 \times U_6(2)$ (the centralizer of an element of order three
    in $G$, see~\cite[Table~5.1]{LSS92}),
\item
    $O^-_{10}(2)$ (the centralizer in $G$ of an element of order three
    in $G.3 \setminus G$, see~\cite[Table~5.1]{LSS92}).
\end{itemize}

Furthermore,
we assume that $2.G$ contains subgroups of the structure $2.F_4(2)$,
the double cover of $F_4(2)$.
This fact is stated for example in~\cite[Section~6]{SW99};
I would be happy to add a reference to a proof of the statement.

By~\cite{BMO16},
we can assume the correctness of the character tables of
$G$, $2.F_4(2)$, $Fi_{22}$, $U_6(2)$, and $O^-_{10}(2)$
that are available in the {\GAP} Character Table Library.
(In particular, the character table of $G$ has been computed
with MAGMA~\cite{Magma} from the smallest faithful permutation representation
of $G$.)

\subsection{Outer automorphisms of $G$}

The orbits of $\alpha$ on the conjugacy classes of $G$ are determined
by the character table of $G$,
since the group of character table automorphisms of $G$ contains a unique
subgroup of order three.

\begin{verbatim}
    gap> t:= CharacterTable( "2E6(2)" );;
    gap> t2:= CharacterTable( "2E6(2).2" );;
    gap> aut:= AutomorphismsOfTable( t );;
    gap> Factors( Size( aut ) );
    [ 2, 2, 2, 2, 2, 2, 2, 2, 3 ]
    gap> syl:= SylowSubgroup( aut, 3 );;
    gap> IsNormal( aut, syl );
    true
    gap> orbs:= Orbits( syl, [ 1 .. NrConjugacyClasses( t ) ] );;
    gap> orbsalpha:= List( Filtered( orbs, l -> Length( l ) <> 1 ), Set );
    [ [ 11, 12, 13 ], [ 16, 17, 18 ], [ 39, 40, 41 ], [ 43, 44, 45 ], 
      [ 46, 47, 48 ], [ 64, 65, 66 ], [ 67, 68, 69 ], [ 75, 76, 77 ], 
      [ 78, 79, 80 ], [ 88, 89, 90 ], [ 91, 92, 93 ], [ 94, 95, 96 ], 
      [ 114, 115, 116 ], [ 117, 118, 119 ] ]
\end{verbatim}

The orbits of $\beta$ on the conjugacy classes of $G$ are determined
by the class fusion of $G$ in $G.2$,
except that we have to choose one of the three possible sets of orbits,
which correspond to the three subgroups of index three in the
automorphism group of $G$.
We choose the case where class number $11$ is fixed and
the classes $12$ and $13$ are swapped.


\begin{verbatim}
    gap> tfust2:= PossibleClassFusions( t, t2 );;
    gap> poss:= Set( List( tfust2, l -> Filtered( InverseMap( l ), IsList ) ) );
    [ [ [ 11, 12 ], [ 16, 17 ], [ 39, 40 ], [ 43, 44 ], [ 46, 47 ], [ 55, 56 ], 
          [ 61, 62 ], [ 64, 65 ], [ 67, 68 ], [ 75, 76 ], [ 78, 79 ], [ 88, 89 ], 
          [ 91, 92 ], [ 94, 95 ], [ 99, 100 ], [ 103, 104 ], [ 109, 110 ], 
          [ 114, 115 ], [ 117, 118 ], [ 123, 124 ], [ 125, 126 ] ], 
      [ [ 11, 13 ], [ 16, 18 ], [ 39, 41 ], [ 43, 45 ], [ 46, 48 ], [ 55, 56 ], 
          [ 61, 62 ], [ 64, 66 ], [ 67, 69 ], [ 75, 77 ], [ 78, 80 ], [ 88, 90 ], 
          [ 91, 93 ], [ 94, 96 ], [ 99, 100 ], [ 103, 104 ], [ 109, 110 ], 
          [ 114, 116 ], [ 117, 119 ], [ 123, 124 ], [ 125, 126 ] ], 
      [ [ 12, 13 ], [ 17, 18 ], [ 40, 41 ], [ 44, 45 ], [ 47, 48 ], [ 55, 56 ], 
          [ 61, 62 ], [ 65, 66 ], [ 68, 69 ], [ 76, 77 ], [ 79, 80 ], [ 89, 90 ], 
          [ 92, 93 ], [ 95, 96 ], [ 99, 100 ], [ 103, 104 ], [ 109, 110 ], 
          [ 115, 116 ], [ 118, 119 ], [ 123, 124 ], [ 125, 126 ] ] ]
    gap> orbsbeta:= poss[3];;
\end{verbatim}

\subsection{Tools for determining the conjugacy classes of $2.G$}%
\label{theclasses}

The aim of this section is to provide {\GAP} functions for determining
the class fusion from $2.G$ to $G$,
that is, determining for which conjugacy classes $g^G$ of $G$
the two elements in $\pi^{-1}( g )$ are conjugate in $2.G$ or not.
In the former case, $\pi^{-1}( g^G )$ forms a single conjugacy class of $2.G$;
we say that $g^G$ does not split.
In the latter case, $\pi^{-1}( g^G )$ consists of two conjugacy classes
of $2.G$; we say that $g^G$ splits.

Let $z$ denote the central involution in $2.G$.

\subsubsection{Elementary criteria}\label{elementarysplitting}

If $g \in G$ has odd order $n$, say, then the two preimages of $g$
under $\pi$ have the orders $n$ and $2 n$;
hence they cannot be conjugate.

\begin{verbatim}
    gap> orders:= OrdersClassRepresentatives( t );;
    gap> mustsplit:= PositionsProperty( orders, IsOddInt );
    [ 1, 5, 6, 7, 23, 33, 34, 51, 52, 55, 56, 83, 86, 87, 97, 98, 103, 104, 107, 
      108, 123, 124, 125, 126 ]
\end{verbatim}

If $g \in G$ is self-centralizing then
the elements in $\pi^{-1}( g )$ are not conjugate in $\pi^{-1}( G )$.


\begin{verbatim}
    gap> selfCentralizingClassesSplit:= function( t, mustsplit )
    >      local centralizers, orders, i;
    > 
    >      centralizers:= SizesCentralizers( t );
    >      orders:= OrdersClassRepresentatives( t );
    >      for i in [ 1 .. Length( centralizers ) ] do
    >        if centralizers[i] = orders[i] and not i in mustsplit then
    >          Print( "#I  class ", i, " splits (self-centralizing)\n" );
    >          AddSet( mustsplit, i );
    >        fi;
    >      od;
    >    end;;
    gap> selfCentralizingClassesSplit( t, mustsplit );
    #I  class 109 splits (self-centralizing)
    #I  class 110 splits (self-centralizing)
    #I  class 120 splits (self-centralizing)
    #I  class 122 splits (self-centralizing)
\end{verbatim}

Let $g, h$ be elements of $G$, and $n$ be an odd integer
such that $h^n = g$ holds.
If the elements in $\pi^{-1}( g )$ are not conjugate in $2.G$
then also the elements in $\pi^{-1}( h )$ are not conjugate in $2.G$.

(This is in fact a generalization of the first criterion if one assumes
that the elements in  $\pi^{-1}( 1 )$ are not conjugate.)

\begin{verbatim}
    gap> oddRootsOfSplittingClassesSplit:= function( t, mustsplit )
    >      local powmaps, found, p, map, i;
    > 
    >      powmaps:= ComputedPowerMaps( t );
    >      repeat
    >        found:= false;
    >        for p in [ 1 .. Length( powmaps ) ] do
    >          if p mod 2 = 1 and IsBound( powmaps[p] ) then
    >            map:= powmaps[p];
    >            for i in [ 1 .. Length( map ) ] do
    >              if map[i] in mustsplit and not i in mustsplit then
    >                Print( "#I  class ", i, " splits (",
    >                       Ordinal( p ), " root of ", map[i], ")\n" );
    >                found:= true;
    >                AddSet( mustsplit, i );
    >              fi;
    >            od;
    >          fi;
    >        od;
    >      until found = false;
    >    end;;
\end{verbatim}

Let $U$ be a subgroup of $G$, and $g \in U$.

If the elements of $\pi^{-1}( g )$ are conjugate in $\pi^{-1}( U )$
then they are conjugate in $\pi^{-1}( G )$.

\begin{verbatim}
    gap> notSplittingClassesOfSubgroupDoNotSplit:= function( 2sfuss, sfust,
    >                                                  mustnotsplit )
    >      local new, i;
    > 
    >      new:= sfust{ PositionsProperty( InverseMap( 2sfuss ), IsInt ) };
    >      for i in Set( new ) do
    >        if not i in mustnotsplit then
    >          Print( "#I  class ", i, " does not split (as in subgroup)\n" );
    >        fi;
    >      od;
    >      UniteSet( mustnotsplit, new );
    >    end;;
\end{verbatim}

If $|C_G(g)| / |C_U(g)|$ is odd
and if the elements of $\pi^{-1}( g )$ are not conjugate in $\pi^{-1}( U )$
then they are not conjugate in $\pi^{-1}( G )$.

\begin{verbatim}
    gap> splittingClassesWithOddCentralizerIndexSplit:= function( s, t,
    >        sfust, 2sfuss, mustsplit )
    >      local inv, scents, tcents, i;
    > 
    >      inv:= InverseMap( 2sfuss );
    >      scents:= SizesCentralizers( s );
    >      tcents:= SizesCentralizers( t );
    >      for i in [ 1 .. Length( sfust ) ] do
    >        if IsList( inv[i] ) and
    >           IsOddInt( tcents[ sfust[i] ] / scents[i] ) then
    >          if not sfust[i] in mustsplit then
    >            Print( "#I  class ", sfust[i],
    >                   " splits (odd centralizer index)\n" );
    >            AddSet( mustsplit, sfust[i] );
    >          fi;
    >        fi;
    >      od;
    >      oddRootsOfSplittingClassesSplit( t, mustsplit );
    >    end;;
\end{verbatim}

\subsubsection{Norms of induced characters}\label{normsinduced}

The next criterion tries to exploit the fact that for a given character
of a subgroup, the norm of the induced character is an integer,
and that there are only few possibilities for the contribution of each
conjugacy class of $G$ to this norm.

Let $U$ be a subgroup of $G$
and $\chi$ be a character of $\pi^{-1}( U )$
with the property $\chi( z ) = - \chi( 1 )$;
then $\chi( gz ) = - \chi( g )$ holds for any $g \in \pi^{-1}( U )$.

Fix representatives $u_1, u_2, \ldots, u_n$ of the conjugacy classes of $U$,
choose $g_i \in \pi^{-1}( u_i )$,
and define the class function $\chi'$ of $U$ by $\chi'( u_i ) = \chi( g_i )$.
Consider the induced character $\psi = \chi^{\pi^{-1}(G)}$.

If $g$ and $gz$ are conjugate in $\pi^{-1}(G)$ then $\psi(g) = 0$ holds.
If $g$ and $gz$ are not conjugate in $\pi^{-1}(G)$
then set
\[
   R( \pi(g) ) =
      \{ i; 1 \leq i \leq n, g_i \sim g \mbox{\rm\ or\ } g_i z \sim g \},
\]
where $\sim$ denotes conjugacy in $\pi^{-1}(G)$.
Note that at most one of $g_i$, $g_i z$ can be conjugate to $g$.

Then we have
\[
   \psi(g) = |C_{\pi^{-1}(G)}(g)|
               \cdot \sum_{i \in R( \pi(g) )}
                 \frac{\pm \chi(g_i)}{|C_{\pi^{-1}(U)}(g_i)|},
\]
where $\pm$ means that there is an appropriate choice of signs
for the summands such that the equation holds.

Since $|C_{\pi^{-1}(G)}(g)| = 2 \cdot |C_G( \pi(g) )|$
if $g$ and $gz$ are not conjugate in $\pi^{-1}(G)$
and hence
$|C_{\pi^{-1}(U)}(g_i)| = 2 \cdot |C_U(\pi(g_i))|$, this can be written as
\[
   \psi(g) = |C_G( \pi(g) )|
               \cdot \sum_{i \in R( \pi(g) )}
                 \frac{\pm \chi'(u_i)}{|C_U(u_i)|}.
\]

Now we set $R( \pi(g) ) = \emptyset$ if $g$ and $gz$ are conjugate,
and get
\begin{eqnarray*}
   [ \psi, \psi ] & = & \frac{1}{|\pi^{-1}(G)|}
                          \cdot \sum_{g \in \pi^{-1}(G)} |\psi(g)|^2 \\
                  & = & \sum_{g \in \pi^{-1}(G)/\sim}
                          \frac{|\psi(g)|^2}{2 \cdot |C_G( \pi(g) )|} \\
                  & = & \frac{1}{2} \cdot \sum_{g \in \pi^{-1}(G)/\sim}
                          |C_G( \pi(g) )|
                          \cdot \left| \sum_{i \in R( \pi(g) )}
                          \frac{\pm \chi'(u_i)}{|C_U(u_i)|} \right|^2 \\
                  & = & \sum_{g \in G/\sim_G}
                          |C_G(g)|
                          \cdot \left| \sum_{i \in R(g)}
                          \frac{\pm \chi'(u_i)}{|C_U(u_i)|} \right|^2.
\end{eqnarray*}

Note that last expression on the right hand side is expressed
in terms of $G$.
The fact that $[ \psi, \psi ]$ is an integer can be turned into
a splitting criterion, as follows.

\begin{itemize}
\item
    Classes of $G$ that are known not to split can be ignored
    in the summation, since they contribute zero.
\item
    For a class $g^G$ that is known to split,
    the possible choices of signs yield at most $2^{|R(g)| - 1}$
    different contributions
    $\left| \sum_{i \in R(g)} \frac{\pm \chi'(u_i)}{|C_U(u_i)|} \right|^2$
    to the norm.
    The number can be much smaller, for example if some $\chi'(u_i)$ are zero.
\item
    For those classes $g^G$ where we do not know whether they split,
    we get the analogous sets of possible contributions from hypothetical
    sets $R(g)$ (if the class splits)
    plus the value zero (if the class does not split).
\item
    For given $\chi$, we compute all those combinations of possible
    contributions for which the summation yields an integer.
    If the contribution of a given class $g^G$ to the norm of $\psi$
    is nonzero in all these cases,
    we conclude that this class must split.
    If the contribution is zero in all these cases and if zero cannot be
    obtained as some
    $\sum_{i \in R(g)} \frac{\pm \chi'(u_i)}{|C_U(u_i)|}$,
    we conclude that this class must not split.
\item
    Our strategy is to check all irreducible characters $\chi$
    of $\pi^{-1}(U)$ with the property $\chi(z) = - \chi(1)$,
    ordered by increasing number of combinations to be checked.
\end{itemize}

This leads to the following {\GAP} functions.

The function \verb|contributionData| takes
the character tables of $U$ and $G$,
a list \verb|inv| that describes the class fusion,
a class function (projective character) $\chi'$ of $U$,
and the list \verb|mustsplit|.
The entry at position $i$ of \verb|inv| is bound if and only if
some class of $U$ maps to class $i$ of $G$;
in this cases the entry is an integer if the the class of $U$
at this position is the unique class of $U$ that fuses into class $i$ of $G$,
and otherwise the entry is the list of class positions in $U$ that fuse into
this class; this data format is returned by the {\GAP} function
\verb|InverseMap|.
Positions of not splitting classes of $G$ should be unbound in \verb|inv|,
since these classes are known not to contribute to the norm of the
induced character.

The function returns a record with the components
\verb|size| (the number of combinations to be checked),
\verb|safepart| (the sum of the known contributions to the norm of the induced
character),
\verb|bound| (the list of positions of those classes for which only a list of
possible contributions is known),
\verb|contrib| (the list of value lists for the classes in \verb|bound|), and
\verb|zeroonlyifnonsplit| (if the $i$-th entry is \verb|true| then contribution $0$
at $i$ for all solutions means that the $i$-th class does not split).

\begin{verbatim}
    gap> contributionData:= function( s, t, inv, chiprime, mustsplit )
    >      local contrib, zeroonlyifnonsplit, safepart, n, tcents,
    >            sclasses, i, j, val, choices, signs, cand;
    > 
    >      contrib:= [];
    >      zeroonlyifnonsplit:= [];
    >      safepart:= 0;
    >      n:= 1;
    >      tcents:= SizesCentralizers( t );
    >      sclasses:= SizesConjugacyClasses( s );
    >      for i in [ 1 .. Length( inv ) ] do
    >        if IsBound( inv[i] ) then
    >          # The subgroup contains elements in the 'i'-th class.
    >          if IsInt( inv[i] ) then
    >            # Only one class of the subgroup fuses into the 'i'-th class.
    >            j:= inv[i];
    >            val:= sclasses[j] * chiprime[j];
    >            val:= tcents[i] / Size(s)^2 * val * GaloisCyc( val, -1 );
    >            if not IsInt( val ) then
    >              if i in mustsplit then
    >                # The summand is known, add it to 'safepart'.
    >                safepart:= safepart + val;
    >              else
    >                # The class may or may not split.
    >                # If it splits then 'val' is the contribution to the norm.
    >                contrib[i]:= [ 0, val ];
    >                zeroonlyifnonsplit[i]:= true;
    >                n:= n * 2;
    >              fi;
    >            fi;
    >          else
    >            # Several classes of the subgroup fuse into the 'i'-th class.
    >            choices:= List( inv[i], j -> sclasses[j] * chiprime[j] );
    >            signs:= Tuples( [ 1, -1 ], Length( choices ) );
    >            cand:= signs * choices;
    >            cand:= tcents[i] / Size(s)^2 *
    >                     Set( List( cand, x -> x * GaloisCyc( x, -1 ) ) );
    >            if not ForAll( cand, IsInt ) then
    >              if Length( cand ) = 1 then
    >                if i in mustsplit then
    >                  # We get a contribution to 'safepart'.
    >                  safepart:= safepart + cand[1];
    >                else
    >                  UniteSet( cand, [ 0 ] );
    >                  contrib[i]:= cand;
    >                  zeroonlyifnonsplit[i]:= true;
    >                  n:= n * Length( cand );
    >                fi;
    >              else
    >                if not i in mustsplit then
    >                  if not 0 in cand then
    >                    UniteSet( cand, [ 0 ] );
    >                    zeroonlyifnonsplit[i]:= true;
    >                  fi;
    >                fi;
    >                contrib[i]:= cand;
    >                n:= n * Length( cand );
    >              fi;
    >            fi;
    >          fi;
    >        fi;
    >      od;
    > 
    >      return rec( safepart:= safepart,
    >                  contrib:= contrib,
    >                  size:= n,
    >                  bound:= Filtered( [ 1 .. Length( contrib ) ],
    >                                    x -> IsBound( contrib[x] ) ),
    >                  zeroonlyifnonsplit:= zeroonlyifnonsplit,
    >                );
    >    end;;
\end{verbatim}

The function \verb|integralContributions| runs over all combinations
given by \verb|r.contrib|,
and returns the list of those vectors whose sum plus \verb|r.safepart|
is an integer.

\begin{verbatim}
    gap> integralContributions:= function( r )
    >      local positions, len, images, number, index, direction, initial,
    >            norm, solutions, i;
    > 
    >      # Initialize the counter and the list of solutions.
    >      positions:= r.bound;
    >      len:= Length( positions );
    >      images:= r.contrib{ positions };
    >      number:= List( images, Length );
    >      index:= ListWithIdenticalEntries( len, 1 );
    >      direction:= ShallowCopy( index );  # 1 means up, -1 means down
    >      initial:= List( images, l -> l[1] );
    >      norm:= r.safepart + Sum( initial );
    >      solutions:= [];
    >      if IsInt( norm ) then
    >        solutions[1]:= initial;
    >      fi;
    > 
    >      while true do
    >        # Increase the counter. (Change only one position in each step.)
    >        i:= 1;
    >        while i <= len and
    >              ( ( index[i] = number[i] and direction[i] = 1 ) or
    >                ( index[i] = 1 and direction[i] = -1 ) ) do
    >          direction[i]:= - direction[i];
    >          i:= i+1;
    >        od;
    > 
    >        if len < i then
    >          # We are done.
    >          return solutions;
    >        fi;
    > 
    >        # Update at position 'i'.
    >        norm:= norm - images[i][ index[i] ];
    >        index[i]:= index[i] + direction[i];
    >        norm:= norm + images[i][ index[i] ];
    > 
    >        if IsInt( norm ) then
    >          # We have found a solution.
    >          Add( solutions,
    >               List( [ 1 .. len ], i -> images[i][ index[i] ] ) );
    >        fi;
    >      od;
    >    end;;
\end{verbatim}

The function \verb|evaluateContributions| takes
the record computed by the function \verb|contributions|
and the list of vectors which yield induced class functions
with integral norms, as computed by the function \verb|integralContributions|,
and extends the lists \verb|mustsplit|, \verb|mustnotsplit| whenever possible.

\begin{verbatim}
    gap> evaluateContributions:= function( r, res, sfust,
    >                                      mustsplit, mustnotsplit )
    >      local param, i, c;
    > 
    >      param:= Parametrized( res );
    >      for i in [ 1 .. Length( r.bound ) ] do
    >        c:= r.bound[i];
    >        if param[i] = 0 then
    >          # If contribution zero cannot arise as a sum of values
    >          # then the class cannot split.
    >          if IsBound( r.zeroonlyifnonsplit[c] ) and
    >             r.zeroonlyifnonsplit[c] = true and
    >             not c in mustnotsplit then
    >            Print( "#I  class ", c,
    >                   " does not split (contribution criterion)\n" );
    >            if c in mustsplit then
    >              Error( "contradiction for class ", c );
    >            fi;
    >            AddSet( mustnotsplit, c );
    >          fi;
    >        elif IsRat( param[i] ) then
    >          if not c in mustsplit then
    >            Print( "#I  class ", c, " splits (contribution criterion)\n" );
    >            if c in mustnotsplit then
    >              Error( "contradiction for class ", c );
    >            fi;
    >            AddSet( mustsplit, c );
    >          fi;
    >        elif IsList( param[i] ) and not 0 in param[i] then
    >          # If no zero occurs then the class must split.
    >          if not c in mustsplit then
    >            Print( "#I  class ", c, " splits (contribution criterion)\n" );
    >            if c in mustnotsplit then
    >              Error( "contradiction for class ", c );
    >            fi;
    >            AddSet( mustsplit, c );
    >          fi;
    >        fi;
    >      od;
    >    end;;
\end{verbatim}

The strategy is implemented as follows.

\begin{verbatim}
    gap> computeContributions:= function( s, t, sfust, classfuns, bound,
    >                                     mustsplit, mustnotsplit )
    >      local inv, i, known, candidates, r, res;
    > 
    >      inv:= InverseMap( sfust );
    > 
    >      repeat
    >        for i in mustnotsplit do
    >          # The induced character is zero at the preimage of 'i',
    >          # there is no contribution to the norm.
    >          Unbind( inv[i] );
    >        od;
    >        known:= [ ShallowCopy( mustsplit ), ShallowCopy( mustnotsplit ) ];
    >        candidates:= List( classfuns,
    >            chi -> contributionData( s, t, inv, chi, mustsplit ) );
    >        candidates:= Filtered( candidates, r -> r.size < bound );
    >        SortParallel( List( candidates, r -> r.size ), candidates );
    >        for r in candidates do
    >          res:= integralContributions( r );
    >          if Length( res ) = 0 then
    >            Error( "no solution" );
    >          fi;
    >          evaluateContributions( r, res, sfust, mustsplit, mustnotsplit );
    >          oddRootsOfSplittingClassesSplit( t, mustsplit );
    >        od;
    >      until known = [ mustsplit, mustnotsplit ];
    >    end;;
\end{verbatim}

Note that the function \verb|computeContributions| does not stop when it has
determined for all classes in question whether they split or not,
but processes all given characters.
The idea is that contradictions because of wrong assumptions should
get chances to be noticed.

\subsection{Subgroups of the type $F_4(2)$ in $G$}\label{f42subs}

There are three possible class fusion of $F_4(2)$ in $G$,
up to symmetries of the character table of $F_4(2)$,
and these fusions form one orbit under the action of $\alpha$.

\begin{verbatim}
    gap> s:= CharacterTable( "F4(2)" );;
    gap> fus:= PossibleClassFusions( s, t );;
    gap> rep:= RepresentativesFusions( s, fus, Group( () ) );;
    gap> Length( rep );
    3
    gap> oneorbit:= orbsalpha[1];
    [ 11, 12, 13 ]
    gap> List( rep, map -> Intersection( map, oneorbit ) );
    [ [ 11 ], [ 12 ], [ 13 ] ]
\end{verbatim}

By our assumption that $G$ contains $F_4(2)$ type subgroups,
this means that all three possibilities are really class fusions
of $F_4(2)$ type subgroups in $G$.
We fix three such subgroups with these fusions,
in one orbit under the action of $\alpha$,
and call the subgroups $M_3$, $M_4$, and $M_5$.
(Note that the three classes of maximal subgroups of the type $F_4(2)$
in $G$ claimed in the list in~\cite[pp.~191]{CCN85}
appear in the positions $3$ to $5$.)

\begin{verbatim}
    gap> m3:= s;;  m3fust:= rep[1];;
    gap> m4:= s;;  m4fust:= rep[2];;
    gap> m5:= s;;  m5fust:= rep[3];;
\end{verbatim}

Concerning the structures of the preimages $\pi^{-1}( M_i )$,
for $3 \leq i \leq 5$,
we first note that the preimages of the three subgroups in $H$
are conjugate under the lift of $\alpha$ to $H$.
Since the Schur multiplier of $F_4(2)$ is cyclic of order two,
the structure must be either $2 \times 2.F_4(2)$ or $2^2 \times F_4(2)$,
and since we have assumed the existence of subgroups of the type
$2.F_4(2)$ in $2.G$,
we may choose the subgroups in such a way that the former case applies.
This implies that one of the preimages in question is a direct product
$2 \times F_4(2)$,
and the other two have the structure $2.F_4(2)$.
Without loss of generality, we may choose the subgroups such that
$\pi^{-1}( M_3 ) \cong 2 \times F_4(2)$ holds.

\begin{verbatim}
    gap> 2m3:= CharacterTable( "Cyclic", 2 ) * m3;;
    gap> 2m4:= CharacterTable( "2.F4(2)" );;
    gap> 2m5:= 2m4;;
\end{verbatim}

Let us apply the criteria from Section~\ref{theclasses}
to the subgroups $M_3$, $M_4$, $M_5$ of $G$.

\begin{verbatim}
    gap> splittingClassesWithOddCentralizerIndexSplit( m3, t, m3fust,
    >        GetFusionMap( 2m3, m3 ), mustsplit );
    #I  class 73 splits (odd centralizer index)
    #I  class 85 splits (odd centralizer index)
    #I  class 101 splits (odd centralizer index)
    #I  class 106 splits (odd centralizer index)
    gap> splittingClassesWithOddCentralizerIndexSplit( m4, t, m4fust,
    >        GetFusionMap( 2m4, m4 ), mustsplit );
    gap> splittingClassesWithOddCentralizerIndexSplit( m5, t, m5fust,
    >        GetFusionMap( 2m5, m5 ), mustsplit );
    gap> mustnotsplit:= [];;
    gap> notSplittingClassesOfSubgroupDoNotSplit( GetFusionMap( 2m4, m4 ),
    >        m4fust, mustnotsplit );
    #I  class 9 does not split (as in subgroup)
    #I  class 12 does not split (as in subgroup)
    #I  class 14 does not split (as in subgroup)
    #I  class 17 does not split (as in subgroup)
    #I  class 20 does not split (as in subgroup)
    #I  class 21 does not split (as in subgroup)
    #I  class 22 does not split (as in subgroup)
    #I  class 44 does not split (as in subgroup)
    #I  class 47 does not split (as in subgroup)
    #I  class 49 does not split (as in subgroup)
    #I  class 58 does not split (as in subgroup)
    #I  class 68 does not split (as in subgroup)
    #I  class 72 does not split (as in subgroup)
    #I  class 79 does not split (as in subgroup)
    #I  class 81 does not split (as in subgroup)
    #I  class 82 does not split (as in subgroup)
    #I  class 92 does not split (as in subgroup)
    gap> notSplittingClassesOfSubgroupDoNotSplit( GetFusionMap( 2m5, m5 ),
    >        m5fust, mustnotsplit );
    #I  class 13 does not split (as in subgroup)
    #I  class 18 does not split (as in subgroup)
    #I  class 45 does not split (as in subgroup)
    #I  class 48 does not split (as in subgroup)
    #I  class 69 does not split (as in subgroup)
    #I  class 80 does not split (as in subgroup)
    #I  class 93 does not split (as in subgroup)
\end{verbatim}

Since $\pi^{-1}( M_3 )$ is a direct product $\langle z \rangle \times M_3$,
the irreducible characters of $M_3$ can be taken as the characters $\chi'$
corresponding to the characters $\chi$
with the property $\chi(z) = - \chi(1)$.
For $M_4$ and $M_5$, we extract the relevant projective characters
from the character table of $2.F_4(2)$.

\begin{verbatim}
    gap> computeContributions( m3, t, m3fust, Irr( m3 ), 10^7,
    >        mustsplit, mustnotsplit );
    #I  class 2 splits (contribution criterion)
    #I  class 24 splits (3rd root of 2)
    #I  class 25 splits (3rd root of 2)
    #I  class 27 splits (3rd root of 2)
    #I  class 99 splits (3rd root of 27)
    #I  class 100 splits (3rd root of 27)
    #I  class 53 splits (5th root of 2)
    #I  class 8 splits (contribution criterion)
    #I  class 63 splits (3rd root of 8)
    #I  class 105 splits (5th root of 8)
    #I  class 15 splits (contribution criterion)
    #I  class 70 splits (3rd root of 15)
    #I  class 4 does not split (contribution criterion)
    #I  class 16 splits (contribution criterion)
    #I  class 78 splits (3rd root of 16)
    #I  class 59 splits (contribution criterion)
    #I  class 30 splits (contribution criterion)
    #I  class 102 splits (3rd root of 30)
    #I  class 3 splits (contribution criterion)
    #I  class 84 splits (contribution criterion)
    #I  class 26 splits (3rd root of 3)
    #I  class 28 splits (3rd root of 3)
    #I  class 54 splits (5th root of 3)
    #I  class 121 splits (5th root of 28)
    #I  class 32 splits (contribution criterion)
    #I  class 11 splits (contribution criterion)
    #I  class 67 splits (contribution criterion)
    #I  class 75 splits (contribution criterion)
    #I  class 117 splits (contribution criterion)
    #I  class 64 splits (3rd root of 11)
    #I  class 71 does not split (contribution criterion)
    gap> proj:= Filtered( Irr( 2m4 ), x -> x[1] <> x[2] );;
    gap> projmap:= ProjectionMap( GetFusionMap( 2m4, m4 ) );;
    gap> proj:= List( proj, x -> x{ projmap } );;
    gap> computeContributions( m4, t, m4fust, proj, 10^7,
    >        mustsplit, mustnotsplit );
    #I  class 118 does not split (contribution criterion)
    #I  class 31 splits (contribution criterion)
    #I  class 29 does not split (contribution criterion)
    gap> computeContributions( m5, t, m5fust, proj, 10^7,
    >        mustsplit, mustnotsplit );
    #I  class 119 does not split (contribution criterion)
    gap> mustsplit;
    [ 1, 2, 3, 5, 6, 7, 8, 11, 15, 16, 23, 24, 25, 26, 27, 28, 30, 31, 32, 33, 
      34, 51, 52, 53, 54, 55, 56, 59, 63, 64, 67, 70, 73, 75, 78, 83, 84, 85, 86, 
      87, 97, 98, 99, 100, 101, 102, 103, 104, 105, 106, 107, 108, 109, 110, 117, 
      120, 121, 122, 123, 124, 125, 126 ]
    gap> mustnotsplit;
    [ 4, 9, 12, 13, 14, 17, 18, 20, 21, 22, 29, 44, 45, 47, 48, 49, 58, 68, 69, 
      71, 72, 79, 80, 81, 82, 92, 93, 118, 119 ]
\end{verbatim}


\subsection{Element orders in $2.G$}\label{elementorders}

We claim that for $g \in G$,
the orders of the elements in $\pi^{-1}( g )$ are equal to the order of $g$
if this order is even.
(It is obvious that for $g$ of odd order, one element in $\pi^{-1}( g )$
has order $|g|$ and the other has order $2 |g|$.)

For that,
it is enough to show that the preimages of involutions in $G$ under $\pi$
are again involutions,
and this follows from the fact that $M_3$ contains elements from all
involution classes in $G$,
since $\pi^{-1}( M_3 )$ is a direct product and hence involutions in $M_3$
lift to involutions in the preimage.

\begin{verbatim}
    gap> invol:= Positions( orders, 2 );
    [ 2, 3, 4 ]
    gap> Difference( invol, m3fust );
    [  ]
\end{verbatim}

\subsection{The class fusion from $\pi^{-1}( M_3 )$}\label{sectinitialfusion}

In order to get more information about the classes of $2.G$,
we apply a different strategy.
Currently $12$ classes of $G$ are left which contain elements of $M_3$
and for which we do not know whether they split.


\begin{verbatim}
    gap> open:= Difference( m3fust, Union( mustsplit, mustnotsplit ) );
    [ 19, 35, 36, 39, 43, 46, 50, 91, 94, 111, 112, 114 ]
\end{verbatim}

In fact, we need not consider all $2^{12}$ combinations from this set,
because the third power map connects some of the classes.

\begin{verbatim}
    gap> orders{ open };
    [ 4, 8, 8, 8, 8, 8, 8, 16, 16, 24, 24, 24 ]
    gap> PowerMap( t, 3 ){ [ 111, 112, 114 ] };
    [ 35, 36, 39 ]
    gap> poss:= Filtered( Combinations( open ),
    >               x -> ( not 35 in x or 111 in x ) and
    >                    ( not 36 in x or 112 in x ) and
    >                    ( not 39 in x or 114 in x ) );;
    gap> Length( poss );
    1728
\end{verbatim}

For the $1728$ cases, we create a preliminary character table head
of $2.G$, and compute possible class fusions from $\pi^{-1}( M_3 )$
into this table.
Again, the idea is to look at the norms of induced characters,
but now we will use several induced characters simultaneously.

We need a few more {\GAP} functions for this step.
The first one creates the character table head object that corresponds
to the proposed set of splitting classes.
We will use the information about element orders
from Section~\ref{elementorders} by setting the arguments \verb|invmustlift|
and \verb|invmaylift| to empty lists.

\begin{verbatim}
    gap> tableHead:= function( t, tosplit, invmustlift, invmaylift )
    >      local tcents, orders, splcentralizers, spl, splorders, mustlift,
    >            maylift, i, pow, ord;
    > 
    >      tcents:= SizesCentralizers( t );
    >      orders:= OrdersClassRepresentatives( t );
    >      splcentralizers:= [];
    >      spl:= [];
    >      splorders:= [];
    >      mustlift:= ShallowCopy( invmustlift );
    >      maylift:= ShallowCopy( invmaylift );
    > 
    >      if invmaylift <> [] or invmustlift <> [] then
    >        for i in [ 2 .. NrConjugacyClasses( t ) ] do
    >          if orders[i] mod 2 = 0 then
    >            pow:= PowerMap( t, orders[i] / 2 )[i];
    >            if pow in invmustlift then
    >              Add( mustlift, i );
    >            elif pow in invmaylift then
    >              Add( maylift, i );
    >            fi;
    >          fi;
    >        od;
    >      fi;
    > 
    >      for i in [ 1 .. NrConjugacyClasses( t ) ] do
    >        ord:= orders[i];
    >        if i in tosplit then
    >          Append( spl, [ i, i ] );
    >          Append( splcentralizers, tcents[i] * [ 2, 2 ] );
    >          if orders[i] mod 2 = 1 then
    >            Append( splorders, [ ord, 2 * ord ] );
    >          elif i in mustlift then
    >            Append( splorders, [ 2 * ord, 2 * ord ] );
    >          elif i in maylift then
    >            Append( splorders, [ [ ord, 2 * ord ], [ ord, 2 * ord ] ] );
    >          else
    >            Append( splorders, [ ord, ord ] );
    >          fi;
    >        else
    >          Add( spl, i );
    >          Add( splcentralizers, tcents[i] );
    >          if i in mustlift then
    >            Add( splorders, 2 * ord );
    >          elif i in maylift then
    >            Add( splorders, [ ord, 2 * ord ] );
    >          else
    >            Add( splorders, ord );
    >          fi;
    >        fi;
    >      od;
    > 
    >      return ConvertToCharacterTableNC( rec(
    >                 UnderlyingCharacteristic:= 0,
    >                 OrdersClassRepresentatives:= splorders,
    >                 SizesCentralizers:= splcentralizers,
    >                 Size:= splcentralizers[1],
    >                 ComputedClassFusions:= [ rec( name:= Identifier( t ),
    >                     map:= spl ) ],
    >               ) );
    >  end;;
\end{verbatim}

The next function creates an approximation of the class fusion from
the character table \verb|2s| of a subgroup $\pi^{-1}(U)$
into the character table \verb|2t| of $\pi^{-1}(G)$,
using that the composition of the class fusions \verb|sfust| from $U$ to $G$
and \verb|2sfuss| from $\pi^{-1}(U)$ to $U$ must be equal to the composition
of the class fusion \verb|2tfust| from $\pi^{-1}(G)$ to $G$ and the desired
class fusion from $\pi^{-1}(U)$ to $\pi^{-1}(G)$.

Moreover, the function takes a list \verb|defined| of class positions in $G$
such that the preimages in $2.G$ have already been defined;
the fusion from one pair of classes in \verb|2s| to a given pair of classes
in \verb|2t| may be chosen if the numbers in the latter pair map to a class
position in \verb|t| that does not occur in \verb|defined|.


It may happen that the function returns \verb|fail|,
because the splitting of classes prescribed by \verb|2tfust|
is not compatible with the embedding of \verb|2s| in \verb|2t|.
In all other cases, the result is an approximation of the class fusion
from \verb|2s| to \verb|2t| in the sense that the entry at position $i$ is either
a class position in \verb|2t| (the position of the unique possible image class
of the $i$-th class of \verb|2s|) or a list of such class positions.

\begin{verbatim}
    gap> initialFusion:= function( 2s, 2t, 2sfuss, 2tfust, sfust, defined )
    >      local fus, comp, pre, imgs;
    > 
    >      # Use element orders and centralizer orders.
    >      fus:= InitFusion( 2s, 2t );
    > 
    >      # Use the commutative diagram.
    >      comp:= CompositionMaps( InverseMap( 2tfust ),
    >                              CompositionMaps( sfust, 2sfuss ) );
    >      if MeetMaps( fus, comp ) <> true then
    >        return fail;
    >      fi;
    > 
    >      # Define classes that are not yet defined.
    >      defined:= ShallowCopy( defined );
    >      for pre in InverseMap( 2sfuss ) do
    >        if IsList( pre ) then
    >          imgs:= fus{ pre };
    >          if imgs[1] = imgs[2] and IsList( imgs[1] ) then
    >            if Intersection( defined, 2tfust{ imgs[1] } ) = [] then
    >              # The classes in preimage and image split, and we may choose.
    >              fus[ pre[1] ]:= imgs[1][1];
    >              fus[ pre[2] ]:= imgs[1][2];
    >              UniteSet( defined, 2tfust{ imgs[1] } );
    >            fi;
    >          fi;
    >        elif IsList( fus[ pre ] ) then
    >          # The class splits in the image but not in the preimage,
    >          # we should have noticed this earlier.
    >          return fail;
    >        fi;
    >      od;
    > 
    >      return fus;
    >    end;;
\end{verbatim}

The next function tries to improve an approximation of a class fusion
\verb|2sfus2t| between the character tables \verb|2s| and \verb|2t|
by computing the class functions of \verb|2t| that can be induced from a
character \verb|chi| of \verb|2s| via some choice of images compatible with \verb|2sfus2t|,
and then computing the subset of class functions which have integral norms.
Whenever this subset allows us to deduce that some possible images
in \verb|2sfus2t| cannot occur then \verb|2sfus2t| gets improved in place.

If no induced class function can have integral norm, the function returns
\verb|false|, otherwise \verb|true| is returned.

\begin{verbatim}
    gap> useInducedClassFunction:= function( 2s, 2t, chi, 2sfuss, 2sfus2t )
    >      local localfus, unknown, i, swaps, inv, pair, poss, choices, choice,
    >            map, ind, para, new;
    > 
    >      # Remove indet. in places where the character is zero.
    >      localfus:= ShallowCopy( 2sfus2t );
    >      unknown:= [];
    >      for i in [ 1 .. Length( 2sfus2t ) ] do
    >        if IsList( 2sfus2t[i] ) then
    >          if chi[i] = 0 then
    >            localfus[i]:= localfus[i][1];
    >          else
    >            Add( unknown, i );
    >          fi;
    >        fi;
    >      od;
    > 
    >      # Collect the possible swaps.
    >      swaps:= [];
    >      inv:= InverseMap( 2sfuss );
    >      for i in [ 1 .. Length( localfus ) ] do
    >        if IsList( localfus[i] ) then
    >          pair:= inv[ 2sfuss[i] ];
    >          Add( swaps, ( pair[1], pair[2] ) );
    >          localfus{ pair }:= localfus[i];
    >        fi;
    >      od;
    > 
    >      # Try all possibilities (hopefully not too many).
    >      poss:= [];
    >      if IsEmpty( swaps ) then
    >        choices:= [ [] ];
    >      else
    >        choices:= IteratorOfCombinations( swaps );
    >      fi;
    >      for choice in choices do
    >        map:= Permuted( localfus, Product( choice, () ) );
    >        ind:= InducedClassFunctionsByFusionMap( 2s, 2t, [ chi ], map )[1];
    >        if IsInt( ScalarProduct( 2t, ind, ind ) ) then
    >          Add( poss, map );
    >        fi;
    >      od;
    > 
    >      if poss = [] then
    >        return false;
    >      fi;
    > 
    >      para:= Parametrized( poss );
    >      new:= Filtered( unknown, i -> IsInt( para[i] ) );
    >      if new <> [] then
    >        2sfus2t{ new }:= para{ new };
    >      fi;
    > 
    >      return true;
    >    end;;
\end{verbatim}

Now we check the possible splittings that can be distinguished
by characters of $\pi^{-1}( M_3 )$.

\begin{verbatim}
    gap> good:= [];;
    gap> ker:= ClassPositionsOfKernel( GetFusionMap( 2m3, m3 ) );;
    gap> testcharsm3:= Filtered( Irr( 2m3 ),
    >        chi -> not IsSubset( ClassPositionsOfKernel( chi ), ker ) );;
    gap> runOneTest:= function( s, 2s, t, 2t, sfust, testchars, defined )
    >      local fus, pos, l, chi;
    >      fus:= initialFusion( 2s, 2t, GetFusionMap( 2s, s ),
    >                           GetFusionMap( 2t, t ), sfust, defined );
    >      # Process the irreducible characters,
    >      # ordered by increasing indeterminateness.
    >      pos:= PositionsProperty( fus, IsList );
    >      testchars:= ShallowCopy( testchars );
    >      l:= - List( testchars, x -> Number( pos, i -> x[i] = 0 ) );
    >      SortParallel( l, testchars );
    >      for chi in testchars do
    >        if useInducedClassFunction( 2s, 2t, chi, GetFusionMap( 2s, s ),
    >                                    fus ) = false then
    >          # This splitting is not possible.
    >          return fail;
    >        fi;
    >      od;
    >      return fus;
    >    end;;
    gap> defined:= [];;
    gap> for choice in poss do
    >      2t:= tableHead( t, Union( mustsplit, choice ), [], [] );
    >      fus:= runOneTest( m3, 2m3, t, 2t, m3fust, testcharsm3, defined );
    >      if fus <> fail then
    >        Add( good, choice );
    >      fi;
    >    od;
    gap> Length( good );
    1
\end{verbatim}


We get a unique solution for the splitting.
Thus we update our lists,
create the table head corresponding to the currently known splitting,
and recompute the class fusion from $\pi^{-1}( M_3 )$ to $2.G$,
which is also uniquely determined.


\begin{verbatim}
    gap> choice:= good[1];
    [ 19, 36, 39, 43, 46, 50, 91, 94, 111, 112, 114 ]
    gap> UniteSet( mustsplit, choice );
    gap> oddRootsOfSplittingClassesSplit( t, mustsplit );
    #I  class 74 splits (3rd root of 19)
    #I  class 76 splits (3rd root of 19)
    #I  class 77 splits (3rd root of 19)
    gap> UniteSet( mustnotsplit, Difference( open, choice ) );
    gap> Difference( [ 1 .. Length( orders ) ],
    >                Union( mustsplit, mustnotsplit ) );
    [ 10, 37, 38, 40, 41, 42, 57, 60, 61, 62, 65, 66, 88, 89, 90, 95, 96, 113, 
      115, 116 ]
    gap> 2t:= tableHead( t, mustsplit, [], [] );;
    gap> NrConjugacyClasses( 2t );
    202
    gap> 2m3fus2t:= runOneTest( m3, 2m3, t, 2t, m3fust, testcharsm3, defined );
    [ 1, 3, 5, 5, 7, 8, 10, 12, 16, 18, 14, 18, 23, 23, 22, 25, 26, 31, 23, 29, 
      31, 32, 33, 34, 36, 44, 40, 38, 42, 47, 40, 44, 46, 49, 51, 55, 53, 71, 71, 
      57, 62, 58, 63, 67, 71, 68, 75, 76, 80, 78, 82, 84, 84, 99, 92, 91, 103, 
      107, 111, 107, 111, 99, 97, 109, 111, 110, 121, 115, 125, 126, 127, 131, 
      129, 135, 133, 144, 145, 140, 140, 148, 150, 156, 158, 166, 166, 170, 168, 
      181, 176, 182, 178, 189, 185, 193, 191, 2, 4, 6, 6, 7, 9, 11, 13, 16, 19, 
      15, 19, 24, 24, 22, 26, 25, 31, 24, 30, 31, 32, 33, 35, 37, 45, 41, 39, 43, 
      48, 41, 45, 46, 50, 52, 56, 54, 72, 72, 57, 63, 59, 62, 68, 72, 67, 75, 77, 
      81, 79, 83, 85, 85, 100, 93, 91, 104, 108, 112, 108, 112, 100, 98, 109, 
      112, 110, 122, 116, 125, 126, 128, 132, 130, 136, 134, 145, 144, 141, 141, 
      149, 151, 157, 159, 167, 167, 171, 169, 182, 177, 181, 179, 190, 186, 194, 
      192 ]
    gap> defined:= Set( m3fust );;
\end{verbatim}

For later use, we compute the characters of $2.G$ that arise by
induction from $\pi^{-1}( M_3 )$.

\begin{verbatim}
    gap> ind2m3:= InducedClassFunctionsByFusionMap( 2m3, 2t, testcharsm3,
    >                 2m3fus2t );;
\end{verbatim}

\subsection{Subgroups of the type $Fi_{22}$ in $G$}\label{fi22subs}


There are three possible class fusion of $Fi_{22}$ in $G$,
up to symmetries of the character table of $Fi_{22}$,
and these fusions form one orbit under the action of $\alpha$.

\begin{verbatim}
    gap> s:= CharacterTable( "Fi22" );;
    gap> fus:= PossibleClassFusions( s, t );;
    gap> rep:= RepresentativesFusions( s, fus, Group( () ) );;
    gap> Length( rep );
    3
    gap> List( rep, map -> Intersection( map, oneorbit ) );
    [ [ 11 ], [ 12 ], [ 13 ] ]
\end{verbatim}

By our assumption that $G$ contains $Fi_{22}$ type subgroups,
this means that all three possibilities are really class fusions
of $Fi_{22}$ type subgroups in $G$.
We fix three such subgroups with these fusions,
in one orbit under the action of $\alpha$,
and call the subgroups $M_7$, $M_8$, and $M_9$.
(Note that the three classes of maximal subgroups of the type $Fi_{22}$
in $G$ claimed in the list in~\cite[pp.~191]{CCN85}
appear in the positions $7$ to $9$.)

\begin{verbatim}
    gap> m7:= s;;  m7fust:= rep[1];;
    gap> m8:= s;;  m8fust:= rep[2];;
    gap> m9:= s;;  m9fust:= rep[3];;
\end{verbatim}

Concerning the structures of the preimages $\pi^{-1}( M_i )$,
for $7 \leq i \leq 9$,
the same arguments as for $F_4(2)$ yield that at least one of them
is a direct product $2 \times Fi_{22}$,
and that the other two have the same structure
--either also $2 \times Fi_{22}$ or double covers $2.Fi_{22}$.

\subsection{The class fusion from $\pi^{-1}( M_7 )$}

The first question about $M_7$ is whether $\pi^{-1}( M_7 )$ has the type
$2.Fi_{22}$ or $2 \times Fi_{22}$.
The former case is ruled out by the fact that no fusion exists,
according to the current knowledge about $2.G$.
(The reason is that class $11$ of $G$,
which splits as a consequence of the embedding of $M_3$,
contains elements in a class of $M_7$ that does not split.)

\begin{verbatim}
    gap> open:= Difference( m7fust, Union( mustsplit, mustnotsplit ) );
    [  ]
    gap> 2s:= CharacterTable( "2.Fi22" );;
    gap> initialFusion( 2s, 2t, GetFusionMap( 2s, m7 ),
    >        GetFusionMap( 2t, t ), m7fust, defined );
    fail
\end{verbatim}

Thus $\pi^{-1}( M_7 )$ has the type $2 \times Fi_{22}$.

\begin{verbatim}
    gap> 2m7:= CharacterTable( "Cyclic", 2 ) * m7;;
\end{verbatim}

We are now in a better situation than in the case of the subgroup $M_3$,
since we can use that the induced characters computed above
must restrict to characters of $\pi^{-1}( M_7 )$.
This criterion is implemented by the {\GAP} function
\verb|FusionsAllowedByRestrictions|, which requires some parameters as its
last argument; we use the following settings.

\begin{verbatim}
    gap> parametersFABR:= rec( maxlen:= 10, minamb:= 1, maxamb:= 10^6,
    >        quick:= false, contained:= ContainedPossibleCharacters );;
\end{verbatim}

The class fusion in $2.G$ is not uniquely determined,
we get two possibilities and also two possible lists of
characters induced from $\pi^{-1}( M_7 )$.

\begin{verbatim}
    gap> ker:= ClassPositionsOfKernel( GetFusionMap( 2m7, m7 ) );;
    gap> testcharsm7:= Filtered( Irr( 2m7 ),
    >        chi -> not IsSubset( ClassPositionsOfKernel( chi ), ker ) );;
    gap> fus:= initialFusion( 2m7, 2t, GetFusionMap( 2m7, m7 ),
    >              GetFusionMap( 2t, t ), m7fust, defined );;
    gap> possfus:= FusionsAllowedByRestrictions( 2m7, 2t, testcharsm7,
    >                  ind2m3, fus, parametersFABR );
    [ [ 1, 4, 5, 7, 8, 12, 10, 12, 18, 24, 26, 29, 31, 34, 37, 43, 47, 40, 39, 
          46, 52, 49, 44, 50, 51, 55, 63, 62, 68, 75, 78, 80, 80, 83, 84, 86, 88, 
          103, 99, 99, 108, 117, 119, 113, 122, 109, 112, 115, 127, 127, 132, 
          135, 140, 140, 152, 154, 157, 158, 167, 170, 172, 174, 181, 182, 194, 
          2, 3, 6, 7, 9, 13, 11, 13, 19, 23, 25, 30, 31, 35, 36, 42, 48, 41, 38, 
          46, 51, 50, 45, 49, 52, 56, 62, 63, 67, 75, 79, 81, 81, 82, 85, 87, 89, 
          104, 100, 100, 107, 118, 120, 114, 121, 109, 111, 116, 128, 128, 131, 
          136, 141, 141, 153, 155, 156, 159, 166, 171, 173, 175, 182, 181, 193 ], 
      [ 1, 4, 5, 7, 8, 12, 10, 12, 18, 24, 26, 29, 31, 34, 37, 43, 47, 40, 39, 
          46, 52, 49, 44, 50, 51, 55, 62, 63, 68, 75, 78, 80, 80, 83, 84, 86, 88, 
          103, 99, 99, 108, 117, 119, 113, 122, 109, 112, 115, 127, 127, 132, 
          135, 140, 140, 152, 154, 157, 158, 167, 170, 172, 174, 182, 181, 194, 
          2, 3, 6, 7, 9, 13, 11, 13, 19, 23, 25, 30, 31, 35, 36, 42, 48, 41, 38, 
          46, 51, 50, 45, 49, 52, 56, 63, 62, 67, 75, 79, 81, 81, 82, 85, 87, 89, 
          104, 100, 100, 107, 118, 120, 114, 121, 109, 111, 116, 128, 128, 131, 
          136, 141, 141, 153, 155, 156, 159, 166, 171, 173, 175, 181, 182, 193 ] ]
    gap> poss2m7fus2t:= possfus;;
    gap> UniteSet( defined, Set( m7fust ) );
    gap> possind2m7:= List( poss2m7fus2t,
    >        map -> Set( InducedClassFunctionsByFusionMap( 2m7, 2t,
    >                        testcharsm7, map ) ) );;
    gap> List( possind2m7, Length );
    [ 63, 63 ]
    gap> Length( Intersection( possind2m7 ) );
    39
\end{verbatim}

\subsection{The class fusion from $\pi^{-1}( M_8 )$}

Now we turn to the subgroup $M_8$.
Again, we have to decide whether its type is
$2.Fi_{22}$ or $2 \times Fi_{22}$.

This time, the latter case is ruled out by the fact that no fusion exists,
according to the current knowledge about $2.G$.
For that, we first observe that the splitting of three classes
in the image of the class fusion from $M_8$ to $G$ is not yet
decided.
We check in each of the eight possible situations that no class fusion
of $2 \times Fi_{22}$ in $2.G$ exists that is compatible with both the
class fusion of $M_8$ in $G$
and the characters of $2.G$ that are induced from $\pi^{-1}( M_3 )$.

\begin{verbatim}
    gap> 2s:= 2m7;
    CharacterTable( "C2xFi22" )
    gap> open:= Difference( m8fust, Union( mustsplit, mustnotsplit ) );
    [ 40, 65, 115 ]
    gap> good:= [];;
    gap> for choice in Combinations( open ) do
    >      2t:= tableHead( t, Union( mustsplit, choice ), [], [] );
    >      fus:= initialFusion( 2s, 2t, GetFusionMap( 2s, m8 ),
    >                GetFusionMap( 2t, t ), m8fust, defined );;
    >      if FusionsAllowedByRestrictions( 2s, 2t, testcharsm7, ind2m3, fus,
    >             parametersFABR ) <> [] then
    >        Add( good, choice );
    >      fi;
    >    od;
    gap> good;
    [  ]
\end{verbatim}

Thus $\pi^{-1}( M_8 )$ has the type $2.Fi_{22}$.

\begin{verbatim}
    gap> 2m8:= CharacterTable( "2.Fi22" );;
\end{verbatim}

This yields information about the splitting of the three classes.

\begin{verbatim}
    gap> notSplittingClassesOfSubgroupDoNotSplit( GetFusionMap( 2m8, m8 ),
    > m8fust, mustnotsplit );
    #I  class 40 does not split (as in subgroup)
    #I  class 65 does not split (as in subgroup)
    #I  class 115 does not split (as in subgroup)
    gap> 2t:= tableHead( t, mustsplit, [], [] );;
\end{verbatim}

The class fusion in $2.G$ is uniquely determined.

\begin{verbatim}
    gap> ker:= ClassPositionsOfKernel( GetFusionMap( 2m8, m8 ) );;
    gap> testcharsm8:= Filtered( Irr( 2m8 ),
    >        chi -> not IsSubset( ClassPositionsOfKernel( chi ), ker ) );;
    gap> fus:= initialFusion( 2m8, 2t, GetFusionMap( 2m8, m8 ),
    >        GetFusionMap( 2t, t ), m8fust, defined );;
    gap> ind:= Concatenation( ind2m3, Intersection( possind2m7 ) );;
    gap> possfus:= FusionsAllowedByRestrictions( 2m8, 2t, testcharsm8, ind,
    >        fus, parametersFABR );
    [ [ 1, 2, 3, 4, 5, 6, 7, 8, 9, 12, 13, 10, 11, 12, 13, 20, 23, 24, 27, 29, 
          30, 31, 34, 35, 36, 37, 42, 43, 47, 48, 40, 41, 38, 39, 46, 51, 52, 49, 
          50, 44, 45, 49, 50, 51, 52, 55, 56, 64, 64, 69, 75, 75, 78, 79, 80, 81, 
          80, 81, 82, 83, 84, 85, 86, 87, 88, 89, 105, 101, 101, 101, 107, 108, 
          116, 115, 120, 119, 113, 114, 123, 109, 111, 112, 117, 118, 127, 128, 
          127, 128, 131, 132, 135, 136, 142, 142, 153, 152, 155, 154, 156, 157, 
          158, 159, 166, 167, 170, 171, 173, 172, 175, 174, 183, 183, 193, 194 ] ]
    gap> 2m8fus2t:= possfus[1];;
    gap> UniteSet( defined, m8fust );
    gap> ind2m8:= InducedClassFunctionsByFusionMap( 2m8, 2t, testcharsm8,
    >                 2m8fus2t );;
\end{verbatim}

\subsection{The class fusion from $\pi^{-1}( M_9 )$}

By the action of the outer automorphism $\beta$,
we know that $\pi^{-1}( M_9 )$ has the type $2.Fi_{22}$.
The images under $\beta$ of the not splitting classes found above
do not split,
and the class fusion from $\pi^{-1}( M_9 )$ is determined analogously
to that from $\pi^{-1}( M_8 )$.

\begin{verbatim}
    gap> Filtered( orbsbeta, l -> Intersection( l, [ 40, 65, 115 ] ) <> [] );
    [ [ 40, 41 ], [ 65, 66 ], [ 115, 116 ] ]
    gap> UniteSet( mustnotsplit, [ 41, 66, 116 ] );
    gap> open:= Difference( m9fust, Union( mustsplit, mustnotsplit ) );
    [  ]
    gap> 2m9:= 2m8;;
    gap> testcharsm9:= testcharsm8;;
    gap> fus:= initialFusion( 2m9, 2t, GetFusionMap( 2m9, m9 ),
    >        GetFusionMap( 2t, t ), m9fust, defined );;
    gap> ind:= Concatenation( ind2m3, Intersection( possind2m7 ), ind2m8 );;
    gap> possfus:= FusionsAllowedByRestrictions( 2m9, 2t, testcharsm9, ind,
    >        fus, parametersFABR );
    [ [ 1, 2, 3, 4, 5, 6, 7, 8, 9, 12, 13, 10, 11, 12, 13, 21, 23, 24, 28, 29, 
          30, 31, 34, 35, 36, 37, 42, 43, 47, 48, 40, 41, 38, 39, 46, 51, 52, 49, 
          50, 44, 45, 49, 50, 51, 52, 55, 56, 65, 65, 70, 75, 75, 78, 79, 80, 81, 
          80, 81, 82, 83, 84, 85, 86, 87, 88, 89, 106, 102, 102, 102, 107, 108, 
          116, 115, 118, 117, 113, 114, 124, 109, 111, 112, 119, 120, 127, 128, 
          127, 128, 131, 132, 135, 136, 143, 143, 153, 152, 155, 154, 156, 157, 
          158, 159, 166, 167, 170, 171, 173, 172, 175, 174, 184, 184, 193, 194 ] ]
    gap> 2m9fus2t:= possfus[1];;
    gap> UniteSet( defined, m9fust );
    gap> ind2m9:= InducedClassFunctionsByFusionMap( 2m9, 2t, testcharsm9,
    >                 2m9fus2t );;
\end{verbatim}

\subsection{The class fusion from $\pi^{-1}( M_4 )$}

The subgroup $M_4$ of $G$ contains elements in one class for which
we do not know yet whether it splits.
As above, we use norms of induced characters in order to determine
the splitting.
In order to speed up \verb|runOneTest|, we do {\bf not} enter the
information about the classes of $G$ whose preimages are already defined.
Thus we cannot use the class fusion that is returned by \verb|runOneTest|,
it may be not compatible with the information which we have already
collected.

\begin{verbatim}
    gap> open:= Difference( m4fust, Union( mustsplit, mustnotsplit ) );
    [ 95 ]
    gap> ker:= ClassPositionsOfKernel( GetFusionMap( 2m4, m4 ) );;
    gap> testcharsm4:= Filtered( Irr( 2m4 ),
    >        chi -> not IsSubset( ClassPositionsOfKernel( chi ), ker ) );;
    gap> good:= [];;
    gap> for choice in Combinations( open ) do
    >      2t:= tableHead( t, Union( mustsplit, choice ), [], [] );
    >      fus:= runOneTest( m4, 2m4, t, 2t, m4fust, testcharsm4, [] );
    >      if fus <> fail then
    >        Add( good, [ choice, 2t ] );
    >      fi;
    >    od;
    gap> List( good, l -> l[1] );
    [ [  ] ]
    gap> UniteSet( mustnotsplit, open );
    gap> 2t:= good[1][2];;
    gap> 2tfust:= GetFusionMap( 2t, t );;
    gap> fus:= initialFusion( 2m4, 2t, GetFusionMap( 2m4, m4 ),
    >                         2tfust, m4fust, defined );;
    gap> ind:= Concatenation( ind2m3, Intersection( possind2m7 ), ind2m8,
    >                         ind2m9 );;
    gap> possfus:= FusionsAllowedByRestrictions( 2m4, 2t, testcharsm4, ind,
    >        fus, parametersFABR );
    [ [ 1, 2, 3, 4, 6, 5, 5, 6, 7, 7, 8, 9, 10, 11, 12, 13, 16, 20, 14, 15, 20, 
          20, 24, 23, 24, 23, 22, 27, 27, 27, 31, 31, 23, 24, 29, 30, 31, 32, 33, 
          34, 35, 36, 37, 45, 44, 41, 40, 38, 39, 42, 43, 48, 47, 40, 41, 44, 45, 
          46, 46, 50, 49, 52, 51, 55, 56, 53, 54, 73, 73, 73, 57, 57, 64, 64, 58, 
          59, 64, 64, 69, 73, 69, 75, 76, 77, 80, 81, 78, 79, 82, 83, 85, 84, 84, 
          85, 101, 101, 92, 93, 91, 105, 108, 107, 112, 111, 108, 107, 112, 111, 
          101, 101, 97, 98, 109, 109, 111, 112, 110, 123, 117, 118, 125, 126, 
          127, 128, 131, 132, 130, 129, 135, 136, 133, 134, 146, 146, 146, 146, 
          142, 142, 148, 149, 150, 151, 156, 157, 159, 158, 167, 166, 167, 166, 
          170, 171, 168, 169, 183, 183, 177, 176, 183, 183, 178, 179, 189, 190, 
          187, 187, 194, 193, 192, 191 ] ]
    gap> 2m4fus2t:= possfus[1];;
    gap> UniteSet( defined, Set( m4fust ) );
    gap> ind2m4:= InducedClassFunctionsByFusionMap( 2m4, 2t, testcharsm4,
    >                 2m4fus2t );;
\end{verbatim}

\subsection{The class fusion from $\pi^{-1}( M_5 )$}

The case $\pi^{-1}( M_5 )$ is analogous.
The classes in the image of the class fusion from $M_5$ to $G$
that are not known to split are the images under the outer automorphisms
$\beta$ of the classes that were considered for $M_4$.
Thus we know that also these classes do not split.

\begin{verbatim}
    gap> open:= Difference( m5fust, Union( mustsplit, mustnotsplit ) );
    [ 96 ]
    gap> UniteSet( mustnotsplit, open );
\end{verbatim}

The class fusion is determined in the same way as for $M_4$.

\begin{verbatim}
    gap> ker:= ClassPositionsOfKernel( GetFusionMap( 2m5, m5 ) );;
    gap> testcharsm5:= Filtered( Irr( 2m5 ),
    >        chi -> not IsSubset( ClassPositionsOfKernel( chi ), ker ) );;
    gap> fus:= initialFusion( 2m5, 2t, GetFusionMap( 2m5, m5 ),
    >              GetFusionMap( 2t, t ), m5fust, defined );;
    gap> ind:= Concatenation( ind2m3, Intersection( possind2m7 ), ind2m8,
    >              ind2m9, ind2m4 );;
    gap> possfus:= FusionsAllowedByRestrictions( 2m5, 2t, testcharsm5,
    >                  ind, fus, parametersFABR );
    [ [ 1, 2, 3, 4, 6, 5, 5, 6, 7, 7, 8, 9, 10, 11, 12, 13, 16, 21, 14, 15, 21, 
          21, 24, 23, 24, 23, 22, 28, 28, 28, 31, 31, 23, 24, 29, 30, 31, 32, 33, 
          34, 35, 36, 37, 45, 44, 41, 40, 38, 39, 42, 43, 48, 47, 40, 41, 44, 45, 
          46, 46, 50, 49, 52, 51, 55, 56, 53, 54, 74, 74, 74, 57, 57, 65, 65, 58, 
          59, 65, 65, 70, 74, 70, 75, 76, 77, 80, 81, 78, 79, 82, 83, 85, 84, 84, 
          85, 102, 102, 92, 93, 91, 106, 108, 107, 112, 111, 108, 107, 112, 111, 
          102, 102, 97, 98, 109, 109, 111, 112, 110, 124, 119, 120, 125, 126, 
          127, 128, 131, 132, 130, 129, 135, 136, 133, 134, 147, 147, 147, 147, 
          143, 143, 148, 149, 150, 151, 156, 157, 159, 158, 167, 166, 167, 166, 
          170, 171, 168, 169, 184, 184, 177, 176, 184, 184, 178, 179, 189, 190, 
          188, 188, 194, 193, 192, 191 ] ]
    gap> 2m5fus2t:= possfus[1];;
    gap> UniteSet( defined, Set( m5fust ) );
    gap> ind2m5:= InducedClassFunctionsByFusionMap( 2m5, 2t, testcharsm5,
    >                 2m5fus2t );;
\end{verbatim}

\subsection{Subgroups of the type $3 \times U_6(2)$ in $G$}\label{3u62subs}

Let $U_{12}$ be a subgroup of the type $3 \times U_6(2)$ in $G$.
The preimage under $\pi$ may have the structure $6 \times U_6(2)$
or $3 \times 2.U_6(2)$,
where $2.U_6(2)$ denotes the double cover of $U_6(2)$.
The latter possibility must occur, since the former does not admit
a class fusion into $2.G$.

\begin{verbatim}
    gap> s:= CharacterTable( "U6(2)" );;
    gap> InitFusion( CharacterTable( "C6" ) * s, 2t );
    fail
\end{verbatim}

(We construct first the character table of $\pi^{-1}( U_{12} )$
and then that of $U_{12}$ as the table of the factor group by the central
subgroup of order two;
this guarantees that the factor fusion between the two tables
is automatically available.)

\begin{verbatim}
    gap> 2s:= CharacterTable( "2.U6(2)" );;
    gap> 2u12:= CharacterTable( "C3" ) * 2s;;
    gap> orders2u12:= OrdersClassRepresentatives( 2u12 );;
    gap> inv:= First( ClassPositionsOfCentre( 2u12 ), 
    >                 i -> orders2u12[i] = 2 );;
    gap> ker:= [ 1, inv ];;
    gap> u12:= 2u12 / ker;;
    gap> testcharsu12:= Filtered( Irr( 2u12 ),
    >        chi -> not IsSubset( ClassPositionsOfKernel( chi ), ker ) );;
\end{verbatim}

The class fusion of $U_{12}$ in $G$ is not uniquely
determined up to symmetries of $U_{12}$;
we get two candidates for the fusion.
(The fusion \textbf{is} unique up to symmetries of $U_{12}$ and $G$,
but our choices of the fusions for the other subgroups may have broken
some symmetries of $G$.)


\begin{verbatim}
    gap> fus:= PossibleClassFusions( u12, t );;
    gap> rep:= RepresentativesFusions( u12, fus, Group( () ) );;
    gap> Length( rep );
    2
    gap> possu12fust:= rep;;
\end{verbatim}

This subgroup yields new information about the classes of $2.G$:
Seven out of the $12$ classes of $G$ for which we do not yet know
whether they split can now be shown not to split.

\begin{verbatim}
    gap> List( possu12fust, map -> Difference( map,
    >                         Union( mustsplit, mustnotsplit ) ) );
    [ [ 10, 37, 57, 60, 61, 62, 113 ], [ 10, 37, 57, 60, 61, 62, 113 ] ]
    gap> possnotsplit:= List( [ 1, 2 ], i -> ShallowCopy( mustnotsplit ) );;
    gap> for i in [ 1, 2 ] do
    >      notSplittingClassesOfSubgroupDoNotSplit(
    >          GetFusionMap( 2u12, u12 ), rep[i], possnotsplit[i] );
    >    od;
    #I  class 10 does not split (as in subgroup)
    #I  class 37 does not split (as in subgroup)
    #I  class 57 does not split (as in subgroup)
    #I  class 60 does not split (as in subgroup)
    #I  class 61 does not split (as in subgroup)
    #I  class 62 does not split (as in subgroup)
    #I  class 113 does not split (as in subgroup)
    #I  class 10 does not split (as in subgroup)
    #I  class 37 does not split (as in subgroup)
    #I  class 57 does not split (as in subgroup)
    #I  class 60 does not split (as in subgroup)
    #I  class 61 does not split (as in subgroup)
    #I  class 62 does not split (as in subgroup)
    #I  class 113 does not split (as in subgroup)
    gap> Set( possnotsplit );
    [ [ 4, 9, 10, 12, 13, 14, 17, 18, 20, 21, 22, 29, 35, 37, 40, 41, 44, 45, 47, 
          48, 49, 57, 58, 60, 61, 62, 65, 66, 68, 69, 71, 72, 79, 80, 81, 82, 92, 
          93, 95, 96, 113, 115, 116, 118, 119 ] ]
    gap> mustnotsplit:= possnotsplit[1];;
\end{verbatim}

We defer the computation of the possible class fusions
of $\pi^{-1}( U_{12} )$ in $2.G$
until we know characters induced from another subgroup,
in Section~\ref{2o10m2subs}.

\subsection{Subgroups of the type $O^-_{10}(2)$ in $G$}\label{2o10m2subs}

Let $M_{10}$ be a subgroup of type $O^-_{10}(2)$ in $G$.
The class fusion of $M_{10}$ in $G$ is not uniquely
determined up to symmetries of $M_{10}$;
we get two candidates for the fusion.
(The fusion \textbf{is} unique up to symmetries of $M_{10}$ and $G$,
but our choices of the fusions for the other subgroups may have broken
some symmetries of $G$.)


\begin{verbatim}
    gap> s:= CharacterTable( "O10-(2)" );;
    gap> fus:= PossibleClassFusions( s, t );;
    gap> rep:= RepresentativesFusions( s, fus, Group( () ) );;
    gap> Length( rep );
    2
    gap> m10:= s;;  possm10fust:= rep;;
\end{verbatim}

The Schur multiplier of $O^-_{10}(2)$ is trivial,
thus $\pi^{-1}( M_{10} )$ is a direct product $2 \times O^-_{10}(2)$.

\begin{verbatim}
    gap> 2m10:= CharacterTable( "Cyclic", 2 ) * m10;;
\end{verbatim}

The subgroup $M_{10}$ contains elements from one class of $G$
for which we do not know yet whether it splits.
This class turns out not to split.

\begin{verbatim}
    gap> List( possm10fust, map -> Difference( map,
    >                         Union( mustsplit, mustnotsplit ) ) );
    [ [ 42 ], [ 42 ] ]
    gap> possnotsplit:= List( [ 1, 2 ], i -> ShallowCopy( mustnotsplit ) );;
    gap> for i in [ 1, 2 ] do
    >      computeContributions( m10, t, possm10fust[i], Irr( m10 ), 10^7,
    >          ShallowCopy( mustsplit ), possnotsplit[i] );
    >    od;
    #I  class 42 does not split (contribution criterion)
    #I  class 42 does not split (contribution criterion)
    gap> Set( possnotsplit );
    [ [ 4, 9, 10, 12, 13, 14, 17, 18, 20, 21, 22, 29, 35, 37, 40, 41, 42, 44, 45, 
          47, 48, 49, 57, 58, 60, 61, 62, 65, 66, 68, 69, 71, 72, 79, 80, 81, 82, 
          92, 93, 95, 96, 113, 115, 116, 118, 119 ] ]
    gap> mustnotsplit:= possnotsplit[1];;
\end{verbatim}

For each of the two possible class fusions of $M_{10}$ in $G$,
we get a unique possible class fusion of $\pi^{-1}( M_{10} )$ in $2.G$.

\begin{verbatim}
    gap> fus:= List( rep, map -> initialFusion( 2m10, 2t,
    >                 GetFusionMap( 2m10, m10 ), GetFusionMap( 2t, t ),
    >                 map, defined ) );;
    gap> ker:= ClassPositionsOfKernel( GetFusionMap( 2m10, m10 ) );;
    gap> testcharsm10:= Filtered( Irr( 2m10 ),
    >        chi -> not IsSubset( ClassPositionsOfKernel( chi ), ker ) );;
    gap> ind:= Concatenation( ind2m3, ind2m4, ind2m5,
    >              Intersection( possind2m7 ), ind2m8, ind2m9 );;
    gap> possfus:= List( fus, map -> FusionsAllowedByRestrictions( 2m10, 2t,
    >                  testcharsm10, ind, map, parametersFABR ) );
    [ [ [ 1, 3, 5, 5, 7, 10, 8, 8, 8, 10, 12, 16, 14, 17, 23, 22, 23, 31, 29, 33, 
              34, 34, 38, 36, 36, 36, 44, 44, 44, 40, 40, 36, 49, 42, 42, 38, 40, 
              42, 47, 40, 44, 47, 47, 46, 51, 53, 58, 57, 60, 66, 75, 76, 78, 78, 
              78, 78, 82, 84, 84, 86, 88, 92, 90, 90, 97, 94, 94, 90, 111, 111, 
              111, 92, 91, 94, 107, 107, 113, 107, 95, 96, 111, 110, 109, 129, 
              133, 133, 133, 135, 133, 135, 135, 148, 150, 153, 155, 153, 155, 
              158, 164, 166, 168, 178, 176, 180, 180, 191, 191, 191, 193, 195, 
              197, 197, 195, 199, 201, 2, 4, 6, 6, 7, 11, 9, 9, 9, 11, 13, 16, 
              15, 17, 24, 22, 24, 31, 30, 33, 35, 35, 39, 37, 37, 37, 45, 45, 45, 
              41, 41, 37, 50, 43, 43, 39, 41, 43, 48, 41, 45, 48, 48, 46, 52, 54, 
              59, 57, 60, 66, 75, 77, 79, 79, 79, 79, 83, 85, 85, 87, 89, 93, 90, 
              90, 98, 94, 94, 90, 112, 112, 112, 93, 91, 94, 108, 108, 114, 108, 
              95, 96, 112, 110, 109, 130, 134, 134, 134, 136, 134, 136, 136, 149, 
              151, 152, 154, 152, 154, 159, 165, 167, 169, 179, 177, 180, 180, 
              192, 192, 192, 194, 196, 198, 198, 196, 200, 202 ] ], 
      [ [ 1, 3, 5, 5, 7, 10, 8, 8, 8, 10, 12, 16, 14, 17, 23, 22, 23, 31, 29, 33, 
              34, 34, 38, 36, 36, 36, 44, 44, 44, 40, 40, 36, 49, 42, 42, 38, 40, 
              42, 47, 40, 44, 47, 47, 46, 51, 53, 58, 57, 60, 66, 75, 76, 78, 78, 
              78, 78, 82, 84, 84, 86, 88, 92, 90, 90, 97, 94, 94, 90, 111, 111, 
              111, 92, 91, 94, 107, 107, 113, 107, 95, 96, 111, 110, 109, 129, 
              133, 133, 133, 135, 133, 135, 135, 148, 150, 155, 153, 155, 153, 
              158, 164, 166, 168, 178, 176, 180, 180, 191, 191, 191, 193, 195, 
              197, 197, 195, 199, 201, 2, 4, 6, 6, 7, 11, 9, 9, 9, 11, 13, 16, 
              15, 17, 24, 22, 24, 31, 30, 33, 35, 35, 39, 37, 37, 37, 45, 45, 45, 
              41, 41, 37, 50, 43, 43, 39, 41, 43, 48, 41, 45, 48, 48, 46, 52, 54, 
              59, 57, 60, 66, 75, 77, 79, 79, 79, 79, 83, 85, 85, 87, 89, 93, 90, 
              90, 98, 94, 94, 90, 112, 112, 112, 93, 91, 94, 108, 108, 114, 108, 
              95, 96, 112, 110, 109, 130, 134, 134, 134, 136, 134, 136, 136, 149, 
              151, 154, 152, 154, 152, 159, 165, 167, 169, 179, 177, 180, 180, 
              192, 192, 192, 194, 196, 198, 198, 196, 200, 202 ] ] ]
    gap> List( possfus, Length );
    [ 1, 1 ]
    gap> poss2m10fus2t:= Concatenation( possfus );;
    gap> Set( possm10fust[1] ) = Set( possm10fust[2] );
    true
    gap> UniteSet( defined, possm10fust[1] );
\end{verbatim}

The set of induced characters is the same
for each of the two fusion candidates.

\begin{verbatim}
    gap> possind2m10:= List( poss2m10fus2t,
    >        map -> Set( InducedClassFunctionsByFusionMap( 2m10, 2t,
    >                        testcharsm10, map ) ) );;
    gap> possind2m10[1] = possind2m10[2];
    true
    gap> ind2m10:= possind2m10[1];;
\end{verbatim}

With the help of these characters, we can show that
for each of the two possible class fusions of $U_{12}$ in $G$,
we get a unique possible class fusion of $\pi^{-1}( U_{12} )$ in $2.G$.

\begin{verbatim}
    gap> fus:= List( possu12fust, map -> initialFusion( 2u12, 2t,
    >                 GetFusionMap( 2u12, u12 ), GetFusionMap( 2t, t ),
    >                 map, defined ) );;
    gap> ind:= Concatenation( ind2m3, ind2m4, ind2m5,
    >              Intersection( possind2m7 ), ind2m8, ind2m9, ind2m10 );;
    gap> possfus:= List( fus, map -> FusionsAllowedByRestrictions( 2u12, 2t,
    >                  testcharsu12, ind, map, parametersFABR ) );
    [ [ [ 1, 2, 3, 4, 5, 6, 7, 8, 9, 12, 13, 10, 11, 16, 17, 18, 19, 20, 21, 23, 
              24, 31, 34, 35, 42, 43, 42, 43, 36, 37, 47, 48, 40, 41, 38, 39, 44, 
              45, 49, 50, 55, 56, 60, 62, 63, 64, 65, 78, 79, 78, 79, 80, 81, 82, 
              83, 86, 87, 88, 89, 91, 91, 90, 90, 95, 96, 103, 104, 105, 106, 
              107, 108, 135, 136, 153, 152, 155, 154, 8, 9, 36, 37, 40, 41, 46, 
              10, 11, 8, 9, 12, 13, 90, 94, 99, 100, 101, 102, 107, 108, 109, 
              135, 136, 36, 37, 36, 37, 38, 39, 40, 41, 44, 45, 42, 43, 47, 48, 
              51, 52, 170, 171, 180, 181, 182, 183, 184, 78, 79, 78, 79, 80, 81, 
              193, 194, 195, 196, 197, 198, 90, 90, 92, 93, 94, 94, 99, 100, 101, 
              102, 111, 112, 133, 134, 153, 152, 155, 154, 8, 9, 36, 37, 40, 41, 
              46, 10, 11, 8, 9, 12, 13, 90, 94, 99, 100, 101, 102, 107, 108, 109, 
              135, 136, 36, 37, 36, 37, 38, 39, 40, 41, 44, 45, 42, 43, 47, 48, 
              51, 52, 170, 171, 180, 181, 182, 183, 184, 78, 79, 78, 79, 80, 81, 
              193, 194, 195, 196, 197, 198, 90, 90, 92, 93, 94, 94, 99, 100, 101, 
              102, 111, 112, 133, 134, 153, 152, 155, 154 ] ], 
      [ [ 1, 2, 3, 4, 5, 6, 7, 8, 9, 12, 13, 10, 11, 16, 17, 18, 19, 20, 21, 23, 
              24, 31, 34, 35, 42, 43, 42, 43, 36, 37, 47, 48, 40, 41, 38, 39, 44, 
              45, 49, 50, 55, 56, 60, 62, 63, 64, 65, 78, 79, 78, 79, 80, 81, 82, 
              83, 86, 87, 88, 89, 91, 91, 90, 90, 95, 96, 103, 104, 105, 106, 
              107, 108, 135, 136, 155, 154, 153, 152, 8, 9, 36, 37, 40, 41, 46, 
              10, 11, 8, 9, 12, 13, 90, 94, 99, 100, 101, 102, 107, 108, 109, 
              135, 136, 36, 37, 36, 37, 38, 39, 40, 41, 44, 45, 42, 43, 47, 48, 
              51, 52, 170, 171, 180, 181, 182, 183, 184, 78, 79, 78, 79, 80, 81, 
              193, 194, 195, 196, 197, 198, 90, 90, 92, 93, 94, 94, 99, 100, 101, 
              102, 111, 112, 133, 134, 155, 154, 153, 152, 8, 9, 36, 37, 40, 41, 
              46, 10, 11, 8, 9, 12, 13, 90, 94, 99, 100, 101, 102, 107, 108, 109, 
              135, 136, 36, 37, 36, 37, 38, 39, 40, 41, 44, 45, 42, 43, 47, 48, 
              51, 52, 170, 171, 180, 181, 182, 183, 184, 78, 79, 78, 79, 80, 81, 
              193, 194, 195, 196, 197, 198, 90, 90, 92, 93, 94, 94, 99, 100, 101, 
              102, 111, 112, 133, 134, 155, 154, 153, 152 ] ] ]
    gap> List( possfus, Length );
    [ 1, 1 ]
    gap> poss2u12fus2t:= Concatenation( possfus );;
    gap> Set( possu12fust[1] ) = Set( possu12fust[2] );
    true
    gap> UniteSet( defined, possu12fust[1] );
\end{verbatim}


The set of induced characters is the same
for each of the two fusion candidates.

\begin{verbatim}
    gap> possind2u12:= List( poss2u12fus2t,
    >        map -> Set( InducedClassFunctionsByFusionMap( 2u12, 2t,
    >                        testcharsu12, map ) ) );;
    gap> possind2u12[1] = possind2u12[2];
    true
    gap> ind2u12:= possind2u12[1];;
\end{verbatim}


\subsection{What do we know up to now about the table of $2.G$?}

We have determined for all except four classes of $G$ whether they split
or not.

\begin{verbatim}
    gap> posssplit:= Difference( [ 1 .. Length( orders ) ],
    >                            Union( mustsplit, mustnotsplit ) );
    [ 38, 88, 89, 90 ]
    gap> NrConjugacyClasses( t );
    126
    gap> NrConjugacyClasses( 2t );
    202
\end{verbatim}

Thus $2.G$ has at least $202$ and at most $206$ classes,
and we have to compute at least $76$ and at most $80$ faithful
irreducible characters.
If some of the above four classes split then the class fusions of the
subgroups have to be adjusted by shifting the image classes
appropriately.

For the subgroups
$U \in \{ M_3, M_4, M_5, M_7, M_8, M_9, M_{10}, U_{12} \}$ of $G$,
we have computed characters of $2.G$ by induction from $\pi^{-1}( U )$.
For the subgroup $\pi^{-1}( M_7 )$ of $2.G$,
we know only two possible sets of induced characters.
Now we are able to eliminate one of these cases.

\begin{verbatim}
    gap> ForAll( ind2u12,
    >        chi -> ForAll( possind2m7[1],
    >                   psi -> IsInt( ScalarProduct( 2t, chi, psi ) ) ) );
    false
\end{verbatim}

This means that we know the class fusion from $\pi^{-1}( M_7 )$ to $2.G$,
and the corresponding induced characters.

\begin{verbatim}
    gap> 2m7fus2t:= poss2m7fus2t[2];;
    gap> ind2m7:= possind2m7[2];;
\end{verbatim}

The known induced characters do not span the subspace of all ``faithful
class functions of $2.G$'', that is, those class functions $\psi$ with
the property $\psi( g ) = - \psi( g z )$,
since all known induced characters are zero on the preimages in $2.G$
of elements of order $19$ in $G$.

\begin{verbatim}
    gap> nothit:= Difference( [ 1 .. Length( orders ) ],
    >                 Flat( [ m3fust, m4fust, m5fust, m7fust, m8fust, m9fust,
    >                         possm10fust, possu12fust ] ) );
    [ 38, 88, 89, 90, 103, 104 ]
    gap> orders{ nothit };
    [ 8, 16, 16, 16, 19, 19 ]
\end{verbatim}

Thus we need some more characters of $2.G$.

\subsection{Additional characters of $2.G$}\label{minuscharacters}

First we compute characters with nonzero values on the classes of element
order $38$,
by inducing from the cyclic subgroups of order $38$ in $2.G$.
For that, we need some power maps of $2.G$ on the classes that contain
elements in these subgroups;
this information is determined by the corresponding power maps of $G$.

\begin{verbatim}
    gap> for p in [ 2 .. Maximum( OrdersClassRepresentatives( 2t ) ) ] do
    >      if IsPrimeInt( p ) then
    >        PowerMap( t, p );
    >        pow:= InitPowerMap( 2t, p );
    >        comp:= CompositionMaps( InverseMap( GetFusionMap( 2t, t ) ),
    >                   CompositionMaps( PowerMap( t, p ),
    >                      GetFusionMap( 2t, t ) ) );
    >        MeetMaps( pow, comp );
    >        ComputedPowerMaps( 2t )[p]:= pow;
    >      fi;
    >    od;
\end{verbatim}

Note that we need only the faithful induced characters.

\begin{verbatim}
    gap> pos:= Positions( OrdersClassRepresentatives( 2t ), 38 );
    [ 161, 163 ]
    gap> indcyc:= Filtered( InducedCyclic( 2t, pos, "all" ),
    >                       chi -> chi[1] = - chi[2] );;
\end{verbatim}

Additionally, we compute also some $5$-th symmetrizations
of induced characters.
For a given character $\chi$, the symmetrization in question
is $\chi^{5-}$, given by $\chi^{5-}(g) = (\chi(g)^5 - \chi(g)) / 5$.

For that, we need the $5$-th power map on all classes of $2.G$.
The compatibility with the $5$-th power maps of $\pi^{-1}( U )$,
for $U \in \{ M_3, M_7, M_{10} \}$,
determines the map uniquely,
and we can compute the symmetrizations.

Note that the four classes of $G$ which might split consist of elements
of $2$-power order; each of these classes maps to itself under the
$5$-th power map, and there are no other $5$-th roots than itself.
Thus the symmetrizations are independent of their potential splitting.

\begin{verbatim}
    gap> pow:= ComputedPowerMaps( 2t )[5];;
    gap> comp:= CompositionMaps( CompositionMaps( 2m3fus2t, PowerMap( 2m3, 5 ) ),
    >                            InverseMap( 2m3fus2t ) );;
    gap> MeetMaps( pow, comp );
    true
    gap> comp:= CompositionMaps( CompositionMaps( 2m7fus2t, PowerMap( 2m7, 5 ) ),
    >                            InverseMap( 2m7fus2t ) );;
    gap> MeetMaps( pow, comp );
    true
    gap> fus:= Parametrized( poss2m10fus2t );;
    gap> comp:= CompositionMaps( CompositionMaps( fus, PowerMap( 2m10, 5 ) ),
    >                            InverseMap( fus ) );;
    gap> MeetMaps( pow, comp );
    true
    gap> ForAll( pow, IsInt );
    true
    gap> Intersection( PowerMap( t, 5 ), posssplit ) = posssplit;
    true
    gap> ForAll( posssplit, i -> Positions( PowerMap( t, 5 ), i ) = [ i ] );
    true
    gap> ind:= Concatenation(
    >              [ ind2m3, ind2m4, ind2m5, ind2m7, ind2m8, ind2m9 ] );;
    gap> minus:= List( ind, chi -> MinusCharacter( chi, pow, 5 ) );;
    gap> ind:= Concatenation(
    >              [ ind2m3, ind2m4, ind2m5, minus, ind2m7, ind2m8, ind2m9,
    >                ind2u12, ind2m10, indcyc ] );;
\end{verbatim}

\subsection{The faithful irreducible characters of $2.G$}\label{thecharacters}

First we apply the LLL algorithm to our characters of $2.G$.
This yields two irreducible characters.
We project the lattice basis to the orthogonal space of the found
irreducibles, apply LLL again, and find two more irreducible characters.

\begin{verbatim}
    gap> lll:= LLL( 2t, ind );;
    gap> Length( lll.irreducibles );
    2
    gap> irr:= Set( lll.irreducibles );;
    gap> red:= Reduced( 2t, irr, lll.remainders );;
    gap> lll:= LLL( 2t, red.remainders );;
    gap> Length( lll.irreducibles );
    2
    gap> UniteSet( irr, lll.irreducibles );
    gap> red:= Reduced( 2t, irr, lll.remainders );;
    gap> lll:= LLL( 2t, red.remainders );;
    gap> Sum( lll.norms );
    772
    gap> lll:= LLL( 2t, lll.remainders, "sort" );;
    gap> Sum( lll.norms );
    729
    gap> Length( lll.irreducibles );
    0
    gap> red:= Reduced( 2t, irr, lll.remainders );;
    gap> lll:= LLL( 2t, red.remainders );;
    gap> Sum( lll.norms );
    710
\end{verbatim}

Now we have a basis for a $72$-dimensional lattice,
and compute all possible orthogonal embeddings of this lattice
into a $76$-dimensional lattice.
(Note that four more classes of $G$ may split.)

\begin{verbatim}
    gap> Length( lll.norms );
    72
    gap> gram:= MatScalarProducts( 2t, lll.remainders, lll.remainders );;
    gap> emb:= OrthogonalEmbeddings( gram, 72 + 4 );;
    gap> Length( emb.solutions );
    1
    gap> Length( emb.solutions[1] );
    72
\end{verbatim}


There is a unique solution, and we get $72$ new irreducibles.
Together with the $126$ irreducibles of the factor group $G$
and the known four faithful irreducibles,
this yields a square matrix of dimension $202$.

\begin{verbatim}
    gap> dec:= Decreased( 2t, lll.remainders,
    >                     emb.vectors{ emb.solutions[1] } );;
    gap> Length( dec.irreducibles );
    72
    gap> UniteSet( irr, dec.irreducibles );
    gap> factchars:= RestrictedClassFunctions( Irr( t ), 2t );;
    gap> irr:= Concatenation( factchars, irr );;
    gap> DimensionsMat( irr );
    [ 202, 202 ]
\end{verbatim}

Since the squares of the degrees add up to the group order,
we are done.

\begin{verbatim}
    gap> Size( 2t ) = Sum( List( irr, chi -> chi[1]^2 ) );
    true
    gap> SetIrr( 2t, List( irr, x -> Character( 2t, x ) ) );
\end{verbatim}

Finally, we check whether the character table for $2.G$
is permutation equivalent to the character table in
{\GAP}'s library of character tables,
which is equal to the table that is shown in~\cite[pp.~192--199]{CCN85}.

For that, we compute the power maps of $2.G$,
using the known power maps of $G$, the element orders of $2.G$,
and the irreducible characters of $2.G$;
it turns out that the power maps are uniquely determined by
this information.

\begin{verbatim}
    gap> for p in [ 2 .. Maximum( OrdersClassRepresentatives( 2t ) ) ] do
    >      if IsPrimeInt( p ) then
    >        poss:= PossiblePowerMaps( 2t, p,
    >                   rec( powermap:= ComputedPowerMaps( 2t )[p] ) );
    >        if Length( poss ) <> 1 then
    >          Error( "problem with ", Ordinal( p ), " power map" );
    >        fi;
    >        ComputedPowerMaps( 2t )[p]:= poss[1];
    >      fi;
    >    od;
    gap> lib:= CharacterTable( "2.2E6(2)" );;
    gap> tr:= TransformingPermutationsCharacterTables( lib, 2t );;
    gap> IsRecord( tr );
    true
\end{verbatim}

The two tables are equivalent only up to a permutation of the classes.

\begin{verbatim}
    gap> tr.columns;
    (25,26)(62,63)(67,68)(121,122)(144,145)(152,153)(154,155)(172,173)(174,
    175)(181,182)
\end{verbatim}

This is due to the choices in the function \verb|initialFusion|,
see Section~\ref{sectinitialfusion}.


\section{The character table of $2.{}^2E_6(2).2$ (March 28th, 2016)}%
\label{section22e622}

In the following, we compute the character table of a bicyclic extension
$2.G.2$ of the simple group $G = {}^2E_6(2)$
with character theoretic methods.

For that, we assume the character tables of $G$, $G.2$,
and $2.G$ (see Section~\ref{sect22e62}),
the existence of subgroups of the structures $2 \times F_4(2) \times 2$
and $3 \times 2.U_6(2).2$ in $2.G.2$,
and the knowledge of their character tables,
see Sections~\ref{2s2} and~\ref{2u2}.

\subsection{Class numbers of $2.G.2$}\label{classnumbers2G2}

The action of the field automorphism $\beta$ of $G$ on the
conjugacy classes of $G$ is determined by the class fusion of $G$ in $G.2$.
The action of the lift of $\beta$ to $2.G$ is uniquely determined by this
permutation, the table of $2.G$, and the factor fusion from $2.G$ to $G$.

\begin{verbatim}
    gap> t:= CharacterTable( "2E6(2)" );;
    gap> t2:= CharacterTable( "2E6(2).2" );;
    gap> 2t:= CharacterTable( "2.2E6(2)" );;
    gap> tfust2:= GetFusionMap( t, t2 );;
    gap> 2tfust:= GetFusionMap( 2t, t );;
    gap> orbsbeta:= Filtered( InverseMap( tfust2 ), IsList );
    [ [ 12, 13 ], [ 17, 18 ], [ 40, 41 ], [ 44, 45 ], [ 47, 48 ], [ 55, 56 ], 
      [ 61, 62 ], [ 65, 66 ], [ 68, 69 ], [ 76, 77 ], [ 79, 80 ], [ 89, 90 ], 
      [ 92, 93 ], [ 95, 96 ], [ 99, 100 ], [ 103, 104 ], [ 109, 110 ], 
      [ 115, 116 ], [ 118, 119 ], [ 123, 124 ], [ 125, 126 ] ]
    gap> beta:= Product( List( orbsbeta, x -> ( x[1], x[2] ) ) );
    (12,13)(17,18)(40,41)(44,45)(47,48)(55,56)(61,62)(65,66)(68,69)(76,77)(79,
    80)(89,90)(92,93)(95,96)(99,100)(103,104)(109,110)(115,116)(118,119)(123,
    124)(125,126)
    gap> aut:= AutomorphismsOfTable( 2t );;
    gap> Size( aut );
    256
    gap> filt:= Filtered( Elements( aut ),
    >               x -> OnTuples( 2tfust, beta ) = Permuted( 2tfust, x ) );;
    gap> Length( filt );
    1
    gap> betalift:= filt[1];;
\end{verbatim}

This determines the classes of $2.G$ in $2.G.2$.
We fix the corresponding class fusion.

\begin{verbatim}
    gap> 2tfus2t2:= [];;
    gap> max:= 0;;
    gap> for i in [ 1 .. NrConjugacyClasses( 2t ) ] do
    >      img:= i^betalift;
    >      if img = i then
    >        # no fusion
    >        max:= max + 1;
    >        2tfus2t2[i]:= max;
    >      elif i < img then
    >        # fusion of two classes
    >        max:= max + 1;
    >        2tfus2t2[i]:= max;
    >        2tfus2t2[ img ]:= max;
    >      fi;
    >    od;
    gap> max;
    174
\end{verbatim}

Now we can compute the number of classes of $2.G.2$.
If $2.G$ has $n$ classes and the lift of $\beta$ has $f$ nontrivial orbits
then $2.G.2$ has $f + 2 (n - 2f) = 2n - 3f$ classes,
of which $n - f$ classes contain the elements of $2.G$,
and $n - 2f$ classes contain the elements outside $2.G$.

In our case, we have $n = 202$ and $f = 28$,
thus $2.G.2$ has $320$ classes, of which $174$ are inner and $146$ are outer.

\begin{verbatim}
    gap> n:= NrConjugacyClasses( 2t );
    202
    gap> f:= NrMovedPoints( betalift ) / 2;
    28
    gap> 2 * n - 3 * f;  n - f;  n - 2 * f;
    320
    174
    146
\end{verbatim}

Since $G.2$ has $84$ outer classes,
we know that the preimages of exactly $146 - 84 = 62$ of these classes
consist of two classes of $2.G.2$ (``the class splits''),
and the preimages of exactly $84 - 62 = 22$ classes
form one class of $2.G.2$ (``the class does not split'').

The number of faithful irreducible characters of $2.G.2$ is $131$.

\begin{verbatim}
    gap> NrConjugacyClasses( t2 ) - Maximum( tfust2 );
    84
    gap> 320 - NrConjugacyClasses( t2 );
    131
\end{verbatim}

Concerning the splitting of classes,
we start with applying the elementary criteria from
Section~\ref{elementarysplitting}.
The splitting of classes inside $2.G$ is of course given by the
factor fusion from $2.G$ to $G$ and the class fusion of $G$ to $G.2$.

\begin{verbatim}
    gap> inv:= InverseMap( 2tfust );;
    gap> nonsplit:= PositionsProperty( inv, IsInt);;
    gap> mustnotsplit:= Set( tfust2{ nonsplit } );
    [ 4, 9, 10, 12, 13, 16, 18, 19, 20, 27, 33, 35, 36, 38, 39, 41, 43, 44, 51, 
      52, 54, 55, 58, 60, 62, 63, 69, 70, 71, 77, 78, 80, 82, 96, 98, 100 ]
    gap> split:= PositionsProperty( inv, IsList );;
    gap> mustsplit:= Set( tfust2{ split } );
    [ 1, 2, 3, 5, 6, 7, 8, 11, 14, 15, 17, 21, 22, 23, 24, 25, 26, 28, 29, 30, 
      31, 32, 34, 37, 40, 42, 45, 46, 47, 48, 49, 50, 53, 56, 57, 59, 61, 64, 65, 
      66, 67, 68, 72, 73, 74, 75, 76, 79, 81, 83, 84, 85, 86, 87, 88, 89, 90, 91, 
      92, 93, 94, 95, 97, 99, 101, 102, 103, 104, 105 ]
    gap> selfCentralizingClassesSplit( t2, mustsplit );
    #I  class 172 splits (self-centralizing)
    #I  class 180 splits (self-centralizing)
    #I  class 181 splits (self-centralizing)
    #I  class 182 splits (self-centralizing)
    #I  class 183 splits (self-centralizing)
    #I  class 184 splits (self-centralizing)
    #I  class 185 splits (self-centralizing)
    #I  class 188 splits (self-centralizing)
    #I  class 189 splits (self-centralizing)
\end{verbatim}

\subsection{Subgroups of the type $2 \times F_4(2) \times 2$ in $2.G.2$}%
\label{2s2}

By~\cite[Lemma~(3.1)]{Str76b},
we can assume that $2.G.2$ contains subgroups $2.S.2$
of the structure $2 \times F_4(2) \times 2$.
More precisely, we choose the isoclinism type of $2.G.2$ in such a way that
subgroups of this structure exist.

In order to use as much known information as possible about $2.S.2$,
we assume a subgroup $S$ of type $F_4(2)$ in $G$ that extends to
a subgroup $S.2$ of type $F_4(2) \times 2$ in $G.2$
and that is an epimorphic image of a subgroup $2.S$ of type $2 \times F_4(2)$
in $2.G$ under the natural epimorphism $\pi$ from $2.G.2$ to $2.G$.

Technically, we assume that all quadrangles in the following cube
describe commutative diagrams.


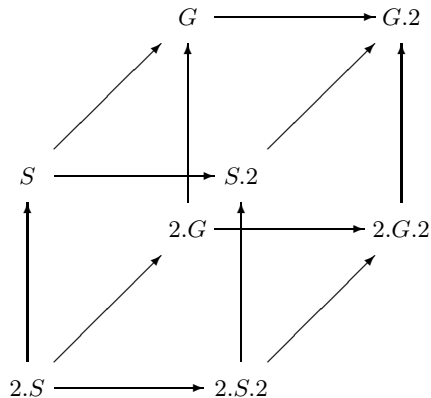
\begin{center}
\setlength{\unitlength}{2pt}
\begin{picture}(70,70)
\put( 0, 0){\makebox(0,0){$2.S$}}
\put(40, 0){\makebox(0,0){$2.S.2$}}
\put( 0,40){\makebox(0,0){$S$}}
\put(40,40){\makebox(0,0){$S.2$}}
\put(30,30){\makebox(0,0){$2.G$}}
\put(70,30){\makebox(0,0){$2.G.2$}}
\put(30,70){\makebox(0,0){$G$}}
\put(70,70){\makebox(0,0){$G.2$}}
\put( 5, 0){\vector(1,0){28}}
\put( 5,40){\vector(1,0){30}}
\put(35,30){\vector(1,0){28}}
\put(35,70){\vector(1,0){30}}
\put( 0, 5){\vector(0,1){30}}
\put(40, 5){\vector(0,1){30}}
\put(30,35){\vector(0,1){30}}
\put(70,35){\vector(0,1){30}}
\put( 5, 5){\vector(1,1){20}}
\put( 5,45){\vector(1,1){20}}
\put(45, 5){\vector(1,1){20}}
\put(45,45){\vector(1,1){20}}
\end{picture}
\end{center}

First we construct the character tables of $S$, $2.S$, $S.2$, and $2.S.2$.

\begin{verbatim}
    gap> c2:= CharacterTable( "C2" );;
    gap> s:= CharacterTable( "F4(2)" );;
    gap> s2:= c2 * s;;
    gap> 2s:= s * c2;;
    gap> 2s2:= s2 * c2;;
    gap> sfuss2:= GetFusionMap( s, s2 );;
    gap> 2sfuss:= GetFusionMap( 2s, s );;
    gap> 2s2fuss2:= GetFusionMap( 2s2, s2 );;
\end{verbatim}

The natural class fusion of $2.S$ in $2.S.2$ is not stored on the tables
if they are constructed this way.
We compute an equivalent table that stores the fusion in question.

\begin{verbatim}
    gap> 2s2alt:= c2 * 2s;;
    gap> Irr( 2s2alt ) = Irr( 2s2 );
    true
    gap> ComputedPowerMaps( 2s2alt ) = ComputedPowerMaps( 2s2 );
    true
    gap> 2sfus2s2:= GetFusionMap( 2s, 2s2alt );;
\end{verbatim}

The diagram on the front side of the cube commutes.

\begin{verbatim}
    gap> CompositionMaps( sfuss2, 2sfuss )
    >        = CompositionMaps( 2s2fuss2, 2sfus2s2 );
    true
\end{verbatim}

Concerning the left side of the cube,
we have to fix suitable class fusions from $S$ to $G$
and from $2.S$ to $2.G$.

\begin{verbatim}
    gap> sfust:= PossibleClassFusions( s, t );;
    gap> rep:= RepresentativesFusions( AutomorphismsOfTable( s ),
    >              sfust, Group( () ) );;
    gap> Length( rep );
    3
    gap> comp:= List( rep, map -> CompositionMaps( InverseMap( 2tfust ),
    >                                 CompositionMaps( map, 2sfuss ) ) );;
    gap> poss:= List( comp, map -> PossibleClassFusions( 2s, 2t,
    >                                  rec( fusionmap:= map ) ) );;
    gap> List( poss, Length );
    [ 1, 0, 0 ]
    gap> sfust:= rep[1];;
    gap> 2sfus2t:= poss[1][1];;
    gap> CompositionMaps( sfust, 2sfuss )
    >        = CompositionMaps( 2tfust, 2sfus2t );
    true
\end{verbatim}

Concerning the top side of the cube,
we have to choose a fusion from $S.2$ into $G.2$.

\begin{verbatim}
    gap> comp:= CompositionMaps( CompositionMaps( tfust2, sfust ),
    >                            InverseMap( sfuss2 ) );;
    gap> poss:= PossibleClassFusions( s2, t2, rec( fusionmap:= comp ) );;
    gap> Length(  poss );
    1
    gap> s2fust2:= poss[1];;
    gap> CompositionMaps( tfust2, sfust )
    >        = CompositionMaps( s2fust2, sfuss2 );
    true
\end{verbatim}

The diagrams on the bottom, right, and back sides of the cube involve the
character table of $2.G.2$, for which we have not yet created an
approximation.

We will use the known diagrams to derive information about the missing ones.

\begin{verbatim}
    gap> splittingClassesWithOddCentralizerIndexSplit( s2, t2, s2fust2,
    >        2s2fuss2, mustsplit );
    #I  class 106 splits (odd centralizer index)
    #I  class 107 splits (odd centralizer index)
    #I  class 114 splits (odd centralizer index)
    #I  class 115 splits (odd centralizer index)
    #I  class 117 splits (odd centralizer index)
    #I  class 133 splits (odd centralizer index)
    #I  class 116 splits (odd centralizer index)
    #I  class 118 splits (odd centralizer index)
    #I  class 119 splits (odd centralizer index)
    #I  class 148 splits (odd centralizer index)
    #I  class 147 splits (odd centralizer index)
    #I  class 156 splits (odd centralizer index)
    #I  class 155 splits (odd centralizer index)
    #I  class 134 splits (odd centralizer index)
    #I  class 140 splits (odd centralizer index)
    #I  class 143 splits (odd centralizer index)
    #I  class 149 splits (odd centralizer index)
    #I  class 176 splits (odd centralizer index)
    #I  class 175 splits (odd centralizer index)
    #I  class 157 splits (odd centralizer index)
    #I  class 174 splits (odd centralizer index)
    #I  class 177 splits (odd centralizer index)
\end{verbatim}

\subsection{Subgroups of the type $3 \times 2.U_6(2).2$ in $2.G.2$}%
\label{2u2}

We have used in Section~\ref{3u62subs} that $G$ contains subgroups of the
structure $3 \times U_6(2)$.
In fact, these groups are the full centralizers of certain elements of order
$3$ in $G$.

\begin{verbatim}
    gap> c3:= CharacterTable( "Cyclic", 3 );;
    gap> h:= c3 * CharacterTable( "U6(2)" );;
    gap> poss:= PossibleClassFusions( h, t );;
    gap> cen:= ClassPositionsOfCentre( h );
    [ 1, 47, 93 ]
    gap> Set( List( poss, x -> x{ cen } ) );
    [ [ 1, 5, 5 ] ]
    gap> SizesCentralizers( t )[5] = Size( h );
    true
\end{verbatim}

Let $x$ be such an element of order $3$ in $G$.
Since the $G$-class of $x$ cannot fuse with another class under the
action of $G.2$,
we know that the centralizer of $x$ in $G.2$ is larger by a factor of two.

We have also used that the preimage in $2.G$ of the centralizer of $x$ in $G$
has the structure $3 \times 2.U_6(2)$.
Again, this subgroup is the full centralizer of an element of order $3$
in $2.G$,
and the centralizer of this element in $2.G.2$ is larger by a factor of two.

In order to determine the structures of these centralizers,
we notice that no $U_6(2)$ type subgroup of $G$ can extend to a
subgroup of the structure $U_6(2) \times 2$ in $G.2$, because no element in
$G.2 \setminus G$ has centralizer order divisible by the order of $U_6(2)$.

\begin{verbatim}
    gap> h:= CharacterTable( "U6(2)" );;
    gap> pos:= PositionsProperty( SizesCentralizers( t2 ),
    >                             x -> x mod Size( h ) = 0 );
    [ 1, 2, 5 ]
    gap> Difference( pos, tfust2 );
    [  ]
\end{verbatim}

Hence $U_6(2)$ extends to $U_6(2).2$,
and since this $U_6(2).2$ type subgroup centralizes $x$ and intersects
trivially with $\langle x \rangle$, the centralizer of $x$ in $G.2$
has the structure $3 \times U_6(2).2$.

Analogously, the preimage of this centralizer in $2.G.2$ has the structure
$3 \times 2.U_6(2).2$.
There are two isomorphism classes of $2.U_6(2).2$ type groups,
we have to decide which one occurs.
The group $U_6(2).2$ contains exactly two classes of involutions that
do not lie in the derived subgroup $U_6(2)$.
The images of these classes in $G.2$ contain elements of subgroups
of the type $F_4(2) \times 2$ which are known to lift to involutions
in $2.G.2$.
This implies that the outer involutions of $U_6(2).2$ lift to involutions
in $2.U_6(2).2$,
and this information determines the isomorphism type of the $2.U_6(2).2$
type subgroup in question.

\begin{verbatim}
    gap> h:= c3 * CharacterTable( "U6(2).2" );;
    gap> possfus:= PossibleClassFusions( h, t2 );;
    gap> inv:= Positions( OrdersClassRepresentatives( h ), 2 );
    [ 2, 3, 4, 38, 39 ]
    gap> outerinv:= Difference( inv, ClassPositionsOfDerivedSubgroup( h ) );
    [ 38, 39 ]
    gap> imgs:= Set( List( possfus, x -> x{ outerinv } ) );
    [ [ 106, 107 ] ]
    gap> List( imgs[1], x -> Positions( s2fust2, x ) );
    [ [ 96, 98 ], [ 97, 99, 100 ] ]
    gap> h:= CharacterTable( "2.U6(2).2" );
    CharacterTable( "2.U6(2).2" )
    gap> Positions( OrdersClassRepresentatives( h ), 2 );
    [ 2, 3, 4, 5, 6, 7, 65, 66, 67, 68 ]
\end{verbatim}

Now we proceed with the subgroups $U = 3 \times U_6(2)$ of $G$,
$U.2 = 3 \times U_6(2).2$ of $G.2$,
$2.U = 3 \times 2.U_6(2)$ of $2.G$,
and $2.U.2 = 3 \times 2.U_6(2).2$ of $2.G.2$
in the same way as with $S$ and its decorations.

\begin{verbatim}
    gap> c3:= CharacterTable( "Cyclic", 3 );;
    gap> 2u:= c3 * CharacterTable( "2.U6(2)" );;
    gap> 2uorders:= OrdersClassRepresentatives( 2u );;
    gap> ker:= First( ClassPositionsOfCentre( 2u ), i -> 2uorders[i] = 2 );
    2
    gap> u:= 2u / [ 1, ker ];;
    gap> 2u2:= c3 * CharacterTable( "2.U6(2).2" );;
    gap> 2u2orders:= OrdersClassRepresentatives( 2u2 );;
    gap> ker:= First( ClassPositionsOfCentre( 2u2 ), i -> 2u2orders[i] = 2 );
    2
    gap> u2:= 2u2 / [ 1, ker ];;
    gap> 2ufusu:= GetFusionMap( 2u, u );;
    gap> 2u2fusu2:= GetFusionMap( 2u2, u2 );;
    gap> poss:= PossibleClassFusions( 2u, 2u2 );;
    gap> rep:= RepresentativesFusions( 2u, poss, Group( () ) );;
    gap> Length( rep );
    1
    gap> 2ufus2u2:= rep[1];;
    gap> ufusu2:= CompositionMaps( 2u2fusu2,
    >              CompositionMaps( 2ufus2u2, InverseMap( 2ufusu ) ) );;
\end{verbatim}

By construction, the diagram on the front side of the cube commutes.

\begin{verbatim}
    gap> CompositionMaps( ufusu2, 2ufusu )
    >        = CompositionMaps( 2u2fusu2, 2ufus2u2 );
    true
\end{verbatim}

Concerning the left side of the cube,
we get two possible class fusions from $2.U$ to $2.G$
and also two possible class fusions from $U$ to $G$.
Again, the diagram on this side commutes by construction.

\begin{verbatim}
    gap> poss:= PossibleClassFusions( 2u, 2t );;
    gap> rep:= RepresentativesFusions( 2u, poss, Group( () ) );;
    gap> Length( rep );
    2
    gap> poss2ufus2t:= rep;;
    gap> possufust:= List( rep, map -> CompositionMaps( 2tfust,
    >        CompositionMaps( map, InverseMap( 2ufusu ) ) ) );;
    gap> Length( Set( possufust ) );
    2
\end{verbatim}

Concerning the top side of the cube,
we have to compute the fusions from $U.2$ to $G.2$ compatibly
with the possible fusions from $U$ to $G$.
Fortunately, we get the same result in both cases.

\begin{verbatim}
    gap> comp1:= CompositionMaps(
    >                CompositionMaps( tfust2, possufust[1] ),
    >                    InverseMap( ufusu2 ) );;
    gap> poss1:= PossibleClassFusions( u2, t2, rec( fusionmap:= comp1 ) );;
    gap> Length( poss1 );
    2
    gap> comp2:= CompositionMaps(
    >                CompositionMaps( tfust2, possufust[2] ),
    >                    InverseMap( ufusu2 ) );;
    gap> poss2:= PossibleClassFusions( u2, t2, rec( fusionmap:= comp2 ) );;
    gap> Length(  poss2 );
    2
    gap> poss1 = poss2;
    true
    gap> rep:= RepresentativesFusions( u2, poss1, Group( () ) );;
    gap> Length( rep );
    1
    gap> u2fust2:= rep[1];;
\end{verbatim}

We use the fusions to derive new information about the class splitting.

\begin{verbatim}
    gap> splittingClassesWithOddCentralizerIndexSplit( u2, t2, u2fust2,
    >        2u2fusu2, mustsplit );
    #I  class 135 splits (odd centralizer index)
    #I  class 139 splits (odd centralizer index)
    #I  class 141 splits (odd centralizer index)
    #I  class 161 splits (odd centralizer index)
    gap> notSplittingClassesOfSubgroupDoNotSplit( 2u2fusu2, u2fust2,
    >        mustnotsplit );
    #I  class 120 does not split (as in subgroup)
    #I  class 126 does not split (as in subgroup)
    #I  class 127 does not split (as in subgroup)
    #I  class 150 does not split (as in subgroup)
    #I  class 151 does not split (as in subgroup)
    #I  class 163 does not split (as in subgroup)
    #I  class 164 does not split (as in subgroup)
    #I  class 165 does not split (as in subgroup)
    #I  class 186 does not split (as in subgroup)
    #I  class 187 does not split (as in subgroup)
\end{verbatim}

\subsection{Norms of induced characters -- a refinement}

The next step is to use the idea developed in Section~\ref{normsinduced}.
However, we are in a better situation now, in the following sense.

Let $2.H.2$ be a subgroup of $2.G.2$,
and let $\chi$ be a character of $2.H.2$ that is faithful on the central
subgroup of order two.
As before, we can compute a set of all possible values for the contribution
of each class of $G.2$ to the norm of the induced character
$\psi = \chi^{2.G.2}$, and considering all combinations of these values
may yield new information on the class splitting.
The number of possible values for a class $C$ of $G.2$ depends on the number
of classes of $H.2$ that fuse into $C$.

Here we can do better,
because we can compute $\psi_{2.G} = (\chi_{2.H})^{2.G}$ and its norm,
hence we have to deal with the combinatorial effort only for the classes
outside $2.G$.
For that, we have to refine the {\GAP} functions from
Section~\ref{normsinduced}.


The functions \verb|contributionData|, \verb|integralContributions|,
and \verb|evaluateContributions| can remain unchanged.
Just \verb|computeContributions| must be changed, in order to deal with the
different situation.

\begin{verbatim}
    gap> computeContributions2:= function( s2, 2s, t2, 2t, s2fust2, 2sfus2s2,
    >                                2sfus2t, tfust2,
    >                                characters, bound,
    >                                mustsplit, mustnotsplit, proj )
    >      local inv, i, known, candidates, r, psi, res;
    > 
    >      inv:= InverseMap( s2fust2 );
    > 
    >      repeat
    >        for i in Union( mustnotsplit, tfust2 ) do
    >          # The induced character is either zero at the preimage of 'i',
    >          # and there is no contribution to the norm,
    >          # or the preimage of 'i' lies inside the subgroup of index 2.
    >          Unbind( inv[i] );
    >        od;
    >        known:= [ ShallowCopy( mustsplit ), ShallowCopy( mustnotsplit ) ];
    >        candidates:= [];
    >        for i in [ 1 .. Length( characters ) ] do
    >          r:= contributionData( s2, t2, inv, characters[i]{ proj },
    >                                mustsplit );
    >          if r.size < bound then
    >            # Restrict the character to 2.U, and induce it to 2.G.
    >            psi:= InducedClassFunctionsByFusionMap( 2s, 2t,
    >                      [ characters[i]{ 2sfus2s2 } ], 2sfus2t )[1];
    >            r.safepart:= r.safepart + ScalarProduct( 2t, psi, psi ) / 2;
    >            Add( candidates, r );
    >          fi;
    >        od;
    >        SortParallel( List( candidates, r -> r.size ), candidates );
    >        for r in candidates do
    >          res:= integralContributions( r );
    >          if Length( res ) = 0 then
    >            Error( "no solution" );
    >          fi;
    >          evaluateContributions( r, res, s2fust2, mustsplit, mustnotsplit );
    >          oddRootsOfSplittingClassesSplit( t2, mustsplit );
    >        od;
    >      until known = [ mustsplit, mustnotsplit ];
    >    end;;
\end{verbatim}

Let us apply this criterion to the subgroup $2.S.2$ of type
$2 \times F_4(2) \times 2$ of $2.G.2$.

\begin{verbatim}
    gap> ker:= ClassPositionsOfKernel( 2s2fuss2 );;
    gap> testchars:= Filtered( Irr( 2s2 ), x -> x[ ker[1] ] <> x[ ker[2] ]  );;
    gap> computeContributions2( s2, 2s, t2, 2t, s2fust2, 2sfus2s2, 2sfus2t,
    >        tfust2, testchars, 10^7, mustsplit, mustnotsplit,
    >        ProjectionMap( 2s2fuss2 ) );
    #I  class 109 splits (contribution criterion)
    #I  class 136 splits (contribution criterion)
    #I  class 158 splits (5th root of 109)
    #I  class 173 splits (7th root of 109)
    #I  class 110 splits (contribution criterion)
    #I  class 108 splits (contribution criterion)
    #I  class 137 splits (contribution criterion)
    #I  class 159 splits (contribution criterion)
    #I  class 138 splits (3rd root of 108)
    #I  class 111 splits (contribution criterion)
    #I  class 112 splits (contribution criterion)
    #I  class 142 splits (3rd root of 111)
    #I  class 144 splits (3rd root of 111)
    #I  class 145 splits (3rd root of 112)
\end{verbatim}



\subsection{Element orders in $2.G.2$}\label{elementorders2G2}

Before we compute the class fusions from the subgroups $2.S.2$ and
$2.U.2$ to $2.G.2$,
we determine the element orders of $2.G.2$.
As in Section~\ref{elementorders}, we have to consider only the classes
of involutions in $G.2$,
and we get that all of them lift to involutions in $2.G.2$
because the same holds for the involutions in $S.2$.

\begin{verbatim}
    gap> orders:= OrdersClassRepresentatives( t2 );;
    gap> invol:= Positions( orders, 2 );;
    gap> Difference( invol, s2fust2 );
    [  ]
\end{verbatim}

\subsection{The class fusion from $2 \times F_4(2) \times 2$}

Currently $14$ classes of $G.2$ are left which contain elements of $S.2$
and for which we do not know whether they split.

\begin{verbatim}
    gap> open:= Difference( s2fust2, Union( mustsplit, mustnotsplit ) );
    [ 113, 121, 122, 123, 125, 128, 129, 131, 132, 146, 153, 154, 162, 166 ]
\end{verbatim}


We are going to use norms of induced characters for determining the splitting.
In fact, we need not consider all $2^{14}$ combinations from the above set,
because the third power map connects some of the classes.

\begin{verbatim}
    gap> PowerMap( t2, 3 ){ [ 146, 162, 166 ] };
    [ 113, 123, 121 ]
    gap> poss:= Filtered( Combinations( open ),
    >               x -> ( not 113 in x or 146 in x ) and
    >                    ( not 123 in x or 162 in x ) and
    >                    ( not 121 in x or 166 in x ) );;
    gap> Length( poss );
    6912
\end{verbatim}

For the $6912$ cases, we create a preliminary character table head
of $2.G.2$, and compute possible class fusions from $2.S.2$
into this table by checking norms of induced characters.

The functions \verb|tableHead| and \verb|useInducedClassFunction|
from Section~\ref{sectinitialfusion} can be used also here,
but we want to use two commutative diagrams and not just one
when we initialize the class fusion.
Therefore, we provide a refined version \verb|initialFusion2|.

\begin{verbatim}
    gap> initialFusion2:= function( 2s2, 2t2,
    >                               2s2fuss2, 2t2fust2, s2fust2,
    >                               2sfus2s2, 2tfus2t2, 2sfus2t,
    >                               defined )
    >      local fus, comp, pre, imgs;
    > 
    >      # Use element orders and centralizer orders.
    >      fus:= InitFusion( 2s2, 2t2 );
    > 
    >      # Use the commutative diagram on the right of the cube.
    >      comp:= CompositionMaps( InverseMap( 2t2fust2 ),
    >                 CompositionMaps( s2fust2, 2s2fuss2 ) );
    >      if not MeetMaps( fus, comp ) then
    >        return fail;
    >      fi;
    > 
    >      # Use the commutative diagram on the bottom of the cube.
    >      comp:= CompositionMaps( 2tfus2t2,
    >                 CompositionMaps( 2sfus2t, InverseMap( 2sfus2s2 ) ) );
    > 
    >      if not MeetMaps( fus, comp ) then
    >        return fail;
    >      fi;
    > 
    >      # Define classes that are not yet defined.
    >      defined:= ShallowCopy( defined );
    >      for pre in InverseMap( 2s2fuss2 ) do
    >        if IsList( pre ) then
    >          imgs:= fus{ pre };
    >          if imgs[1] = imgs[2] and IsList( imgs[1] )
    >             and Intersection( defined, imgs[1] ) = [] then
    >            # The classes in preimage and image split, and we may choose.
    >            fus[ pre[1] ]:= imgs[1][1];
    >            fus[ pre[2] ]:= imgs[1][2];
    >            UniteSet( defined, imgs[1] );
    >          fi;
    >        elif IsList( fus[ pre ] ) then
    >          # The class splits in the image but not in the preimage,
    >          # something must be wrong.
    >          return fail;
    >        fi;
    >      od;
    > 
    >      return fus;
    >    end;;
\end{verbatim}

Now we check the possible splittings that can be distinguished
by characters of $2.S.2$.

\begin{verbatim}
    gap> good:= [];;
    gap> ker:= ClassPositionsOfKernel( GetFusionMap( 2s2, s2 ) );;
    gap> testcharss:= Filtered( Irr( 2s2 ),
    >        chi -> not IsSubset( ClassPositionsOfKernel( chi ), ker ) );;
    gap> runOneTest2:= function( s2, 2s2, t2, 2t2, s2fust2,
    >                            2sfus2s2, 2tfus2t2, 2sfus2t,
    >                            testchars, defined )
    >      local fus, pos, l, chi;
    > 
    >      fus:= initialFusion2( 2s2, 2t2,
    >                            GetFusionMap( 2s2, s2 ),
    >                            GetFusionMap( 2t2, t2 ),
    >                            s2fust2,
    >                            2sfus2s2, 2tfus2t2, 2sfus2t,
    >                            defined );
    > 
    >      # Process the irreducible characters,
    >      # ordered by increasing indeterminateness.
    >      pos:= PositionsProperty( fus, IsList );
    >      testchars:= ShallowCopy( testchars );
    >      l:= - List( testchars, x -> Number( pos, i -> x[i] = 0 ) );
    >      SortParallel( l, testchars );
    >      for chi in testchars do
    >        if useInducedClassFunction( 2s2, 2t2, chi,
    >               GetFusionMap( 2s2, s2 ), fus ) = false then
    >          # This splitting is not possible.
    >          return fail;
    >        fi;
    >      od;
    >      return fus;
    >    end;;
    gap> defined:= [];;
    gap> for choice in poss do
    >      2t2:= tableHead( t2, Union( mustsplit, choice ), [], [] );
    >      fus:= runOneTest2( s2, 2s2, t2, 2t2, s2fust2,
    >                   2sfus2s2, 2tfus2t2, 2sfus2t,
    >                testcharss, defined );
    >      if fus <> fail then
    >        Add( good, choice );
    >      fi;
    >    od;
    gap> Length( good );
    1
\end{verbatim}


We get a unique solution for the splitting.
Thus we update our lists.

\begin{verbatim}
    gap> choice:= good[1];
    [ 113, 121, 122, 123, 125, 131, 146, 153, 154, 162, 166 ]
    gap> UniteSet( mustsplit, choice );
    gap> UniteSet( mustnotsplit, Difference( open, choice ) );
    gap> oddRootsOfSplittingClassesSplit( t2, mustsplit );
    #I  class 168 splits (3rd root of 122)
    #I  class 171 splits (3rd root of 131)
\end{verbatim}

Because of the two splitting classes which were found in the end,
we have to create a new character table object,
and the fusion from $2.S.2$ must be adjusted.

\begin{verbatim}
    gap> 2t2:= tableHead( t2, mustsplit, [], [] );;
    gap> 2s2fus2t2:= runOneTest2( s2, 2s2, t2, 2t2, s2fust2,
    >                       2sfus2s2, 2tfus2t2, 2sfus2t,
    >                       testcharss, defined );;
    gap> NrConjugacyClasses( 2t2 );
    320
    gap> defined:= Set( 2s2fus2t2 );;
    gap> inds:= Set( InducedClassFunctionsByFusionMap( 2s2, 2t2, testcharss,
    >                                                  2s2fus2t2 ) );;
\end{verbatim}

Now we have found all class splittings, according to the considerations
in Section~\ref{classnumbers2G2}.

\subsection{The class fusion from $3 \times 2.U_6(2).2$}

Since no further class splittings can occur,
we compute the possible class fusions from $2.U.2$ to $2.G.2$ that are
compatible with the restrictions of the characters induced from $2.S.2$.
The first approximation given by the commutative diagrams is the same for
both possible fusions from $2.U$ to $2.G$,
and altogether $32$ possible fusions are compatible with the restrictions.

\begin{verbatim}
    gap> poss2u2fus2t2:= List( poss2ufus2t, map ->
    >                              initialFusion2( 2u2, 2t2,
    >                                  GetFusionMap( 2u2, u2 ),
    >                                  GetFusionMap( 2t2, t2 ),
    >                                  u2fust2,
    >                                  2ufus2u2, 2tfus2t2, map,
    >                                  defined ) );;
    gap> Length( Set( poss2u2fus2t2 ) );
    1
    gap> fus:= poss2u2fus2t2[1];;
    gap> ker:= ClassPositionsOfKernel( GetFusionMap( 2u2, u2 ) );;
    gap> testcharsu:= Filtered( Irr( 2u2 ),
    >        chi -> not IsSubset( ClassPositionsOfKernel( chi ), ker ) );;
    gap> possfus:= FusionsAllowedByRestrictions( 2u2, 2t2, testcharsu,
    >                  inds, fus, parametersFABR );;
    gap> List( possfus, Indeterminateness );
    [ 16, 16 ]
\end{verbatim}


Exactly two of these maps yield integral norms of characters induced from
$2.U.2$.
We compute the two possible lists of induced characters.

\begin{verbatim}
    gap> good:= [];;
    gap> for map in Set( Concatenation( List( possfus, ContainedMaps ) ) ) do
    >      indu:= Set( InducedClassFunctionsByFusionMap( 2u2, 2t2,
    >                      testcharsu, map ) );
    >      if ForAll( indu, x -> IsInt( ScalarProduct( 2t2, x, x ) ) ) then
    >        Add( good, map );
    >      fi;
    >    od;
    gap> Length( good );
    2
    gap> poss2u2fus2t2:= good;;
    gap> indu:= List( good,
    >                 map -> Set( InducedClassFunctionsByFusionMap( 2u2, 2t2,
    >                                 testcharsu, map ) ) );;
    gap> List( indu, Length );
    [ 98, 98 ]
    gap> Length( Intersection( indu ) );
    34
    gap> Set( poss2u2fus2t2[1] ) = Set( poss2u2fus2t2[2] );
    true
    gap> UniteSet( defined, Set( poss2u2fus2t2[1] ) );
\end{verbatim}

\subsection{Approximations for some power maps of $2.G.2$}

In order to compute characters of $2.G.2$ which distinguish the preimages
of splitting classes that contain no elements of $S.2$ or $U.2$,
we will compute some symmetrizations of the known characters
(cf. Section~\ref{minuscharacters}).
For that, we need approximations of the $3$-rd and the $5$-th power map
of $2.G.2$.

We start with the $5$-th power map.

\begin{verbatim}
    gap> 2t2fust2:= GetFusionMap( 2t2, t2 );;
    gap> p:= 5;;
    gap> pow5:= InitPowerMap( 2t2, p );;
    gap> comp:= CompositionMaps( InverseMap( 2t2fust2 ),
    >            CompositionMaps( PowerMap( t2, p ), 2t2fust2 ) );;
    gap> MeetMaps( pow5, comp );
    true
    gap> comp:= CompositionMaps( 2tfus2t2,
    >            CompositionMaps( PowerMap( 2t, p ),
    >                InverseMap( 2tfus2t2 ) ) );;
    gap> MeetMaps( pow5, comp );
    true
    gap> comp:= CompositionMaps( 2s2fus2t2,
    >            CompositionMaps( PowerMap( 2s2, p ),
    >                InverseMap( 2s2fus2t2 ) ) );;
    gap> MeetMaps( pow5, comp );
    true
    gap> para:= Parametrized( poss2u2fus2t2 );;
    gap> comp:= CompositionMaps( para,
    >            CompositionMaps( PowerMap( 2u2, p ),
    >                InverseMap( para ) ) );;
    gap> MeetMaps( pow5, comp );
    true
    gap> Indeterminateness( pow5 );
    4096
\end{verbatim}

We may {\bf choose} the images under the $5$-th power map of those classes
whose images in $G.2$ have not been defined yet and which have element order
divisible by $5$.

\begin{verbatim}
    gap> ambig:= PositionsProperty( pow5, IsList );
    [ 283, 284, 287, 288, 307, 308, 309, 310, 317, 318, 319, 320 ]
    gap> ambig:= Filtered( ambig,
    >                i -> OrdersClassRepresentatives( 2t2 )[i] mod 5 = 0 );
    [ 309, 310, 319, 320 ]
    gap> Intersection( ambig, defined );
    [  ]
    gap> pow5{ ambig };
    [ [ 204, 205 ], [ 204, 205 ], [ 227, 228 ], [ 227, 228 ] ]
    gap> pow5{ ambig }:= [ 204, 205, 227, 228 ];;
    gap> UniteSet( defined, ambig );
\end{verbatim}

We cannot improve the $5$-th power map with the currently known induced
characters.
This means that we can compute the symmetrizations $\chi^{5-}$ for all
these characters $\chi$ with the known approximation of the $5$-th power map.

The $3$-rd power map is initialized analogously.
Additionally we use the fact that the $3$-rd and the $5$-th power maps
commute.

\begin{verbatim}
    gap> p:= 3;;
    gap> pow3:= InitPowerMap( 2t2, p );;
    gap> comp:= CompositionMaps( InverseMap( 2t2fust2 ),
    >            CompositionMaps( PowerMap( t2, p ), 2t2fust2 ) );;
    gap> MeetMaps( pow3, comp );
    true
    gap> comp:= CompositionMaps( 2tfus2t2,
    >            CompositionMaps( PowerMap( 2t, p ),
    >                InverseMap( 2tfus2t2 ) ) );;
    gap> MeetMaps( pow3, comp );
    true
    gap> comp:= CompositionMaps( 2s2fus2t2,
    >            CompositionMaps( PowerMap( 2s2, p ),
    >                InverseMap( 2s2fus2t2 ) ) );;
    gap> MeetMaps( pow3, comp );
    true
    gap> para:= Parametrized( poss2u2fus2t2 );;
    gap> comp:= CompositionMaps( para,
    >            CompositionMaps( PowerMap( 2u2, p ),
    >                InverseMap( para ) ) );;
    gap> MeetMaps( pow3, comp );
    true
    gap> Indeterminateness( pow3 );
    65536
    gap> CommutativeDiagram( pow3, pow5, pow5, pow3 );
    rec( imp1 := [ 309, 310, 319, 320 ], imp2 := [  ], imp3 := [  ], imp4 := [  ] 
     )
    gap> Indeterminateness( pow3 );
    4096
    gap> ambig:= PositionsProperty( pow3, IsList );
    [ 239, 240, 273, 274, 283, 284, 287, 288, 307, 308, 317, 318 ]
    gap> ambig:=  Difference( ambig, defined );
    [ 283, 284, 287, 288, 307, 308, 317, 318 ]
    gap> pow3{ ambig };
    [ [ 206, 207 ], [ 206, 207 ], [ 218, 219 ], [ 218, 219 ], [ 231, 232 ], 
      [ 231, 232 ], [ 317, 318 ], [ 317, 318 ] ]
    gap> pow3{ [ 283, 284, 287, 288, 307, 308 ] }:=
    >              [ 206, 207, 218, 219, 231, 232 ];;
    gap> CommutativeDiagram( pow3, pow5, pow5, pow3 );
    rec( imp1 := [  ], imp2 := [  ], imp3 := [ 283, 284, 287, 288, 307, 308 ], 
      imp4 := [  ] )
\end{verbatim}

We will need more information about the $3$-rd power map,
at the cost of case distinctions.
For each set of candidates of characters induced from $2.U.2$,
we get an approximation of the $3$-rd power map.

\begin{verbatim}
    gap> poss3:= [];;
    gap> for possindu in indu do
    >      testchars:= Concatenation( inds, possindu );
    >      Add( poss3, PowerMapsAllowedBySymmetrizations( 2t2, testchars,
    >                      testchars, StructuralCopy( pow3 ), p,
    >                      parametersFABR ) );
    >    od;
    gap> List( poss3, Length );
    [ 1, 1 ]
    gap> poss3:= List( poss3, l -> l[1] );;
    gap> List( poss3, Indeterminateness );
    [ 4, 4 ]
\end{verbatim}

For other purposes,
we will also need an approximation of the $2$-nd power map
of $2.G.2$.

\begin{verbatim}
    gap> p:= 2;;
    gap> pow2:= InitPowerMap( 2t2, p );;
    gap> comp:= CompositionMaps( InverseMap( 2t2fust2 ),
    >               CompositionMaps( PowerMap( t2, p ), 2t2fust2 ) );;
    gap> MeetMaps( pow2, comp );
    true
    gap> comp:= CompositionMaps( 2tfus2t2,
    >               CompositionMaps( PowerMap( 2t, p ),
    >                   InverseMap( 2tfus2t2 ) ) );;
    gap> MeetMaps( pow2, comp );
    true
    gap> comp:= CompositionMaps( 2s2fus2t2,
    >               CompositionMaps( PowerMap( 2s2, p ),
    >                   InverseMap( 2s2fus2t2 ) ) );;
    gap> MeetMaps( pow2, comp );
    true
    gap> para:= Parametrized( poss2u2fus2t2 );;
    gap> comp:= CompositionMaps( para,
    >               CompositionMaps( PowerMap( 2u2, p ),
    >                   InverseMap( para ) ) );;
    gap> MeetMaps( pow2, comp );
    true
    gap> Indeterminateness( pow2 );
    131072
    gap> CommutativeDiagram( pow3, pow2, pow2, pow3 );
    rec( imp1 := [  ], imp2 := [  ], 
      imp3 := [ 283, 284, 287, 288, 307, 308, 319, 320 ], imp4 := [  ] )
    gap> Indeterminateness( pow2 );
    512
    gap> CommutativeDiagram( pow5, pow2, pow2, pow5 );
    rec( imp1 := [  ], imp2 := [  ], imp3 := [ 309, 310 ], imp4 := [  ] )
    gap> Indeterminateness( pow2 );
    128
\end{verbatim}

\subsection{The faithful irreducible characters of $2.G.2$}

A few faithful irreducible characters of $2.G.2$ are obtained by
induction from $2.G$.

\begin{verbatim}
    gap> indt:= Set( InducedClassFunctionsByFusionMap( 2t, 2t2,
    >               Filtered( Irr( 2t ), x -> x[1] <> x[2] ), 2tfus2t2 ) );;
    gap> irr:= Filtered( indt, x -> ScalarProduct( 2t2, x, x ) = 1 );;
    gap> Length( irr ); 
    7
\end{verbatim}

A few more characters are obtained by the LLL algorithm.

\begin{verbatim}
    gap> red:= Reduced( 2t2, irr, Concatenation( indt, inds ) );;
    gap> Length( red.irreducibles );
    0
    gap> lll:= LLL( 2t2, red.remainders, 99/100 );;
    gap> Length( lll.irreducibles );
    4
    gap> UniteSet( irr, lll.irreducibles );
    gap> Length( irr ); 
    11
    gap> missing:= NrConjugacyClasses( 2t2 ) - NrConjugacyClasses( t2 )
    >                  - Length( irr );
    120
    gap> redindu:= List( indu, l -> ReducedCharacters( 2t2, irr, l ) );;
    gap> List( redindu, r -> r.irreducibles );
    [ [  ], [  ] ]
    gap> redindu:= List( redindu, r -> r.remainders );;
\end{verbatim}

We have two cases to consider,
one for each of the two sets of characters induced from $2.U.2$;
note that we have to prescribe the $3$-rd power map that fits to the
chosen list of induced characters.

First we try the first set of induced characters and the first power map
candidate.

\begin{verbatim}
    gap> inducand:= redindu[1];;
    gap> testchars:= Concatenation( irr, red.remainders, inducand );;
    gap> minus5:= List( testchars, x -> MinusCharacter( x, pow5, 5 ) );;
    gap> minus3:= List( testchars, x -> MinusCharacter( x, poss3[1], 3 ) );;
    gap> minus:= Reduced( 2t2, irr, Concatenation( minus5, minus3 ) );;
    gap> Length( minus.irreducibles );
    0
    gap> lll2:= LLL( 2t2, Concatenation( lll.remainders, inducand,
    >                         minus.remainders ), 99/100 );;
    gap> Length( lll2.irreducibles );
    0
    gap> Length( lll2.norms );
    119
    gap> gram:= MatScalarProducts( 2t2, lll2.remainders, lll2.remainders );;
    gap> emb:= OrthogonalEmbeddings( gram, missing );;
    gap> Length( emb.solutions );
    1
    gap> dec:= Decreased( 2t2, lll2.remainders, emb.vectors{ emb.solutions[1] } );;
    gap> Length( dec.irreducibles );
    118
    gap> Length( dec.remainders );
    1
\end{verbatim}

We have found all except two irreducibles.
What are the two missing characters?
Of course they are the two extensions of one irreducible character of $2.G$.

\begin{verbatim}
    gap> redt:= Reduced( 2t2, Concatenation( irr, dec.irreducibles ), indt );;
    gap> Length( redt.remainders );
    1
    gap> redt.remainders[1] in indt;
    true
\end{verbatim}

The two missing irreducibles are the only ones that distinguish the two
classes of element order $56$.

\begin{verbatim}
    gap> split:= Union( Filtered( InverseMap( 2t2fust2 ), IsList ) );;
    gap> Filtered( split, i -> ForAll( Concatenation( irr, dec.irreducibles ),
    >                                  x -> x[i] = 0 ) );
    [ 317, 318 ]
    gap> Positions( OrdersClassRepresentatives( 2t2 ), 56 );
    [ 317, 318 ]
    gap> SizesCentralizers( 2t2 )[317];
    112
\end{verbatim}

This means that the values of the two missing irreducibles on these classes
are $\pm x$ such that $2 |x|^2 = |C_{2.G.2}(g)|/2 = 56$ holds,
where $g$ is an element of order $56$ in $2.G.2$.
As a consequence, we know that either $x = \pm 2 \sqrt{7}$ or
$x = \pm 2 \sqrt{-7}$.
Moreover, the orthogonality relations imply that the classes of element
order $56$ are the only classes outside $2.G$ where two missing irreducible
characters do not vanish.

We have to decide which of the two cases occurs,
that is, whether the classes of element order $56$ are real or not.
For that, we use symmetrizations $\chi^{2-}$, so we need the approximation
of the $2$-nd power map of $2.G.2$ computed above.

\begin{verbatim}
    gap> factirr:= List( Irr( t2 ), x -> x{ 2t2fust2 } );;
    gap> poss2:= PowerMapsAllowedBySymmetrizations( 2t2, factirr,
    >                dec.irreducibles, StructuralCopy( pow2 ), 2,
    >                parametersFABR );;
    gap> Length( poss2 );
    1
    gap> Indeterminateness( poss2[1] );
    1
    gap> cand:= ShallowCopy( redt.remainders[1] / 2 );;
    gap> cand{ [ 317, 318 ] }:= [ 1, -1 ] * ( 2 * Sqrt(-7) );;
    gap> minus2:= MinusCharacter( cand, poss2[1], 2 );;
    gap> ForAll( Flat( MatScalarProducts( 2t2, factirr, [ minus2 ] ) ), IsInt );
    false
    gap> cand{ [ 317, 318 ] }:= [ 1, -1 ] * ( 2 * Sqrt(7) );;
    gap> minus2:= MinusCharacter( cand, poss2[1], 2 );;
    gap> ForAll( Flat( MatScalarProducts( 2t2, factirr, [ minus2 ] ) ), IsInt );
    true
\end{verbatim}

This proves that the classes of element order $56$ are real.
We set the list of irreducibles in the character table of $2.G.2$,
compute the power maps, and compare the table with the {\ATLAS} table
from {\GAP}'s library.
It turns out that this first solution is equivalent to the library table.

\begin{verbatim}
    gap> cand2:= ShallowCopy( cand );;
    gap> cand2{ [ 317, 318 ] }:= [ -1, 1 ] * ( 2 * Sqrt(7) );;
    gap> SetIrr( 2t2, Concatenation( factirr, irr, dec.irreducibles,
    >                     [ cand, cand2 ] ) );
    gap> for p in [ 2 .. Maximum( OrdersClassRepresentatives( 2t2 ) ) ] do
    >      if IsPrimeInt( p ) then
    >        if p = 2 then
    >          poss:= PossiblePowerMaps( 2t2, p,
    >                     rec( powermap:= poss2[1] ) );
    >        elif p = 3 then
    >          poss:= PossiblePowerMaps( 2t2, p,
    >                     rec( powermap:= StructuralCopy( pow3 ) ) );
    >        elif p = 5 then
    >          poss:= PossiblePowerMaps( 2t2, p,
    >                     rec( powermap:= StructuralCopy( pow5 ) ) );
    >        else
    >          poss:= PossiblePowerMaps( 2t2, p );
    >        fi;
    >        if Length( poss ) <> 1 then
    >          Error( "not expected" );
    >        fi;
    >        ComputedPowerMaps( 2t2 )[p]:= poss[1];
    >      fi;
    >    od;
    gap> lib:= CharacterTable( "2.2E6(2).2" );;
    gap> tr:= TransformingPermutationsCharacterTables( lib, 2t2 );;
    gap> tr.columns;
    (177,178)(179,180)(183,184)(185,186)(189,190)(195,196)(199,200)(201,202)(204,
    205)(206,207)(208,209)(218,219)(223,224)(225,226)(227,228)(231,232)(233,
    234)(235,236)(241,242)(243,244)(245,246)(247,248)(253,254)(258,259)(260,
    261)(266,267)(268,269)(273,274)(275,276)(280,281)(283,284)(287,288)(293,
    294)(299,300)(307,308)(309,310)
\end{verbatim}

Now we repeat the above arguments to the second set of induced characters.
Again, we get one complete character table which is equivalent to the
library table.
This proves the correctness of the {\ATLAS} table.

\begin{verbatim}
    gap> 2t2:= tableHead( t2, mustsplit, [], [] );;
    gap> inducand:= redindu[2];;
    gap> testchars:= Concatenation( irr, red.remainders, inducand );;
    gap> minus5:= List( testchars, x -> MinusCharacter( x, pow5, 5 ) );;
    gap> minus3:= List( testchars, x -> MinusCharacter( x, poss3[2], 3 ) );;
    gap> minus:= Reduced( 2t2, irr, Concatenation( minus5, minus3 ) );;
    gap> Length( minus.irreducibles );
    0
    gap> lll2:= LLL( 2t2, Concatenation( lll.remainders, inducand,
    >                         minus.remainders ), 99/100 );;
    gap> Length( lll2.irreducibles );
    0
    gap> Length( lll2.norms );
    119
    gap> gram:= MatScalarProducts( 2t2, lll2.remainders, lll2.remainders );;
    gap> emb:= OrthogonalEmbeddings( gram, missing );;
    gap> Length( emb.solutions );
    1
    gap> dec:= Decreased( 2t2, lll2.remainders, emb.vectors{ emb.solutions[1] } );;
    gap> Length( dec.irreducibles );
    118
    gap> Length( dec.remainders );
    1
    gap> redt:= Reduced( 2t2, Concatenation( irr, dec.irreducibles ), indt );;
    gap> Length( redt.remainders );
    1
    gap> redt.remainders[1] in indt;
    true
    gap> poss2:= PowerMapsAllowedBySymmetrizations( 2t2, factirr,
    >                dec.irreducibles, StructuralCopy( pow2 ), 2,
    >                parametersFABR );;
    gap> Length( poss2 );
    1
    gap> Indeterminateness( poss2[1] );
    1
    gap> cand:= ShallowCopy( redt.remainders[1] / 2 );;
    gap> cand{ [ 317, 318 ] }:= [ 1, -1 ] * ( 2 * Sqrt(-7) );;
    gap> minus2:= MinusCharacter( cand, poss2[1], 2 );;
    gap> ForAll( Flat( MatScalarProducts( 2t2, factirr, [ minus2 ] ) ), IsInt );
    false
    gap> cand{ [ 317, 318 ] }:= [ 1, -1 ] * ( 2 * Sqrt(7) );;
    gap> minus2:= MinusCharacter( cand, poss2[1], 2 );;
    gap> ForAll( Flat( MatScalarProducts( 2t2, factirr, [ minus2 ] ) ), IsInt );
    true
    gap> cand2:= ShallowCopy( cand );;
    gap> cand2{ [ 317, 318 ] }:= [ -1, 1 ] * ( 2 * Sqrt(7) );;
    gap> SetIrr( 2t2, Concatenation( factirr, irr, dec.irreducibles,
    >                     [ cand, cand2 ] ) );
    gap> for p in [ 2 .. Maximum( OrdersClassRepresentatives( 2t2 ) ) ] do
    >      if IsPrimeInt( p ) then
    >        if p = 2 then
    >          poss:= PossiblePowerMaps( 2t2, p,
    >                     rec( powermap:= poss2 ) );
    >        elif p = 3 then
    >          poss:= PossiblePowerMaps( 2t2, p,
    >                     rec( powermap:= StructuralCopy( pow3 ) ) );
    >        elif p = 5 then
    >          poss:= PossiblePowerMaps( 2t2, p,
    >                     rec( powermap:= StructuralCopy( pow5 ) ) );
    >        else
    >          poss:= PossiblePowerMaps( 2t2, p );
    >        fi;
    >        if Length( poss ) <> 1 then
    >          Error( "not expected" );
    >        fi;
    >        ComputedPowerMaps( 2t2 )[p]:= poss[1];
    >      fi;
    >    od;
    gap> tr:= TransformingPermutationsCharacterTables( lib, 2t2 );;
    gap> tr.columns;
    (177,178)(179,180)(183,184)(185,186)(189,190)(195,196)(199,200)(201,202)(204,
    205)(206,207)(208,209)(218,219)(223,224)(225,226)(227,228)(231,232)(233,
    234)(235,236)(239,240)(241,242)(243,244)(245,246)(247,248)(253,254)(258,
    259)(260,261)(266,267)(268,269)(275,276)(280,281)(283,284)(287,288)(293,
    294)(299,300)(307,308)(309,310)
\end{verbatim}

\section{The character table of $2.B$ (May 8th, 2016)}%
\label{section2B}

In the following, we compute the character table of the double cover
$2.B$ of the sporadic simple Baby Monster group $B$,
with character theoretic methods.

For that, we assume the character table of $B$,
the existence of subgroups of the structures $2^2.{}^2E_6(2).2$
and $2 \times Th$ in $2.B$,
and the knowledge of their character tables;
see Sections~\ref{V42E622in2B} and~\ref{2xThin2B}.

\subsection{Assumptions}\label{assumptions2B}

\subsection{Subgroups of the type $2^2.{}^2E_6(2).2$ in $2.B$}%
\label{V42E622in2B}

We assume that the sporadic simple Monster group $M$ contains an involution
$z$ whose centralizer in $M$ has the structure $2.B$,
and that $M$ contains an elementary abelian subgroup $E$ of order four
that contains $z$ and whose normalizer in $M$ has the structure
$2^2.{}^2E_6(2).S_3$.
The three involutions in $E$ are conjugate in $N_M(E)$,
thus $N_{2.B}(E) = 2.B \cap N$ has the structure $2^2.{}^2E_6(2).2$.
The factor group $N_{2.B}(E) / \langle z \rangle$ is an
involution centralizer in the Baby Monster group $B$,
thus it is isomorphic with the group of the structure
$2.{}^2E_6(2).2$ whose character table is constructed in
Section~\ref{section22e622}.

We denote $N_{2.B}(E)$ by $2.S$
and $N_{2.B}(E) / \langle z \rangle$ by $S$.

We will need the character table of $2.S$.
It can be constructed with character theoretic methods, as follows.

\begin{itemize}
\item
    The character table of the index two subgroup $2.S'$ of $2.S$
    can be constructed from the character table of $2.{}^2E_6(2)$
    and the action of the outer automorphism of order three on the
    classes of ${}^2E_6(2)$,
    using that this automorphism lifts to an automorphism of $2.S'$
    that permutes the three central involutions transitively;
    the general method for this kind of character table construction
    is described in the section
    ``Character Tables of Groups of the Structure $2^2.G$'' of~\cite{Auto},
    and the case of $2.S'$ is treated in the section
    ``Examples for the Type $2^2.G$''.
\item
    Now the character table of $2.S$ can be constructed
    from the character tables of its subgroup $2.S'$ and its factor group
    of the structure $2.{}^2E_6(2).2$,
    using the general method described in \cite{Bre11}.
    The case of $2.S$ is treated in the section
    ``Examples for the Type $M.G.A$'' of \cite{Auto}.
\end{itemize}

The character tables of $S$ and $2.S$ is available in
{\GAP}'s character table library,
via the names \verb|"2.2E6(2).2"| and \verb|"2^2.2E6(2).2"|.

\begin{verbatim}
    gap> t:= CharacterTable( "B" );;
    gap> s:= CharacterTable( "2.2E6(2).2" );;
    gap> 2s:= CharacterTable( "2^2.2E6(2).2" );;
    gap> 2sfuss:= GetFusionMap( 2s, s );;
\end{verbatim}

First we determine the class fusion from $S$ to $B$.
Up to the action of those table automorphisms of the character table of $S$
that lift to the character table of $2.S$,
we get two candidates.

\begin{verbatim}
    gap> fus:= PossibleClassFusions( s, t );;
    gap> Length( fus );
    16
    gap> n:= NrConjugacyClasses( 2s );;
    gap> indperm:= pi -> PermList( CompositionMaps( 2sfuss,
    >        CompositionMaps( ListPerm( pi, n ), InverseMap( 2sfuss ) ) ) );;
    gap> ind:= Group( List( GeneratorsOfGroup( AutomorphismsOfTable( 2s ) ),
    >                       indperm ) );;
    gap> Size( ind );
    16
    gap> rep:= RepresentativesFusions( ind, fus, Group( () ) );;
    gap> Length( rep );
    2
\end{verbatim}


Let $\pi$ denote the epimorphism from $2.B$ to $B$.
We know that the preimage of the centre of $S$ under $\pi$
is elementary abelian.
Since the central involution of $S$ lies in the class \verb|2A| of $B$,
this means that the preimages of this class in $2.B$ are involutions.
Only the first fusion candidate is compatible with this condition.

\begin{verbatim}
    gap> List( rep, map -> map{ ClassPositionsOfCentre( s ) } );
    [ [ 1, 2 ], [ 1, 2 ] ]
    gap> List( rep, map -> Positions( map, 2 ) );
    [ [ 2, 4, 175 ], [ 2, 4, 176 ] ]
    gap> pos:= List( [ 175, 176 ], i -> Positions( 2sfuss, i ) );
    [ [ 251 ], [ 252 ] ]
    gap> OrdersClassRepresentatives( 2s ){ [ 251, 252 ] };
    [ 2, 4 ]
    gap> sfust:= rep[1];;
\end{verbatim}

Next we use the subgroups $S$ and $2.S$ to determine the class fusion
from $2.B$ to $B$, that is, to decide for which classes of $B$ the
preimage under $\pi$ forms one class of $2.B$ or splits into two
classes of $2.B$.
The tools for that are introduced in Section~\ref{theclasses}.

\begin{verbatim}
    gap> orders:= OrdersClassRepresentatives( t );;
    gap> mustsplit:= PositionsProperty( orders, IsOddInt );
    [ 1, 6, 7, 18, 19, 31, 46, 47, 54, 75, 81, 82, 91, 98, 109, 112, 113, 128, 
      131, 145, 146, 151, 155, 160, 172, 173, 177 ]
    gap> selfCentralizingClassesSplit( t, mustsplit );
    #I  class 158 splits (self-centralizing)
    #I  class 159 splits (self-centralizing)
    #I  class 165 splits (self-centralizing)
    #I  class 169 splits (self-centralizing)
    #I  class 170 splits (self-centralizing)
    #I  class 171 splits (self-centralizing)
    #I  class 176 splits (self-centralizing)
    #I  class 182 splits (self-centralizing)
    #I  class 183 splits (self-centralizing)
    #I  class 184 splits (self-centralizing)
    gap> splittingClassesWithOddCentralizerIndexSplit( s, t, sfust,
    >        2sfuss, mustsplit );
    #I  class 110 splits (odd centralizer index)
    #I  class 68 splits (odd centralizer index)
    #I  class 73 splits (odd centralizer index)
    #I  class 79 splits (odd centralizer index)
    #I  class 94 splits (odd centralizer index)
    #I  class 136 splits (odd centralizer index)
    #I  class 140 splits (odd centralizer index)
    gap> mustnotsplit:= [];;
    gap> notSplittingClassesOfSubgroupDoNotSplit( 2sfuss, sfust, mustnotsplit );
    #I  class 2 does not split (as in subgroup)
    #I  class 4 does not split (as in subgroup)
    #I  class 5 does not split (as in subgroup)
    #I  class 8 does not split (as in subgroup)
    #I  class 9 does not split (as in subgroup)
    #I  class 10 does not split (as in subgroup)
    #I  class 11 does not split (as in subgroup)
    #I  class 13 does not split (as in subgroup)
    #I  class 14 does not split (as in subgroup)
    #I  class 15 does not split (as in subgroup)
    #I  class 17 does not split (as in subgroup)
    #I  class 20 does not split (as in subgroup)
    #I  class 21 does not split (as in subgroup)
    #I  class 23 does not split (as in subgroup)
    #I  class 24 does not split (as in subgroup)
    #I  class 25 does not split (as in subgroup)
    #I  class 27 does not split (as in subgroup)
    #I  class 29 does not split (as in subgroup)
    #I  class 30 does not split (as in subgroup)
    #I  class 32 does not split (as in subgroup)
    #I  class 33 does not split (as in subgroup)
    #I  class 34 does not split (as in subgroup)
    #I  class 35 does not split (as in subgroup)
    #I  class 36 does not split (as in subgroup)
    #I  class 37 does not split (as in subgroup)
    #I  class 38 does not split (as in subgroup)
    #I  class 39 does not split (as in subgroup)
    #I  class 40 does not split (as in subgroup)
    #I  class 41 does not split (as in subgroup)
    #I  class 44 does not split (as in subgroup)
    #I  class 48 does not split (as in subgroup)
    #I  class 50 does not split (as in subgroup)
    #I  class 52 does not split (as in subgroup)
    #I  class 55 does not split (as in subgroup)
    #I  class 56 does not split (as in subgroup)
    #I  class 57 does not split (as in subgroup)
    #I  class 58 does not split (as in subgroup)
    #I  class 59 does not split (as in subgroup)
    #I  class 61 does not split (as in subgroup)
    #I  class 62 does not split (as in subgroup)
    #I  class 63 does not split (as in subgroup)
    #I  class 65 does not split (as in subgroup)
    #I  class 66 does not split (as in subgroup)
    #I  class 67 does not split (as in subgroup)
    #I  class 69 does not split (as in subgroup)
    #I  class 70 does not split (as in subgroup)
    #I  class 71 does not split (as in subgroup)
    #I  class 72 does not split (as in subgroup)
    #I  class 76 does not split (as in subgroup)
    #I  class 77 does not split (as in subgroup)
    #I  class 78 does not split (as in subgroup)
    #I  class 80 does not split (as in subgroup)
    #I  class 83 does not split (as in subgroup)
    #I  class 84 does not split (as in subgroup)
    #I  class 85 does not split (as in subgroup)
    #I  class 86 does not split (as in subgroup)
    #I  class 87 does not split (as in subgroup)
    #I  class 88 does not split (as in subgroup)
    #I  class 92 does not split (as in subgroup)
    #I  class 93 does not split (as in subgroup)
    #I  class 95 does not split (as in subgroup)
    #I  class 97 does not split (as in subgroup)
    #I  class 99 does not split (as in subgroup)
    #I  class 101 does not split (as in subgroup)
    #I  class 102 does not split (as in subgroup)
    #I  class 103 does not split (as in subgroup)
    #I  class 114 does not split (as in subgroup)
    #I  class 115 does not split (as in subgroup)
    #I  class 116 does not split (as in subgroup)
    #I  class 117 does not split (as in subgroup)
    #I  class 118 does not split (as in subgroup)
    #I  class 119 does not split (as in subgroup)
    #I  class 120 does not split (as in subgroup)
    #I  class 122 does not split (as in subgroup)
    #I  class 123 does not split (as in subgroup)
    #I  class 124 does not split (as in subgroup)
    #I  class 126 does not split (as in subgroup)
    #I  class 129 does not split (as in subgroup)
    #I  class 130 does not split (as in subgroup)
    #I  class 132 does not split (as in subgroup)
    #I  class 133 does not split (as in subgroup)
    #I  class 134 does not split (as in subgroup)
    #I  class 135 does not split (as in subgroup)
    #I  class 137 does not split (as in subgroup)
    #I  class 138 does not split (as in subgroup)
    #I  class 139 does not split (as in subgroup)
    #I  class 141 does not split (as in subgroup)
    #I  class 149 does not split (as in subgroup)
    #I  class 150 does not split (as in subgroup)
    #I  class 152 does not split (as in subgroup)
    #I  class 153 does not split (as in subgroup)
    #I  class 154 does not split (as in subgroup)
    #I  class 156 does not split (as in subgroup)
    #I  class 157 does not split (as in subgroup)
    #I  class 161 does not split (as in subgroup)
    #I  class 163 does not split (as in subgroup)
    #I  class 166 does not split (as in subgroup)
    #I  class 167 does not split (as in subgroup)
    #I  class 168 does not split (as in subgroup)
    #I  class 175 does not split (as in subgroup)
    #I  class 178 does not split (as in subgroup)
    #I  class 179 does not split (as in subgroup)
    #I  class 180 does not split (as in subgroup)
    #I  class 181 does not split (as in subgroup)
    gap> proj:= Filtered( Irr( 2s ), x -> x[1] <> x[2] );;
    gap> projmap:= ProjectionMap( 2sfuss );;
    gap> proj:= List( proj, x -> x{ projmap } );;
    gap> computeContributions( s, t, sfust, proj, 10^6,
    >        mustsplit, mustnotsplit );
    #I  class 3 splits (contribution criterion)
    #I  class 12 splits (contribution criterion)
    #I  class 42 splits (contribution criterion)
    #I  class 22 splits (3rd root of 3)
    #I  class 26 splits (3rd root of 3)
    #I  class 60 splits (3rd root of 12)
    #I  class 64 splits (3rd root of 12)
    #I  class 125 splits (3rd root of 42)
    #I  class 49 splits (5th root of 3)
    #I  class 51 splits (5th root of 3)
    #I  class 105 splits (5th root of 12)
    #I  class 143 splits (5th root of 26)
    #I  class 144 splits (5th root of 26)
    #I  class 111 splits (11th root of 3)
    #I  class 104 does not split (contribution criterion)
    #I  class 28 splits (contribution criterion)
    #I  class 142 splits (5th root of 28)
\end{verbatim}


\subsection{Element orders in $2.B$}\label{elementorders2B}

All classes of involutions in $B$ contain elements of $S$.
We see that exactly the preimages under $\pi$ of the elements
in the third involution class \verb|2C| have order four,
the other involutions of $B$ lift to involutions in $2.B$.

\begin{verbatim}
    gap> invol:= Positions( OrdersClassRepresentatives( t ), 2 );
    [ 2, 3, 4, 5 ]
    gap> IsSubset( sfust, invol );
    true
    gap> invols:= List( invol, i -> Positions( sfust, i ) );
    [ [ 2, 4, 175 ], [ 3, 5 ], [ 176, 177 ], [ 6, 7, 178 ] ]
    gap> preim2s:= List( invols,
    >        l -> PositionsProperty( 2sfuss, x -> x in l ) );
    [ [ 3, 6, 251 ], [ 4, 5, 7, 8 ], [ 252, 253 ], [ 9, 10, 254 ] ]
    gap> List( preim2s, l -> OrdersClassRepresentatives( 2s ){ l } );
    [ [ 2, 2, 2 ], [ 2, 2, 2, 2 ], [ 4, 4 ], [ 2, 2, 2 ] ]
    gap> invmustlift:= [ 4 ];;
    gap> invmaylift:= [];;
\end{verbatim}

\subsection{The class fusion from $2^2.{}^2E_6(2).2$}

For two classes of $B$ that contain elements of $S$,
we have not yet decided whether the preimages under $\pi$ consist
of one or two classes of $2.B$.
We try to use norms of induced characters for a decision.

\begin{verbatim}
    gap> open:= Difference( sfust, Union( mustsplit, mustnotsplit ) );
    [ 89, 90 ]
    gap> defined:= [];;
    gap> ker:= ClassPositionsOfKernel( GetFusionMap( 2s, s ) );
    [ 1, 2 ]
    gap> testcharss:= Filtered( Irr( 2s ),
    >        chi -> not IsSubset( ClassPositionsOfKernel( chi ), ker ) );;
    gap> good:= [];;
    gap> for choice in Combinations( open ) do
    >      2t:= tableHead( t, Union( mustsplit, choice ),
    >                      invmustlift, invmaylift ); 
    >      fus:= runOneTest( s, 2s, t, 2t, sfust, testcharss, defined );
    >      if fus <> fail then
    >        Add( good, rec( choice:= choice, table:= 2t, map:= fus ) );
    >      fi;
    >    od;
    gap> List( good, x -> x.choice );
    [ [ 89 ], [ 90 ] ]
\end{verbatim}

We see that exactly one of the two classes must split.
The first possibility can be excluded using other scalar products
between induced characters;
for that, we check all four possible fusions for this case.

\begin{verbatim}
    gap> List( good, x -> Indeterminateness( x.map ) );
    [ 4, 4 ]
    gap> goodfus:= [];;
    gap> for map in ContainedMaps( good[1].map ) do
    >      2t:= good[1].table;
    >      ind:= InducedClassFunctionsByFusionMap( 2s, 2t, testcharss, map );
    >      if ForAll( Flat( MatScalarProducts( 2t, ind, ind ) ), IsInt ) then
    >        Add( goodfus, map );
    >      fi;
    >    od;
    gap> Length( goodfus );
    0
    gap> 2t:= good[2].table;;
    gap> 2tfust:= GetFusionMap( 2t, t );;
    gap> AddSet( mustnotsplit, 89 );
    gap> AddSet( mustsplit, 90 );
\end{verbatim}

Two of the remaining four fusion candidates get excluded
by the same criterion.

\begin{verbatim}
    gap> inds:= [];;
    gap> for map in ContainedMaps( good[2].map ) do
    >      ind:= InducedClassFunctionsByFusionMap( 2s, 2t, testcharss, map );
    >      if ForAll( Flat( MatScalarProducts( 2t, ind, ind ) ), IsInt ) then
    >        Add( inds, rec( 2tfust:= 2tfust, characters:= ind, map:= map ) );
    >      fi;
    >    od;
    gap> Length( inds );
    2
\end{verbatim}

The second candidate turns out to admit no $2$-nd power map.
This determines the class fusion from $2.S$ and the characters induced
from irreducible characters of $2.S$.

\begin{verbatim}
    gap> p:= 2;;
    gap> 2tfust:= GetFusionMap( 2t, t );;
    gap> for testinds in inds do
    >      pow:= InitPowerMap( 2t, p );
    >      Congruences( 2t, testinds.characters, pow, p, false );
    >      TransferDiagram( pow, 2tfust, PowerMap( t, p ) );
    >      TransferDiagram( PowerMap( 2s, p ), testinds.map, pow );
    >      factirr:= List( Irr( t ), x -> x{ 2tfust } );
    >      testinds.pow2:= PowerMapsAllowedBySymmetrizations( 2t, factirr,
    >                          testinds.characters, pow, p, parametersFABR );
    >    od;
    gap> List( inds, x -> Length( x.pow2 ) );
    [ 1, 0 ]
    gap> inds:= inds[1];;
\end{verbatim}

Currently $17$ classes of $B$ are left for which we do not know yet
whether they split under $\pi^{-1}$ or not.

\begin{verbatim}
    gap> Length( Difference( [ 1 .. NrConjugacyClasses( t ) ],
    >            Union( mustsplit, mustnotsplit ) ) );
    17
\end{verbatim}

We consider another subgroup of $B$ for getting more information
about these classes.

\subsection{Subgroups of the type $2 \times Th$ in $2.B$}%
\label{2xThin2B}

The sporadic simple Monster group contains a subgroup of the structure
$S_3 \times Th$, where the involutions in the direct factor $S_3$
are conjugate to $z$, see \cite[p.~207]{Wilson87}.
Thus $2.B$ contains a subgroup of the structure $2 \times Th$.


We apply our criteria to this subgroup.

\begin{verbatim}
    gap> th:= CharacterTable( "Th" );;
    gap> poss:= PossibleClassFusions( th, t );;
    gap> Length( poss );
    2
    gap> rep:= RepresentativesFusions( th, poss, Group( () ) );;
    gap> Length( rep );
    1
    gap> thfust:= rep[1];;
    gap> 2th:= CharacterTable( "Cyclic", 2 ) * th;;
    gap> 2thfusth:= GetFusionMap( 2th, th );;
    gap> splittingClassesWithOddCentralizerIndexSplit( th, t, thfust,
    >        2thfusth, mustsplit );
    #I  class 96 splits (odd centralizer index)
\end{verbatim}

We have to adjust the class fusion and the induced characters from $2.S$
to the splitting of class $96$.
Since such cases will occur several times,
we write a small auxiliary function.

\begin{verbatim}
    gap> splitFusionAndCharacters:= function( r, t, tosplit_in_t )
    >      local 2tfust, inv, tosplit_in_2t, result, shift, j, i, spl;
    > 
    >      2tfust:= r.2tfust;
    >      inv:= InverseMap( 2tfust );
    >      tosplit_in_2t:= inv{ tosplit_in_t };
    >      if ForAny( inv{ tosplit_in_t }, IsList ) then
    >        Error( "the classes in ",
    >               Filtered( tosplit_in_t, i -> IsList( inv[i] ) ),
    >               " were already split" );
    >      elif ForAny( r.characters,
    >                   x -> not IsZero( x{ tosplit_in_2t } ) ) then
    >        Error( "all characters must vanish on the classes to be split" );
    >      fi;
    > 
    >      # Adjust the characters of '2t'.
    >      spl:= Concatenation( [ 1 .. Length( r.characters[1] ) ],
    >                           tosplit_in_2t );
    >      Sort( spl );
    >      result:= rec( characters:= List( r.characters, chi -> chi{ spl } ),
    >                    2tfust:= 2tfust{ spl } );
    > 
    >      if IsBound( r.map ) then
    >        if Intersection( tosplit_in_2t, r.map ) <> [] then
    >          Error( "the classes to be split must not occur in the subgroup" );
    >        fi;
    > 
    >        # Adjust the fusion from a subgroup to '2t'.
    >        Add( tosplit_in_2t, Length( 2tfust ) + 1 );
    >        shift:= [];
    >        for j in [ 1 .. tosplit_in_2t[1] - 1 ] do
    >          shift[j]:= 0;
    >        od;
    >        for i in [ 1 .. Length( tosplit_in_2t ) - 1 ] do
    >          for j in [ tosplit_in_2t[i] .. tosplit_in_2t[ i+1 ] - 1 ] do
    >            shift[j]:= i;
    >          od;
    >        od;
    >        result.map:= r.map + shift{ r.map };
    >      fi;
    > 
    >      return result;
    > end;;
    gap> inds:= splitFusionAndCharacters( inds, t, [ 96 ] );;
\end{verbatim}

Four classes of $B$ for which we do not yet know whether they split
contain elements of $Th$.

\begin{verbatim}
    gap> open:= Difference( thfust, Union( mustsplit, mustnotsplit ) );
    [ 45, 53, 108, 127 ]
    gap> defined:= Set( sfust );;
    gap> ker:= ClassPositionsOfKernel( GetFusionMap( 2th, th ) );
    [ 1, 49 ]
    gap> testcharsth:= Filtered( Irr( 2th ),
    >               x -> not IsSubset( ClassPositionsOfKernel( x ), ker ) );;
    gap> good:= [];;
    gap> for choice in Combinations( open ) do
    >      2t:= tableHead( t, Union( mustsplit, choice ),
    >                      invmustlift, invmaylift );
    >      fus:= runOneTest( th, 2th, t, 2t, thfust, testcharsth, defined );
    >      if fus <> fail then
    >        Add( good, rec( choice:= choice, map:= fus, table:= 2t ) );
    >      fi;
    >    od;
    gap> List( good, x -> x.choice );
    [ [ 45, 53, 127 ], [ 53 ], [ 53, 127 ] ]
\end{verbatim}

We see that the class $53$ splits, and the class $108$ does not split,
but the current information does not allow us to decide which of the
three possible solutions is the right one.

\begin{verbatim}
    gap> UniteSet( mustsplit, [ 53 ] );
    gap> UniteSet( mustnotsplit, [ 108 ] );
\end{verbatim}


Thus we consider all three cases,
adjust the list of characters induced from $2.S$ accordingly,
and use them to determine the class fusion from $2 \times Th$.
Using also the compatibility with the $2$-nd power map,
each possible splitting of classes yields a unique class fusion
from $2 \times Th$ and thus a unique list of induced characters
from $2.S$ and from $2 \times Th$.

\begin{verbatim}
    gap> good2:= [];;
    gap> for r in good do
    >      splitinds:= splitFusionAndCharacters( inds, t, r.choice );
    >      2t:= r.table;;
    >      fus:= StructuralCopy( r.map );;
    >      2tfust:= GetFusionMap( 2t, t );;
    >      factirr:= List( Irr( t ), x -> x{ 2tfust } );;
    >      possfus:= FusionsAllowedByRestrictions( 2th, 2t, testcharsth,
    >                   splitinds.characters, fus, parametersFABR );;
    >      for paramap in possfus do
    >        for map in ContainedMaps( paramap ) do
    >          indth:= Set( InducedClassFunctionsByFusionMap( 2th, 2t,
    >                           testcharsth, map ) );
    >          if ForAll( Flat( MatScalarProducts( 2t, indth, indth ) ),
    >                     IsInt ) and
    >             ForAll( Flat( MatScalarProducts( 2t, indth,
    >                               splitinds.characters ) ), IsInt ) then
    >            # Use the 2-nd power map.
    >            ind:= Concatenation( splitinds.characters, indth );
    >            pow:= InitPowerMap( 2t, p );
    >            if Congruences( 2t, ind, pow, p, false ) = true and
    >               TransferDiagram( pow, 2tfust, PowerMap( t, p ) ) <> fail and
    >               TransferDiagram( PowerMap( 2th, p ), map, pow ) <> fail and
    >               TransferDiagram( PowerMap( 2s, p ), splitinds.map,
    >                                pow ) <> fail then
    >              poss:= PowerMapsAllowedBySymmetrizations( 2t, factirr, ind,
    >                         pow, p, parametersFABR );
    >              if Length( poss ) <> 0 then
    >                r.pow2:= poss;
    >                Add( good2, rec( table:= 2t,
    >                                 choice:= r.choice,
    >                                 2thfus2t:= map,
    >                                 ind:= ind,
    >                                 2sfus2t:= splitinds.map ) );
    >              fi;
    >            fi;
    >          fi;
    >        od;
    >      od;
    >    od;
    gap> List( good2, x -> x.choice );
    [ [ 45, 53, 127 ], [ 53 ], [ 53, 127 ] ]
\end{verbatim}

\subsection{Additional characters of $2.B$}

First we consider the case where only class $53$ splits.
The cases where also the classes $45$ or $127$ split
will be dealt with later.

\begin{verbatim}
    gap> 2t:= good2[2].table;;
    gap> 2tfust:= GetFusionMap( 2t, t );;
    gap> 2thfus2t:= good2[2].2thfus2t;;
    gap> 2sfus2t:= good2[2].2sfus2t;;
    gap> ind:= good2[2].ind;;
    gap> factirr:= List( Irr( t ), x -> x{ 2tfust } );;
\end{verbatim}

Several classes of $B$ that are known to split into two classes of $2.B$
are not contained in the conjugates of the subgroups $2.S$ and $2 \times Th$.

\begin{verbatim}
    gap> nothit_in_t:= Difference( mustsplit, Union( thfust, sfust ) );;
    gap> nothit_in_2t:= PositionsProperty( 2tfust, i -> i in nothit_in_t );
    [ 66, 67, 136, 137, 147, 148, 149, 150, 162, 163, 166, 167, 186, 187, 188, 
      189, 217, 218, 222, 223, 224, 225, 226, 227, 228, 229, 230, 231, 234, 235, 
      236, 237, 242, 243 ]
    gap> orders_2t:= OrdersClassRepresentatives( 2t );;
    gap> orders_2t{ nothit_in_2t };
    [ 10, 10, 20, 20, 23, 46, 23, 46, 24, 24, 25, 50, 30, 30, 30, 30, 40, 40, 44, 
      44, 46, 46, 46, 46, 47, 94, 47, 94, 104, 104, 55, 110, 60, 60 ]
\end{verbatim}

In order to get characters of $2.B$ which distinguish the preimages of these
classes, we want to induce from the cyclic subgroups.
For that, we need the values of the power maps on these classes.


\begin{verbatim}
    gap> powermaps:= ComputedPowerMaps( 2t );;
    gap> primes:= Filtered( [ 1 .. Maximum( orders_2t ) ], IsPrimeInt );;
    gap> for p in primes do
    >      pow:= InitPowerMap( 2t, p );
    >      if TransferDiagram( pow, 2tfust, PowerMap( t, p ) ) = fail or
    >         TransferDiagram( PowerMap( 2th, p ), 2thfus2t, pow ) = fail or
    >         TransferDiagram( PowerMap( 2s, p ), 2sfus2t, pow ) = fail or
    >         ConsiderSmallerPowerMaps( 2t, pow, p, false ) <> true or
    >         Congruences( 2t, ind, pow, p, false ) <> true then
    >        Error( "contradiction" );
    >      fi;
    >      poss:= PowerMapsAllowedBySymmetrizations( 2t, ind, ind,
    >                 pow, p, parametersFABR );
    >      if Length( poss ) = 1 then
    >        powermaps[p]:= poss[1];
    >      else
    >        powermaps[p]:= pow;
    >      fi;
    >    od;
\end{verbatim}


Now the power maps for several of the interesting classes are known.

\begin{verbatim}
    gap> known:= Filtered( nothit_in_2t,
    >                i -> ForAll( powermaps, map -> IsInt( map[i] ) ) );
    [ 147, 148, 149, 150, 166, 167, 228, 229, 230, 231, 236, 237 ]
\end{verbatim}

Note that the ambiguities in the current power map information
concern the question whether the $p$-th power of an element is conjugate
to some element $g$ or to $g z$,
since the power maps of the factor group $B$ must be respected.
In several cases, we can improve the power map information using the
following criterion.

Suppose that $g$ and $g z$ are not conjugate
and that $p$ is an odd prime such that $g^p$ is conjugate
to either $g$ or $g z$.
Suppose that $k$ is an odd number (dividing the order of $g$)
such that $g^k$ and $g^k z$ are not conjugate.
Then $g^p$ is conjugate to $g$ if and only if $(g^k)^p$ is conjugate to $g^k$.


\begin{verbatim}
    gap> List( Difference( nothit_in_2t, known ),
    >          i -> Number( powermaps, map -> IsList( map[i] ) ) );
    [ 22, 22, 25, 25, 26, 26, 26, 26, 26, 26, 27, 27, 26, 26, 27, 27, 27, 27, 27, 
      27, 27, 27 ]
    gap> oddprimes:= Difference( primes, [ 2 ] );;
    gap> inv:= InverseMap( 2tfust );;
    gap> for i in [ 1 .. Length( inv ) ] do
    >      if IsList( inv[i] ) then
    >        # the classes of g and g*z
    >        pair:= inv[i];
    >        for p in oddprimes do
    >          if powermaps[p]{ pair } = [ pair, pair ] then
    >            # the p-th powers of g and g*z are conj. to g or g*z
    >            for k in Filtered( oddprimes,
    >                               x -> orders_2t[ pair[1] ] mod x = 0 ) do
    >              img:= powermaps[k][ pair[1] ];
    >              if IsList( img ) then
    >                if powermaps[p]{ img } = img then
    >                  # the p-th power of g^k is conj. to g^k 
    >                  powermaps[p]{ pair }:= pair;
    >                elif powermaps[p]{ img } = Reversed( img ) then
    >                  # the p-th power of g^k is conj. to g^k*z
    >                  powermaps[p]{ pair }:= pair{ [ 2, 1 ] };
    >                fi;
    >              fi;
    >            od;
    >          fi;
    >        od;
    >      fi;
    >    od;
    gap> List( Difference( nothit_in_2t, known ),
    >          i -> Number( powermaps, map -> IsList( map[i] ) ) );
    [ 1, 1, 1, 1, 1, 1, 17, 17, 17, 17, 27, 27, 26, 26, 18, 18, 18, 18, 27, 27, 
      2, 2 ]
\end{verbatim}

We may choose the $5$-th powers of the classes $66$, $67$, $136$, $137$,
and the $3$-rd powers of the classes $162$, $163$.

\begin{verbatim}
    gap> pos:= [ 66, 67, 136, 137 ];;
    gap> powermaps[5]{ pos };
    [ [ 4, 5 ], [ 4, 5 ], [ 16, 17 ], [ 16, 17 ] ]
    gap> powermaps[5]{ pos }:= [ 4, 5, 16, 17 ];;
    gap> pos:= [ 162, 163 ];;
    gap> powermaps[3]{ pos };
    [ [ 53, 54 ], [ 53, 54 ] ]
    gap> powermaps[3]{ pos }:= [ 53, 54 ];;
\end{verbatim}

Concerning the classes $242$ and $243$ (of element order $60$),
the $3$-rd and $5$-th powers are not yet determined.
We may choose the $3$-rd powers, this determines the $5$-th powers.

\begin{verbatim}
    gap> pos:= [ 242, 243 ];;
    gap> powermaps[3]{ pos };
    [ [ 136, 137 ], [ 136, 137 ] ]
    gap> powermaps[5]{ pos };
    [ [ 78, 79 ], [ 78, 79 ] ]
    gap> powermaps[3]{ [ 78, 79 ] };
    [ 16, 17 ]
    gap> powermaps[3]{ pos }:= [ 136, 137 ];;
    gap> powermaps[5]{ pos }:= [ 78, 79 ];;
\end{verbatim}

Let us look at the four relevant classes of element order $46$.
They are the preimages of two Galois conjugate classes of $B$.
The current power maps would admit the possibilities that all four classes
are Galois conjugate or that they are two pairs of generators of
nonconjugate cyclic subgroups.
The former possibility is excluded by the fact that
the number field generated by $46$-th roots of unity does not
contain a field of degree four.
Thus the character values on the four classes lie in the unique quadratic
subfield, which is the field generated by $\sqrt{-23}$.
We can set the correct power map values with the function
\verb|setGaloisInfo| that has been introduced in Section~\ref{sectpowermapsTh}.

\begin{verbatim}
    gap> pos:= Intersection( Difference( nothit_in_2t, known ),
    >                        Positions( orders_2t, 46 ) );
    [ 224, 225, 226, 227 ]
    gap> powermaps[5]{ pos };
    [ [ 226, 227 ], [ 226, 227 ], [ 224, 225 ], [ 224, 225 ] ]
    gap> setGaloisInfo( powermaps, [ 224, 226 ], orders_2t, primes, Sqrt(-23) );
    gap> setGaloisInfo( powermaps, [ 225, 227 ], orders_2t, primes, Sqrt(-23) );
\end{verbatim}

Now only the $23$-rd power map is not uniquely determined on the classes
in question.
Since we have not yet distinguished the two preimages of the corresponding
classes of the factor group $B$ (all induced characters computed up to now
are zero on these classes), we may choose them.

\begin{verbatim}
    gap> powermaps[23]{ pos };                                            
    [ [ 4, 5 ], [ 4, 5 ], [ 4, 5 ], [ 4, 5 ] ]
    gap> powermaps[23]{ pos }:= [ 4, 5, 4, 5 ];;
    gap> ForAll( List( powermaps, x -> x{ pos } ), IsPositionsList );
    true
\end{verbatim}

Let us look at the four relevant classes of element order $30$.
They are the preimages of two Galois conjugate classes of $B$,
with character values in the field generated by $\sqrt{-15}$.

\begin{verbatim}
    gap> pos:= Intersection( nothit_in_2t, Positions( orders_2t, 30 ) );
    [ 186, 187, 188, 189 ]
    gap> Field( Flat( List( factirr, x -> x{ pos } ) ) )
    >    = Field( Rationals, [ Sqrt( -15 ) ] );
    true
\end{verbatim}

The $3$-rd powers of the classes $186$, $187$ are the classes $66$ and $67$,
which are known to be rational from the above computations.
This means that the classes $186$ and $187$ are not Galois conjugate.
Thus there are two pairs of Galois conjugate classes,
and we may choose that the class $186$ is Galois conjugate to $188$.
Furthermore, we may choose that the $3$-rd powers of the class $186$ lie
in the class $66$.
This determines the $5$-th powers.

\begin{verbatim}
    gap> List( powermaps, x -> x[66] );
    [ , 25, 66,, 4,, 66,,,, 66,, 66,,,, 66,, 66,,,, 66,,,,,, 66,, 66,,,,,, 66,,,, 
      66,, 66,,,, 66,,,,,, 66,,,,,, 66,, 66,,,,,, 66,,,, 66,, 66,,,,,, 66,,,, 66,,
      ,,,, 66,,,,,,,, 66,,,, 66,, 66,,,, 66,, 66 ]
    gap> setGaloisInfo( powermaps, [ 186, 188 ], orders_2t, primes, Sqrt(-15) );
    gap> setGaloisInfo( powermaps, [ 187, 189 ], orders_2t, primes, Sqrt(-15) );
    gap> powermaps[3]{ pos };
    [ [ 66, 67 ], [ 66, 67 ], [ 66, 67 ], [ 66, 67 ] ]
    gap> powermaps[3]{ pos }:= [ 66, 67, 66, 67 ];;
    gap> powermaps[5]{ pos };
    [ [ 34, 35 ], [ 34, 35 ], [ 34, 35 ], [ 34, 35 ] ]
    gap> powermaps[3]{ [ 34, 35 ] };
    [ 4, 5 ]
    gap> powermaps[5]{ [ 66, 67 ] };
    [ 4, 5 ]
    gap> powermaps[5]{ pos }:= [ 34, 35, 34, 35 ];;
    gap> ForAll( List( powermaps, x -> x{ pos } ), IsPositionsList ); #  true
    true
\end{verbatim}

Let us look at the two classes of element order $44$.
They are the preimages of a unique (and hence rational) class of $B$.
Thus the character values on the classes of element order $44$ lie in
a quadratic subfield of the field of $44$-th roots of unity;
so the possible fields are the rational number field and the fields
generated by $\sqrt{-1}$, $\sqrt{11}$, $\sqrt{-11}$.
For all four cases, we compute the candidates of characters induced from
the cyclic subgroup of order $44$,
and check their norms and scalar products of $2$-nd symmetrizations
with the irreducible characters of $B$.
It turns out that the field of character values on the classes in question
is generated by $\sqrt{-11}$.

\begin{verbatim}
    gap> pos:= Intersection( nothit_in_2t, Positions( orders_2t, 44 ) );
    [ 222, 223 ]
    gap> vals:= List( [ 1, -1, 11, -11 ], Sqrt );;
    gap> good:= [];;
    gap> for val in vals do
    >      setGaloisInfo( powermaps, pos, orders_2t, primes, val );
    >      indcyc:= InducedCyclic( 2t, pos, "all" );
    >      if ForAll( indcyc, x -> IsInt( ScalarProduct( 2t, x, x ) ) ) then
    >        minus:= MinusCharacter( indcyc[1], powermaps[2], 2 );
    >        if ForAll( List( factirr, x -> ScalarProduct( 2t, x, minus ) ), 
    >                   IsInt ) then
    >          Add( good, val );
    >        fi;
    >      fi;
    >    od;
    gap> good = [ Sqrt( -11 ) ];
    true
    gap> setGaloisInfo( powermaps, pos, orders_2t, primes, good[1] );
\end{verbatim}

Let us look at the two classes of element order $104$.
As in the previous case,
they are the preimages of a unique (and hence rational) class of $B$.
Here the candidates for the fields of character values are generated by
$\sqrt{n}$, for $n \in \{ \pm 1, \pm 2, \pm 13, \pm 26 \}$.
We proceed as above, and get the unique solution $\sqrt{-26}$.

\begin{verbatim}
    gap> pos:= Intersection( nothit_in_2t, Positions( orders_2t, 104 ) );
    [ 234, 235 ]
    gap> vals:= List( [ 1, -1, 2, -2, 13, -13, 26, -26 ], Sqrt );;
    gap> good:= [];;
    gap> for val in vals do
    >      setGaloisInfo( powermaps, pos, orders_2t, primes, val );
    >      indcyc:= InducedCyclic( 2t, pos, "all" );
    >      if ForAll( indcyc, x -> IsInt( ScalarProduct( 2t, x, x ) ) ) then
    >        minus:= MinusCharacter( indcyc[1], powermaps[2], 2 );
    >        if ForAll( List( factirr, x -> ScalarProduct( 2t, x, minus ) ), 
    >                   IsInt ) then
    >          Add( good, val );
    >        fi;
    >      fi;
    >    od;
    gap> good = [ Sqrt( -26 ) ];
    true
    gap> setGaloisInfo( powermaps, pos, orders_2t, primes, good[1] );
\end{verbatim}

Now exactly two classes are left for which we want to determine the
power map values; these classes have element order $40$.

\begin{verbatim}
    gap> pos:= Intersection( nothit_in_2t, Positions( orders_2t, 40 ) );
    [ 217, 218 ]
\end{verbatim}

First we observe that the $2$-nd power map for these classes is determined
by the fact that the $2$-nd power map commutes with the $5$-th power map.

\begin{verbatim}
    gap> powermaps[2]{ pos };
    [ [ 136, 137 ], [ 136, 137 ] ]
    gap> TransferDiagram( powermaps[5], powermaps[2], powermaps[5] );
    rec( impbetween := [ 131, 139, 217, 218 ], impinside1 := [  ], 
      impinside2 := [  ] )
    gap> IsPositionsList( powermaps[2] );
    true
\end{verbatim}


As in the above cases, we have to determine the field of character values
on these classes, in order to write down the power maps.
Here the candidates for the fields of character values are generated by
$\sqrt{n}$, for $n \in \{ \pm 1, \pm 2, \pm 5, \pm 10 \}$.
We proceed as above, but unfortunately we get two solutions.

\begin{verbatim}
    gap> vals:= List( [ 1, -1, 2, -2, 5, -5, 10, -10 ], Sqrt );;
    gap> good:= [];;
    gap> for val in vals do
    >      setGaloisInfo( powermaps, pos, orders_2t, primes, val );
    >      indcyc:= InducedCyclic( 2t, pos, "all" );
    >      if ForAll( indcyc, x -> IsInt( ScalarProduct( 2t, x, x ) ) ) then
    >        minus:= MinusCharacter( indcyc[1], powermaps[2], 2 );
    >        if ForAll( List( factirr, x -> ScalarProduct( 2t, x, minus ) ),
    >                   IsInt ) then
    >          Add( good, val );
    >        fi;
    >      fi;
    >    od;
    gap> good = [ Sqrt( 5 ), Sqrt( -5 ) ];
    true
\end{verbatim}

\subsection{The irreducible characters of $2.B$}

Now we can induce characters from the interesting cyclic subgroups,
except that we get two sets of candidates induced from the subgroups
of order $40$, depending on the power maps.

We try the first case, where the field of character values on class $217$
is generated by $\sqrt{5}$.
As in the preceding sections, we apply the LLL algorithm to the set
of induced characters, and compute the possible orthogonal embeddings.
Note that the maximal possible dimension is $63 + 14 = 77$.

\begin{verbatim}
    gap> indcyc:= InducedCyclic( 2t, Difference( nothit_in_2t, pos ), "all" );;
    gap> indcyc:= Reduced( 2t, factirr, indcyc ).remainders;;
    gap> setGaloisInfo( powermaps, pos, orders_2t, primes, Sqrt( 5 ) );
    gap> indcyc40r5:= InducedCyclic( 2t, pos, "all" );;
    gap> indcyc40r5:= Reduced( 2t, factirr, indcyc40r5 ).remainders;;
    gap> testind:= Concatenation( ind, indcyc, indcyc40r5 );;
    gap> lll:= LLL( 2t, testind, 99/100 );;
    gap> Length( lll.norms );
    63
    gap> Length( mustsplit );
    63
    gap> Length( Difference( [ 1 .. NrConjugacyClasses( t ) ],
    >            Union( mustsplit, mustnotsplit ) ) );
    14
    gap> gram:= MatScalarProducts( 2t, lll.remainders, lll.remainders );;
    gap> emb:= OrthogonalEmbeddings( gram, 63 + 14 );;
    gap> List( emb.solutions, Length );
    [ 63, 63, 63, 65, 65, 65 ]
    gap> dec:= List( emb.solutions,
    >              x -> Decreased( 2t, lll.remainders, emb.vectors{ x } ) );;
    gap> Positions( dec, fail );
    [ 1, 2, 3, 6 ]
\end{verbatim}

We got six possible embeddings,
of which two are compatible with the situation that the vectors
are virtual characters;
these two solutions would require (at least) $65$
faithful irreducible characters.

\begin{verbatim}
    gap> dec:= Filtered( dec, x -> x <> fail );;
    gap> List( dec, r -> Length( r.irreducibles ) );
    [ 61, 61 ]
\end{verbatim}

We claim that these solutions cannot lead to the character table of $2.B$.
Suppose that the $61$ vectors are irreducible characters of $2.B$.
Then there are four more irreducible characters $\chi_i$,
$i \in \{ 1, 2, 3, 4 \}$, such that $n = \sum_{i=1}^4 \chi_i(1)^2$
is the difference of $|B|$ and the sum of the squares of the degrees
of the $61$ irreducible characters.

\begin{verbatim}
    gap> degreesum:= List( dec,
    >        r -> Sum( List( r.irreducibles, x -> x[1]^2 ) ) );;
    gap> Set( degreesum );
    [ 4154780380522839827726467072000000 ]
    gap> n:= Size( t ) - degreesum[1];
    1100703586363451113472000000
\end{verbatim}

Reducing the input characters with the $61$ irreducible characters
yields a character $\psi$ of norm $2$ and degree $\sqrt{2n}$.

\begin{verbatim}
    gap> red:= List( dec, r -> Reduced( 2t, r.irreducibles, testind ) );;
    gap> norm2:= List( red, r -> First( r.remainders,
    >                              x -> ScalarProduct( 2t, x, x ) = 2 ) );;
    gap> norm2[1] = norm2[2];
    true
    gap> norm2[1][1] = Sqrt( 2 * n );
    true
\end{verbatim}

Without loss of generality, let $\psi = \chi_1 + \chi_2$.
Then the general inequality $2 (a^2 + b^2) \geq (a+b)^2$
yields $\chi_1(1)^2 + \chi_2(1)^2 \geq \psi(1)^2 / 2 = n$.
This contradicts the conditions $\chi_3(1) \geq 0$ and $\chi_4(1) \geq 0$.

Let us try the other possibility for the power maps.

\begin{verbatim}
    gap> setGaloisInfo( powermaps, pos, orders_2t, primes, Sqrt( -5 ) );
    gap> indcyc40i5:= InducedCyclic( 2t, pos, "all" );;
    gap> indcyc40i5:= Reduced( 2t, factirr, indcyc40i5 ).remainders;;
    gap> testind:= Concatenation( ind, indcyc, indcyc40i5 );;
    gap> lll:= LLL( 2t, testind, 99/100 );;
    gap> gram:= MatScalarProducts( 2t, lll.remainders, lll.remainders );;
    gap> emb:= OrthogonalEmbeddings( gram, 63 + 14 );;
    gap> List( emb.solutions, Length );
    [ 63, 63, 63, 65, 65, 65 ]
    gap> dec:= List( emb.solutions,
    >              x -> Decreased( 2t, lll.remainders, emb.vectors{ x } ) );;
    gap> Positions( dec, fail );
    [ 3, 6 ]
    gap> dec:= Filtered( dec, x -> x <> fail );;
    gap> List( dec, r -> Length( r.irreducibles ) );
    [ 63, 63, 61, 61 ]
\end{verbatim}

Here we get four solutions, two in dimension $63$ and two in dimension $65$.
The latter two solutions are excluded with the same argument as above.

\begin{verbatim}
    gap> dec1:= dec{ [ 1, 2 ] };;
    gap> dec2:= dec{ [ 3, 4 ] };;
    gap> degreesum:= List( dec2,
    >        r -> Sum( List( r.irreducibles, x -> x[1]^2 ) ) );;
    gap> Set( degreesum );
    [ 4154780380522839827726467072000000 ]
    gap> n:= Size( t ) - degreesum[1];
    1100703586363451113472000000
    gap> red:= List( dec2, r -> Reduced( 2t, r.irreducibles, testind ) );;
    gap> norm2:= List( red, r -> First( r.remainders,
    >                              x -> ScalarProduct( 2t, x, x ) = 2 ) );;
    gap> Length( Set( norm2 ) );
    1
    gap> norm2[1][1] = Sqrt( 2 * n );
    true
\end{verbatim}

Thus we are left with the two solutions in dimension $63$,
for which all irreducibles are known.
Each of them leads to a character table that is equivalent to the
character table that is stored for $2.B$ in {\GAP}'s character table library.

\begin{verbatim}
    gap> HasIrr( 2t );
    false
    gap> SetIrr( 2t, Concatenation( factirr, dec1[1].irreducibles ) );
    gap> lib:= CharacterTable( "2.B" );;
    gap> TransformingPermutationsCharacterTables( 2t, lib );
    rec( columns := (4,5)(29,30)(34,35)(63,64)(66,67)(122,123)(125,126)(145,
        146)(186,187)(188,189)(224,225)(226,227), 
      group := <permutation group with 17 generators>, 
      rows := (185,213,243,236,191,205,225,229,247,199,238,235,221,192,214,230,
        201,186,211,227,226,208,198,242,218,193,209,220,196,210,241,203,189,224,
        219,197,204,223,240,190,215,195,233)(187,206,212,234,194,244,217,202,245,
        200,232,222,207,216)(188,239)(231,246,237) )
    gap> ResetFilterObj( 2t, HasIrr );
    gap> SetIrr( 2t, Concatenation( factirr, dec1[2].irreducibles ) );
    gap> TransformingPermutationsCharacterTables( 2t, lib );
    rec( columns := (4,5)(29,30)(34,35)(63,64)(66,67)(122,123)(125,126)(145,
        146)(186,187)(188,189)(222,223)(224,225)(226,227), 
      group := <permutation group with 17 generators>, 
      rows := (185,213,243,190,215,195,233)(186,211,227,226,208,198,242,218,193,
        209,220,196,210,241,203,189,224,219,197,204,223,240,236,235,221,192,214,
        230,201)(187,206,212,234,194,244,217,202,245,200,232,222,207,216)(188,
        239)(191,205,225,229,247,199,238)(231,246,237) )
\end{verbatim}

It remains to show that the two possible class splittings cannot occur
where the classes $127$ or $45$ of $B$ split in $2.B$.

For that, it suffices to split the relevant classes
in the characters that are induced from cyclic subgroups.
First we deal with the case where additionally exactly class $127$ splits.
As above, we have to consider two possible sets of induced characters.



\begin{verbatim}
    gap> good2[3].choice;
    [ 53, 127 ]
    gap> pos:= Positions( 2tfust, 127 );
    [ 165 ]
    gap> 2t:= good2[3].table;;
    gap> 2tfust:= GetFusionMap( 2t, t );;
    gap> 2thfus2t:= good2[3].2thfus2t;;
    gap> 2sfus2t:= good2[3].2sfus2t;;
    gap> ind:= good2[3].ind;;
    gap> factirr:= List( Irr( t ), x -> x{ 2tfust } );;
    gap> UniteSet( mustsplit, good2[3].choice );
    gap> spl:= SortedList( Concatenation( [ 1 .. 247 ], pos ) );;
    gap> testind:= Concatenation( good2[3].ind,
    >        List( Concatenation( indcyc, indcyc40r5 ), x -> x{ spl } ) );;
    gap> lll:= LLL( 2t, testind, 99/100 );;
    gap> Length( lll.norms );
    64
    gap> Length( mustsplit );
    64
    gap> Length( Difference( [ 1 .. NrConjugacyClasses( t ) ],
    >            Union( mustsplit, mustnotsplit ) ) );
    13
    gap> gram:= MatScalarProducts( 2t, lll.remainders, lll.remainders );;
    gap> emb:= OrthogonalEmbeddings( gram, 64+13 );;
    gap> List( emb.solutions, Length );
    [ 64, 64, 64, 65, 65, 65, 66, 66, 66, 67, 67, 67, 67, 67, 67, 67, 67, 67, 69, 
      69, 69 ]
    gap> dec:= List( emb.solutions,
    >            x -> Decreased( 2t, lll.remainders, emb.vectors{ x } ) );;
    gap> Positions( dec, fail );
    [ 1, 2, 3, 6, 7, 8, 9, 16, 17, 18, 21 ]
    gap> dec:= Filtered( dec, x -> x <> fail );;
    gap> List( dec, r -> Length( r.irreducibles ) );
    [ 61, 61, 63, 61, 61, 63, 61, 61, 61, 61 ]
    gap> degreesum:= List( dec,
    >        r -> Sum( List( r.irreducibles, x -> x[1]^2 ) ) );;
    gap> degreesumset:= Set( degreesum );
    [ 4154780380522839827726467072000000, 4154781481226426191177580544000000 ]
    gap> n:= Size( t ) - degreesumset;
    [ 1100703586363451113472000000, 0 ]
    gap> List( degreesumset, x -> Positions( degreesum, x ) );
    [ [ 1, 2, 4, 5, 7, 8, 9, 10 ], [ 3, 6 ] ]
    gap> dec:= dec{ Positions( degreesum, degreesumset[1] ) };;
    gap> red:= List( dec, r -> Reduced( 2t, r.irreducibles, testind ) );;
    gap> norm2:= List( red, r -> First( r.remainders,
    >                              x -> ScalarProduct( 2t, x, x ) = 2 ) );;
    gap> Length( Set( norm2 ) );
    1
    gap> norm2[1][1] = Sqrt( 2 * n[1] );
    true
\end{verbatim}

We see that all possible embeddings would require dimension
at least $65$.
In those cases where $63$ irreducible characters are known,
the sum of squares of their degrees is already too large.
And in those cases where $61$ irreducible characters are known,
there is a character of norm two (as above)
such that no space is left for the remaining irreducibles.

Next we consider the second set of induced characters.

\begin{verbatim}
    gap> testind:= Concatenation( good2[3].ind,
    >        List( Concatenation( indcyc, indcyc40i5 ), x -> x{ spl } ) );;
    gap> lll:= LLL( 2t, testind, 99/100 );;
    gap> gram:= MatScalarProducts( 2t, lll.remainders, lll.remainders );;
    gap> emb:= OrthogonalEmbeddings( gram, 64+13 );;
    gap> List( emb.solutions, Length );
    [ 64, 64, 64, 65, 65, 65, 66, 66, 66, 67, 67, 67, 67, 67, 67, 67, 67, 67, 69, 
      69, 69 ]
    gap> dec:= List( emb.solutions,
    >            x -> Decreased( 2t, lll.remainders, emb.vectors{ x } ) );;
    gap> Positions( dec, fail );
    [ 1, 2, 3, 6, 7, 8, 9, 16, 17, 18, 21 ]
    gap> dec:= Filtered( dec, x -> x <> fail );;
    gap> List( dec, r -> Length( r.irreducibles ) );
    [ 61, 61, 63, 61, 61, 63, 61, 61, 61, 61 ]
    gap> degreesum:= List( dec,
    >        r -> Sum( List( r.irreducibles, x -> x[1]^2 ) ) );;
    gap> degreesumset:= Set( degreesum );
    [ 4154780380522839827726467072000000, 4154781481226426191177580544000000 ]
    gap> n:= Size( t ) - degreesumset;
    [ 1100703586363451113472000000, 0 ]
    gap> List( degreesumset, x -> Positions( degreesum, x ) );
    [ [ 1, 2, 4, 5, 7, 8, 9, 10 ], [ 3, 6 ] ]
    gap> dec:= dec{ Positions( degreesum, degreesumset[1] ) };;
    gap> red:= List( dec, r -> Reduced( 2t, r.irreducibles, testind ) );;
    gap> norm2:= List( red, r -> First( r.remainders,
    >                              x -> ScalarProduct( 2t, x, x ) = 2 ) );;
    gap> Length( Set( norm2 ) );
    1
    gap> norm2[1][1] = Sqrt( 2 * n[1] );
    true
\end{verbatim}

All solutions for this case can be excluded by the same arguments.
Thus \verb|good2[3]| does not lead to a character table.

Finally, we consider the two sets of induced characters in the case
that the two classes $45$ and $127$ of $B$ split in $2.B$.


\begin{verbatim}
    gap> good2[1].choice;
    [ 45, 53, 127 ]
    gap> 2tfust:= GetFusionMap( good2[2].table, t );;
    gap> pos:= Union( Positions( 2tfust, 45 ), Positions( 2tfust, 127 ) );
    [ 57, 165 ]
    gap> 2t:= good2[1].table;;
    gap> 2tfust:= GetFusionMap( 2t, t );;
    gap> 2thfus2t:= good2[1].2thfus2t;;
    gap> 2sfus2t:= good2[1].2sfus2t;;
    gap> ind:= good2[1].ind;;
    gap> factirr:= List( Irr( t ), x -> x{ 2tfust } );;
    gap> UniteSet( mustsplit, good2[1].choice );
    gap> spl:= SortedList( Concatenation( [ 1 .. 247 ], pos ) );;
    gap> testind:= Concatenation( good2[1].ind,
    >        List( Concatenation( indcyc, indcyc40r5 ), x -> x{ spl } ) );;
    gap> lll:= LLL( 2t, testind, 99/100 );;
    gap> Length( lll.norms );
    65
    gap> Length( mustsplit );
    65
    gap> Length( Difference( [ 1 .. NrConjugacyClasses( t ) ],
    >            Union( mustsplit, mustnotsplit ) ) );
    12
    gap> gram:= MatScalarProducts( 2t, lll.remainders, lll.remainders );;
    gap> emb:= OrthogonalEmbeddings( gram, 65+12 );;
    gap> List( emb.solutions, Length );
    [ 66, 66, 66, 66, 66, 66, 66, 66, 66, 67, 67, 67, 67, 67, 67, 67, 67, 67, 67, 
      67, 67, 67, 67, 67, 67, 67, 67, 68, 68, 68, 69, 69, 69, 69, 69, 69, 69, 69, 
      69, 69, 69, 69, 69, 69, 69, 69, 69, 69, 71, 71, 71 ]
    gap> dec:= List( emb.solutions,
    >            x -> Decreased( 2t, lll.remainders, emb.vectors{ x } ) );;
    gap> Positions( dec, fail );
    [ 1, 2, 3, 4, 5, 6, 7, 8, 9, 12, 13, 14, 15, 24, 25, 26, 27, 28, 29, 30, 33, 
      34, 35, 36, 45, 46, 47, 48, 51 ]
    gap> dec:= Filtered( dec, x -> x <> fail );;
    gap> List( dec, r -> Length( r.irreducibles ) );
    [ 63, 63, 63, 61, 61, 61, 63, 61, 61, 61, 63, 63, 63, 61, 61, 61, 63, 61, 61, 
      61, 61, 61 ]
    gap> degreesum:= List( dec,
    >        r -> Sum( List( r.irreducibles, x -> x[1]^2 ) ) );;
    gap> degreesumset:= Set( degreesum );
    [ 4154780380522839827726467072000000, 4154781481226426191177580544000000 ]
    gap> Size( t ) - degreesumset;
    [ 1100703586363451113472000000, 0 ]
    gap> Positions( degreesum, degreesumset[2] );
    [ 1, 2, 3, 7, 11, 12, 13, 17 ]
    gap> dec:= dec{ Positions( degreesum, degreesumset[1] ) };;
    gap> red:= List( dec, r -> Reduced( 2t, r.irreducibles, testind ) );;
    gap> norm2:= List( red, r -> First( r.remainders,
    >                              x -> ScalarProduct( 2t, x, x ) = 2 ) );;
    gap> Length( Set( norm2 ) );
    1
    gap> norm2[1][1];
    0
\end{verbatim}

In the situation of the first set of induced characters,
we get no character table.
The only difference in the arguments is that some solutions
are excluded by the existence of a ``character'' of norm two
and degree zero.

As a last step, we consider the second set of induced characters.
The arguments are exactly the same as above.

\begin{verbatim}
    gap> testind:= Concatenation( good2[1].ind,
    >        List( Concatenation( indcyc, indcyc40i5 ), x -> x{ spl } ) );;
    gap> lll:= LLL( 2t, testind, 99/100 );;
    gap> gram:= MatScalarProducts( 2t, lll.remainders, lll.remainders );;
    gap> emb:= OrthogonalEmbeddings( gram, 65+12 );;
    gap> List( emb.solutions, Length );
    [ 66, 66, 66, 66, 66, 66, 66, 66, 66, 67, 67, 67, 67, 67, 67, 67, 67, 67, 67, 
      67, 67, 67, 67, 67, 67, 67, 67, 68, 68, 68, 69, 69, 69, 69, 69, 69, 69, 69, 
      69, 69, 69, 69, 69, 69, 69, 69, 69, 69, 71, 71, 71 ]
    gap> dec:= List( emb.solutions,
    >            x -> Decreased( 2t, lll.remainders, emb.vectors{ x } ) );;
    gap> Positions( dec, fail );
    [ 1, 2, 3, 4, 5, 6, 7, 8, 9, 12, 13, 14, 15, 24, 25, 26, 27, 28, 29, 30, 33, 
      34, 35, 36, 45, 46, 47, 48, 51 ]
    gap> dec:= Filtered( dec, x -> x <> fail );;
    gap> List( dec, r -> Length( r.irreducibles ) );
    [ 63, 63, 63, 61, 61, 61, 63, 61, 61, 61, 63, 63, 63, 61, 61, 61, 63, 61, 61, 
      61, 61, 61 ]
    gap> degreesum:= List( dec,
    >        r -> Sum( List( r.irreducibles, x -> x[1]^2 ) ) );;
    gap> degreesumset:= Set( degreesum );
    [ 4154780380522839827726467072000000, 4154781481226426191177580544000000 ]
    gap> Size( t ) - degreesumset;
    [ 1100703586363451113472000000, 0 ]
    gap> Positions( degreesum, degreesumset[2] );
    [ 1, 2, 3, 7, 11, 12, 13, 17 ]
    gap> dec:= dec{ Positions( degreesum, degreesumset[1] ) };;
    gap> red:= List( dec, r -> Reduced( 2t, r.irreducibles, testind ) );;
    gap> norm2:= List( red, r -> First( r.remainders,
    >                              x -> ScalarProduct( 2t, x, x ) = 2 ) );;
    gap> Length( Set( norm2 ) );
    1
    gap> norm2[1][1];
    0
\end{verbatim}

Thus also \verb|good2[1]| does not lead to a character table.

\bibliographystyle{amsalpha}

\begin{thebibliography}{CCN{\etalchar{+}}85}

\bibitem[BCP97]{Magma}
W.~Bosma, J.~J. Cannon, and C.~Playoust, \emph{The {Magma} algebra system. {I}.
  {T}he user language}, J. Symbolic Comput. \textbf{24} (1997),
  235{\textendash}265.

\bibitem[BMO]{BMO16}
T.~Breuer, G.~Malle, and E.~A. O'Brien, \emph{Reliability and reproducibility
  of {A}tlas information}, arXiv:1603.08650.

\bibitem[Bre]{Auto}
T.~Breuer, \emph{Using table automorphisms for constructing character tables in
  \textsf{GAP}}, \href
  {http://www.math.rwth-aachen.de/~Thomas.Breuer/ctbllib/doc/ctblcons.pdf}
  {\texttt{http://www.math.rwth-aachen.de/}\discretionary
  {}{}{}\texttt{\texttt{\symbol{126}}Thomas.Breuer/}\discretionary
  {}{}{}\texttt{ctbllib/}\discretionary {}{}{}\texttt{doc/}\discretionary
  {}{}{}\texttt{ctblcons.pdf}}.

\bibitem[Bre11]{Bre11}
\bysame, \emph{Computing character tables of groups of type {$M.G.A$}}, LMS J.
  Comput. Math. \textbf{14} (2011), 173{\textendash}178.

\bibitem[Bre13]{CTblLib}
\bysame, \emph{The \textsf{GAP} {C}haracter {T}able {L}ibrary, {V}ersion
  1.2.2}, \href {http://www.math.rwth-aachen.de/~Thomas.Breuer/ctbllib}
  {\texttt{http://www.math.rwth-aachen.de/}\discretionary
  {}{}{}\texttt{\texttt{\symbol{126}}Thomas.Breuer/}\discretionary
  {}{}{}\texttt{ctbllib}}, Mar 2013, \textsf{GAP} package.

\bibitem[CCN{\etalchar{+}}85]{CCN85}
J.~H. Conway, R.~T. Curtis, S.~P. Norton, R.~A. Parker, and R.~A. Wilson,
  \emph{Atlas of finite groups}, Oxford University Press, Eynsham, 1985,
  Maximal subgroups and ordinary characters for simple groups, With
  computational assistance from J. G. Thackray.

\bibitem[Dem72]{Dem72}
Ulrich Dempwolff, \emph{On extensions of an elementary abelian group of order
  {$2^{5}$} by {${\rm GL}(5,\,2)$}}, Rend. Sem. Mat. Univ. Padova \textbf{48}
  (1972), 359--364 (1973).

\bibitem[GAP16]{GAP483}
The GAP{\nobreakspace}Group, \emph{{GAP{\textendash}Groups, Algorithms, and
  Programming, Version 4.8.3}}, 2016, \href {http://www.gap-system.org}
  {\texttt{http://www.gap-system.org}}.

\bibitem[Jan76]{Jan76}
Zvonimir Janko, \emph{A new finite simple group of order {$86\cdot 775\cdot
  571\cdot 046\cdot 077\cdot 562\cdot 880$} which possesses {$M\sb{24}$} and
  the full covering group of {$M\sb{22}$} as subgroups}, J. Algebra \textbf{42}
  (1976), no.~2, 564--596.

\bibitem[LLL82]{LLL82}
A.~K. Lenstra, H.~W. Lenstra,~Jr., and L.~Lov{\a'a}sz, \emph{Factoring
  polynomials with rational coefficients}, Math. Ann. \textbf{261} (1982),
  no.~4, 515{\textendash}534.

\bibitem[LSS92]{LSS92}
M.~W. Liebeck, J.~Saxl, and G.~M. Seitz, \emph{Subgroups of maximal rank in
  finite exceptional groups of {L}ie type}, Proc. London Math. Soc. (3)
  \textbf{65} (1992), 297{\textendash}325.

\bibitem[Par77]{Par77}
David Parrott, \emph{On {T}hompson's simple group}, J. Algebra \textbf{46}
  (1977), no.~2, 389--404.

\bibitem[Ple95]{Ple90}
W.~Plesken, \emph{Solving {$XX^{\rm tr} = A$} over the integers}, Linear
  Algebra Appl. \textbf{226/228} (1995), 331--344.

\bibitem[Smi76]{Smi76c}
P.~E. Smith, \emph{A simple subgroup of {$M?$} and {$E_8(3)$}}, Bull. London
  Math. Soc. \textbf{8} (1976), no.~2, 161{\textendash}165.

\bibitem[Str76]{Str76b}
G.~Stroth, \emph{A characterization of {F}ischer's sporadic simple group of the
  order {$2\sp{41}\cdot 3\sp{13}\cdot 5\sp{6}\cdot 7\sp{2}\cdot 11\cdot 13\cdot
  17\cdot 19\cdot 23\cdot 31\cdot 47$}}, J. Algebra \textbf{40} (1976), no.~2,
  499--531.

\bibitem[Sul99]{SW99}
R.~A. Suleiman, I. A. I. amd~Wilson, \emph{Construction of exceptional covers
  of generic groups}, Math. Proc. Cambridge Philos. Soc. \textbf{125} (1999),
  31--38.

\bibitem[Wil87]{Wilson87}
R.~A. Wilson, \emph{Some subgroups of the {B}aby {M}onster}, Invent. Math.
  \textbf{89} (1987), no.~1, 197{\textendash}218.

\end{thebibliography}

\newcommand{\etalchar}[1]{$^{#1}$}
\providecommand{\bysame}{\leavevmode\hbox to3em{\hrulefill}\thinspace}
\providecommand{\href}[2]{#2}

\end{document}